\documentclass[opre,sglanonrev]{informs4}

\RequirePackage{tgtermes}
\RequirePackage{newtxtext}
\RequirePackage{newtxmath}
\RequirePackage{bm}
\RequirePackage{endnotes}

\OneAndAHalfSpacedXI

\usepackage{graphicx}
\usepackage{booktabs}

\newcommand{\bblambda}{\boldsymbol{\lambda}}
\newcommand{\bbR}{\mathbb{R}}

\newcommand{\bbD}{\boldsymbol{D}}

\newcommand{\bbE}{\mathbb{E}}

\usepackage[dvipsnames]{xcolor}

\usepackage[colorlinks=true,linkcolor=blue,citecolor=blue,urlcolor=blue]{hyperref}

\usepackage{algorithm}
\usepackage{algpseudocode}
\usepackage{tikz}
\usetikzlibrary{arrows.meta, positioning,calc}
\usepackage{enumitem}
\usepackage{sansmath}

\usepackage{subcaption}

\usepackage{natbib}
 \bibpunct[, ]{(}{)}{,}{a}{}{,}
 \def\bibfont{\small}

\EquationsNumberedThrough

\TheoremsNumberedThrough     

\ECRepeatTheorems

\MANUSCRIPTNO{}

\makeatletter

\def\theARTICLETOP{}

\RRHSecondLine{}
\LRHSecondLine{}

\def\theARTICLEABSTRACT{
  \HOOKb
  \vspace*{18pt}
  \noindent
  \begin{minipage}[t]{\textwidth}\parindent1em
    \ABSfont
    \noindent\theABSTRACT\endgraf
    \vskip5pt
    \theFUNDING
    \theKEYWORDS
    \theSUBJECTCLASS
    \theAREAOFREVIEW
    \theMSCCLASS
    \theORMSCLASS
    \if@BLINDREV\else\theHISTORY\fi
    \noindent\hrulefill
  \end{minipage}
  \vspace*{0pt}
}

\makeatother

\begin{document}

\RUNAUTHOR{Can Er and Mo Liu}

\RUNTITLE{Decision-Focused Bias Correction for Fluid Approximation}

\TITLE{Decision-Focused Bias Correction for Fluid Approximation}

\ARTICLEAUTHORS{

\AUTHOR{Can Er}
\AFF{Department of Statistics and Operations Research,
University of North Carolina at Chapel Hill, \EMAIL{caner@unc.edu}}

\AUTHOR{Mo Liu}
\AFF{Department of Statistics and Operations Research,
University of North Carolina at Chapel Hill, \EMAIL{mo\_liu@unc.edu}}
 
}

\ABSTRACT{
We revisit the multi-period newsvendor network problem, in which demands from multiple customers are correlated and jointly time-varying. Due to the curse of dimensionality associated with estimating the full joint demand distribution, we consider fluid approximation, a widely used approach for solving two-stage stochastic optimization problems such as large-scale service-system design. However, replacing the underlying random distribution (e.g., the demand distribution) with its mean (e.g., the time-varying average arrival rate) introduces bias in performance estimation and can lead to suboptimal decisions.

In this paper, we investigate how to identify an alternative point statistic, not necessarily the mean, such that substituting this statistic into the two-stage newsvendor network problem yields an optimal decision. We refer to this statistic as the decision-corrected point estimate (a time-varying arrival rate). Although the critical fractile is well known to be the decision-corrected point forecast for the single-item newsvendor problem, counterexamples show that such a point statistic may not exist for newsvendor networks. We establish necessary and sufficient conditions for the existence of such a corrected point estimate and propose an algorithm for computing it. Numerical experiments on real data demonstrate that using the proposed decision-corrected point forecast in fluid approximation achieves substantially lower cost than traditional fluid approximation and sample average approximation benchmarks.
}


\maketitle

\section{Introduction}\label{sec:Intro}

As one of the most fundamental problems in supply chain management, the single-item, single-period newsvendor problem can be formulated as
\begin{equation}\label{eq:nv-intro}
    \min_{b \geq 0} \; cb + \mathbb{E}_{D \sim F}\left[
    \min_{x\geq 0} \; p(D - x)
    \;\;\text{s.t.}\;\; x \le b,\; x \le D
    \right],
\end{equation}
where $b$, $x$, $D$, $p$, and $c$ denote the ordering quantity, satisfied demand, realized demand, unit lost-sales cost, and unit ordering cost, respectively. When $c < p$, the optimal capacity has the well-known closed-form solution $b^* = F^{-1}(1-c/p)$,
which is the $(1-c/p)$-quantile of the demand distribution.

This paper studies a multi-resource, multi-customer, multi-period version of the newsvendor problem by extending the variables in \eqref{eq:nv-intro} to multiple dimensions. The problem is formulated as
\begin{equation}\label{eq:multi-period-intro}
    \min_{b \geq 0} \; c^{\top} b \, T + \mathbb{E}_{\{D_t\}_{t=1}^T\sim \mathcal{D}}\!\left[
    \min_{x_t \geq 0} \; \sum_{t=1}^{T} p^{\top}(D_{t} - Rx_{t})
    \;\;\text{s.t.}\;\; Ax_{t} \leq b,\; Rx_{t} \leq D_{t} \;\; \forall t
    \right].\tag{P1}
\end{equation}

In \eqref{eq:multi-period-intro}, there are $n$ customer classes, $m$ resource types, and $T$ time periods. The decision vector $b\in \bbR^m$ represents the resource capacity levels, such as staffing capacities, which remain constant over the planning horizon. The vector $c\in\bbR^m$ denotes the unit procurement cost of these resources per period. In the inner problem, for each period $t$, the demand vector $D_t\in\bbR^n$ is assigned to resource pools through the matching decision $x_t$, whose dimension depends on the connectivity between customer classes and resource pools. The routing matrix $R$ and the capacity requirement matrix $A$ encode which resource pools can serve which customer classes and how much resource capacity is required for each unit of demand. See Section \ref{sec:ProblemSetup} for a detailed explanation of the formulation.

Problem \eqref{eq:multi-period-intro} is important for two reasons. First, it extends the classical newsvendor model to a newsvendor network with time-varying demand; see, e.g., \cite{mieghem2002newsvendor,van2003commissioned}. By specifying the matrices $R$ and $A$, the model captures how managers allocate capacities across heterogeneous resources while facing stochastic and fluctuating demand, heterogeneous resource usage rates, and different service capabilities. In a general service network (such as call centers and healthcare facilities), the optimal design needs to balance customer demand against limited capacity. These systems involve multiple customer classes, heterogeneous resources or agents with diverse skills, and time-varying demand. Operational decisions are typically made in two stages: in the first stage, managers determine the staffing/procurement capacity (how much resource to purchase); in the second stage, real-time routing decisions are made to assign customers to available resources. The total cost comprises the first-stage procurement cost and the second-stage cost associated with unmet demand.

Second, the two-stage structure in \eqref{eq:multi-period-intro} is closely related to service-system design through the lens of fluid approximation for large queueing systems, as discussed below.

\subsection{Connection between \eqref{eq:multi-period-intro} and fluid approximation}

The second motivation for studying Problem \eqref{eq:multi-period-intro} comes from the fluid approximation for large queueing systems; see, e.g., \cite{harrison2005method,bassamboo2010optimal}. Under this interpretation, the capacity vector $b$ represents the number of servers in each server pool, $p$ represents the unit abandonment rate multiplied by the unit abandonment cost, matrix $R$ represents the service-capability matrix, and matrix $A$ represents the reciprocals of the service rates. The online matching between servers and customers is characterized by the matrix $x_t$. See \cite{kerimov2024dynamic} for details. The inner problem provides a full-information, hindsight lower bound on the lost-service or abandonment cost in each period. The discrete-time formulation can also be extended to continuous time by modeling the fluid dynamics of the arrival process; see, for example, \cite{bai2022fluid} and Appendix \ref{append:equivalence}.

The fluid approximation is one of the most common approaches in the design of large-scale service systems, including emergency facilities, call centers, hospitals, and supply chains. These systems are challenging to optimize because they involve uncertainty in customer arrivals, patient types, service times, job requirements, consumer behavior, and other random factors.
The fluid approximation method replaces the joint demand distribution with a time-varying demand arrival rate. In particular, it treats the arrival rate as the exact realization of demand, so the problem of determining the number of servers in each pool becomes a deterministic optimization problem \citep{harrison2005method,whitt2002stochastic}. 

Using the average behavior of the random variables is computationally attractive because it reduces a complex stochastic system to a deterministic one in which uncertain quantities are replaced by their means. This approximation is intuitive because the mean often captures substantial information about the system without requiring the full demand distribution, as discussed in \cite{bassamboo2006design}. Fluid approximation is widely adopted and often performs well in large-scale service-system design when the uncertainty variation is small \citep{mills2013resource,kim2019modeling,hu2021prediction,kerimov2025optimality}. Beyond service systems, fluid approximation is also widely used to study the evolution of cellular systems, diseases, and other stochastic systems; see, for example, \cite{banerjee2022heavy,pakdaman2010fluid,OlveraCravioto2011Transition}.

\subsection{Why a decision-corrected point forecast?}

The main challenge in solving \eqref{eq:multi-period-intro} is that the future joint demand distribution $\mathcal{D}$ across all periods and products is unknown. If $\mathcal{D}$ were known, then \eqref{eq:multi-period-intro} would be convex in the capacity decision $b$ and could be solved tractably using convex optimization through sample average approximation (SAA). In practice, however, future demand is unknown and may exhibit nonstationarity, seasonality, and dependence on contextual features. Estimating the full joint distribution $\mathcal{D}$ is often intractable, since demand can be correlated across products and time periods. Treating demand as independent across products or periods is usually unrealistic, while estimating the full joint distribution suffers from rapidly growing statistical complexity as the number of products and periods increases. This statistical burden is also reflected numerically in Figure \ref{fig:total} in Section \ref{sec:numerical}, where, in the presence of temporal effects, a benchmark SAA method requires a training set roughly 20 times larger than that of the proposed method to achieve the same cost.

Given the difficulty of estimating the future joint demand distribution, a natural and practically important question is whether a fluid-approximation idea can still be used for future capacity decisions. In particular, can we use a time-varying arrival rate as a surrogate for the unknown distribution, and plug it into Problem \eqref{eq:multi-period-intro} to determine the optimal capacity decision?

Although it is conventional to use the per-period average demand (such as the hourly average demand) as the fluid approximation input (i.e., the Poisson arrival rate within that hour), this approximation remains an approximation: Replacing a random distribution by its expectation removes substantial information about uncertainty.  The system performance evaluated under mean demand, such as cost or throughput, is not necessarily equal to the expected performance under the original demand distribution. Consequently, the system design that is optimal under the mean demand arrival rate need not be optimal for the true stochastic system. Traditional literature, such as \cite{bassamboo2009data}, justifies fluid approximation by showing that this gap vanishes asymptotically as the scale of the system grows. In finite-scale systems, however, the gap can be non-negligible or even substantial, especially when demand variance is large and the mean demand is small. In such settings, standard fluid approximation can lead to biased and suboptimal service-system designs. To the best of our knowledge, relatively little work studies this gap systematically or develops methods to correct it. 

In this paper, we aim to identify an alternative statistic of the time-varying joint demand distribution. This statistic need not be the mean, but substituting it into the two-stage optimization problem should consistently recover the optimal capacity decision. We refer to this statistic as the \textit{decision-corrected demand arrival rate}. Specifically, we study the following question:

\textit{\bfseries Given the two-stage stochastic optimization problem \eqref{eq:multi-period-intro}, does a decision-corrected demand arrival rate exist for any demand distribution? If so, how can it be forecast from data?}

Such a decision-corrected statistic may not exist; see Examples \ref{counterexample} and \ref{ex:t1-stochastic}. When the decision-corrected arrival rate exists, the pipeline for solving \eqref{eq:outer_problem} and determining future capacity sizing is illustrated in Figure \ref{fig:workflow}. In Figure \ref{fig:workflow}, historical arrival data are used to construct the decision-corrected arrival rate using the methods introduced in Section \ref{sec:correction}. The resulting historical decision-corrected arrival rates are then used to forecast future decision-corrected arrival rates based on prediction models, such as time-series models. In the final step of Figure \ref{fig:workflow}, minimizing the total cost under the forecasted decision-corrected arrival rates can be shown to be equivalent to minimizing the expected cost under the true demand distribution. This workflow provides a way to determine optimal staffing levels while avoiding the need to estimate a high-dimensional future joint demand distribution. 

\begin{figure}[htbp]
\centering
\begin{tikzpicture}[
    node distance=1.2cm and 2.5cm, 
    every node/.style={align=center},
    box/.style={
        draw,
        rectangle,
        rounded corners,
        minimum height=1.2cm,
        minimum width=2.8cm
    },
    arrow/.style={->, thick}
]

\node[box] (arrival) {Historical \\ demand data};

\node[box, right=1.9cm of arrival] (p2) {Decision-corrected \\arrival rate in history};

\node[box, below=of p2] (p1) {Decision-corrected arrival rate in future};

\node[box, right=1.3cm of p1] (queue) {Two-stage staffing problem \eqref{eq:outer_problem}\\ (Subsumes multi-resource\\ multi-customer newsvendor)};

\draw[arrow] (arrival) -- node[above]{Abstract}node[below]{} (p2);
\draw[arrow] (p2) -- node[right]{forecast} (p1);
\draw[arrow] (p1) -- node[above]{Plug-in}node[below]{} (queue);

\end{tikzpicture}
\caption{Using forecasted decision-corrected arrival rates for future capacity decisions.}
\label{fig:workflow}
\end{figure}

The pipeline in Figure \ref{fig:workflow} also falls within the \textit{predict-then-optimize} paradigm. Traditional forecasting approaches estimate demand without accounting for the subsequent two-stage decision-making problem in which the forecast is used. However, in the workflow illustrated in Figure \ref{fig:workflow}, a more accurate prediction, in the sense of being closer to the average demand arrivals, does not necessarily lead to a better procurement decision. Therefore, the proposed data-driven method for constructing the decision-corrected arrival rate provides a principled way to align prediction with decision quality. The main contributions of this paper are summarized as follows:

\begin{itemize}

    \item We are the first to study whether a decision-corrected point statistic exists for the multi-period newsvendor network problem in \eqref{eq:outer_problem}. Counterexamples show that such a point statistic may fail to exist even in a single-period network with two customer classes and three resources.

    \item We are the first to provide necessary and sufficient conditions on the parameters of the service network for the universal existence of the decision-corrected arrival rate for any demand distribution.

    \item From a data-dependent perspective, we provide necessary and sufficient conditions on the observed demand data for the distribution-dependent existence of the decision-corrected arrival rate. We propose an algorithm to verify the existence.

    \item When the decision-corrected arrival rate exists, we develop an algorithm to construct it directly from data and establish the consistency of the proposed method. This correction method provides \textit{insights into how to move beyond averaged demand and align arrival-rate estimation with downstream decisions in service system design.}
    
    \item Numerical experiments using real patient arrival data from a hospital demonstrate that the prediction model based on the decision-corrected arrival rate achieves a lower total cost than three benchmarks: the traditional mean-based fluid approximation, a stationary SAA benchmark that pools historical data, and a nonstationary SAA benchmark that incorporates time-series forecasting.
\end{itemize}

The remainder of the paper is organized as follows. Section \ref{sec:LitRev} reviews related work in fluid approximations and predict-then-optimize. Section \ref{sec:model} introduces the two-stage model, establishes the notation, presents the fluid model that serves as our reference point, and illustrates the decision bias with a simple example. Section \ref{sec:existence} develops the notion of a decision-corrected arrival rate and establishes the conditions for its existence and non-existence, building on an equivalence to a deterministic expanded-scenario formulation, and verifying KKT optimality conditions. Section \ref{sec:correction} presents a data-guided method for constructing decision-corrected arrival rates that applies to general service networks, and shows that for decomposable networks the construction simplifies to a closed-form demand quantile. Section \ref{sec:numerical} reports numerical experiments comparing the decision-corrected approach with the traditional fluid approximation and sample average approximation benchmarks. Conclusions are given in Section \ref{sec:Conclusion}. Omitted proofs, additional discussions, and details for the numerical experiments appear in the appendices.

\section{Literature Review}\label{sec:LitRev}

The decision-focused bias-correction method in this paper is proposed to solve the two-stage stochastic optimization problem. To position our paper, we first review the connection between the fluid approximation methods and the predict-then-optimize regime in the general two-stage stochastic optimization problems. Since one of the motivations of fluid approximation is to solve the staffing problem in service systems, we next review the literature that uses the popular fluid approximation method for solving such problems.
To the best of our knowledge, we are the first to address the bias-correction method for fluid approximation.

\subsection{Predict-then-optimize methods} 

 Using fluid approximations for decision making, our setting falls within the predict-then-optimize paradigm, where estimated arrival rates are first learned and then used to optimize service system design. This framework subsumes both the newsvendor problem and two-stage resource allocation problems.

\paragraph{Point forecast for the newsvendor problem.}
Incorporating operational costs into demand estimation for the single-product newsvendor problem has been widely studied. For example, \cite{liyanage2005practical}, \cite{feng2022developing}, and \cite{feng2025contextual} consider point statistics in the setting of operational data analytics approaches, while \cite{iceo} studies polynomial approximation methods. \cite{homem2024forecasting} characterizes optimal pointwise forecasts for general two-stage stochastic optimization problems. However, their formulation does not cover the multi-resource newsvendor problem, as shown in Example \ref{ex:t1-stochastic}.
SAA methods are studied in \cite{levi2015data} and \cite{ban2019big} for non-contextual and contextual settings, respectively, with finite-sample error analyses in \cite{siegel2021profit, siegel2023data}. However, the existing literature focuses exclusively on single-product settings. In contrast, our work is the first to study how the structure of multi-product newsvendor problems affects the existence of decision-corrected estimators.

\paragraph{Two-stage resource allocation and scenario reduction.}
Arrival-rate estimation has also been extensively studied in healthcare and call center systems; see, for example, \cite{green2013nursevendor}, \cite{gans2015parametric},  \cite{ibrahim2016modeling}, and \cite{hu2021prediction, hu2022optimal}. Likelihood-based estimation methods are considered in \cite{oreshkin2016rate}, while demand uncertainty and safety-stock design are studied in \cite{bassamboo2010capacity}. Data-driven staffing policies for multi-class systems are proposed in \cite{bassamboo2009data}. 
More recent work considers joint estimation and allocation under fluid approximations \citep{kanoria2024blind, zheng2025joint,zhong2025learning}. Fluid-based policies have also been shown to outperform traditional heuristics in applications such as mass-casualty triage \citep{mills2013resource} and limited-inventory allocation across multiple customer classes \citep{ge2025optimal, tang2025split}. Relatedly, \cite{bertsimas2023optimization} studies scenario reduction for SAA from a sampling perspective, showing that the number of scenarios can be reduced by over 90\% while achieving comparable performance.

\paragraph{General decision-focused learning. } In the general decision-focused learning literature, papers have designed different methods to obtain a prediction that results in the optimal downstream decision-making problem (e.g., \cite{elmachtoub2022smart,huang2024decision,liu2023active,albert2025post}). The vast majority of this literature considers uncertainty only in the objective function, because uncertainty in the feasible region (i.e., the constraints), even when it appears only in the right-hand side of the linear constraints, introduces significantly greater challenges (see Appendix \ref{append:pitfall}). Closer to our setting, \cite{estes2023smart} handle right-hand side uncertainty in two-stage linear programs with side information by training a linear regression predictor under a smart predict-then-optimize loss that approximates the regret of the first-stage decision induced by the predicted right-hand side after second-stage recourse, while \cite{hu2023unknownconstraint} handle unknown constraint parameters in mixed integer linear programs by training an end-to-end predictor under a post-hoc regret loss that penalizes correcting an infeasible first-stage solution by a second-stage recourse optimization. These works are methodological in nature: they propose training procedures that reduce decision error on average. Our contribution is complementary and structural. Rather than fitting a predictor to a particular loss, we ask whether there exists a point statistic of the demand distribution whose plug-in fluid solution is itself optimal for the two-stage problem, and we provide necessary and sufficient conditions for its existence together with a data-driven construction.

To the best of our knowledge, we are the first to propose a decision-unbiased method to estimate parameters in uncertain feasible regions, formed by the resource allocation problem.

\subsection{Fluid approximation methods in service system design}
Although the fluid approximations studied in this paper do not directly capture the dynamics in the queueing network, they serve as a good approximation when determining the capacity of servers. \cite{harrison2005method} is the seminal paper in this area that considers the call center staffing problem. It uses fluid approximations to simplify multidimensional scheduling problems with customer abandonment. Particularly, this multi-period staffing problem reduces to a multi-product newsvendor-type problem under fluid approximation. This line of work was extended to doubly–stochastic arrivals with an LP–based staffing and routing method in \cite{bassamboo2006design}, where asymptotic optimality in heavy traffic is established. \cite{gamarnik2006validity} shows that the fluid model approximates the steady state in the heavy traffic regime in generalized Jackson networks. \cite{bassamboo2010accuracy} further establishes the asymptotic optimality and quantifies the gap between the fluid approximation and the optimal capacity decisions in systems with impatient customers.  Existence and uniqueness of a fluid model for time-varying many-server queues with abandonment are established in \cite{kang2014existence}. \cite{banerjee2022heavy} further demonstrates that heavy-traffic approximations must account for service-time heterogeneity.

Another important class of decisions based on fluid approximation arises in dynamic matching. \cite{gurvich2015dynamic} studies pairing customers with servers over time under compatibility constraints, where fluid approximation is used to derive effective policies.
\cite{kerimov2024dynamic} further considers a dynamic matching model and demonstrates that a periodic-clearing policy—periodically solving an integer linear program achieves a constant regret. \cite{kerimov2025optimality} further shows that suitably designed greedy policies (e.g., longest-queue and static priority based on residual networks) based on fluid models achieve small all-time regret.

Given the aforementioned broad applicability of fluid approximations, our decision-corrected point estimation provides new insights into how to design and estimate such fluid models across a range of applications.

\section{Two-Stage Decision-Making Model}\label{sec:model}

In Section \ref{sec:ProblemSetup}, we present the formulation of our two-stage model for determining the capacity and allocations under random demand arrivals. To solve this problem, in Section \ref{sec:motivation_fa}, we review the common fluid approximation and demonstrate why it is biased for decision-making.

\subsection{Problem setup}\label{sec:ProblemSetup}

We formally introduce Problem \eqref{eq:outer_problem}, which serves multiple customer types using multiple resource or agent pools. Each customer type, or class, may be served by one or more resource pools, and each resource pool may be able to serve more than one customer class. This structure naturally arises in settings such as call centers, healthcare systems, and technical support environments.

\begin{figure}[H]
  \centering
  \begin{subfigure}[t]{0.49\linewidth}
    \centering
    \includegraphics[width=\linewidth]{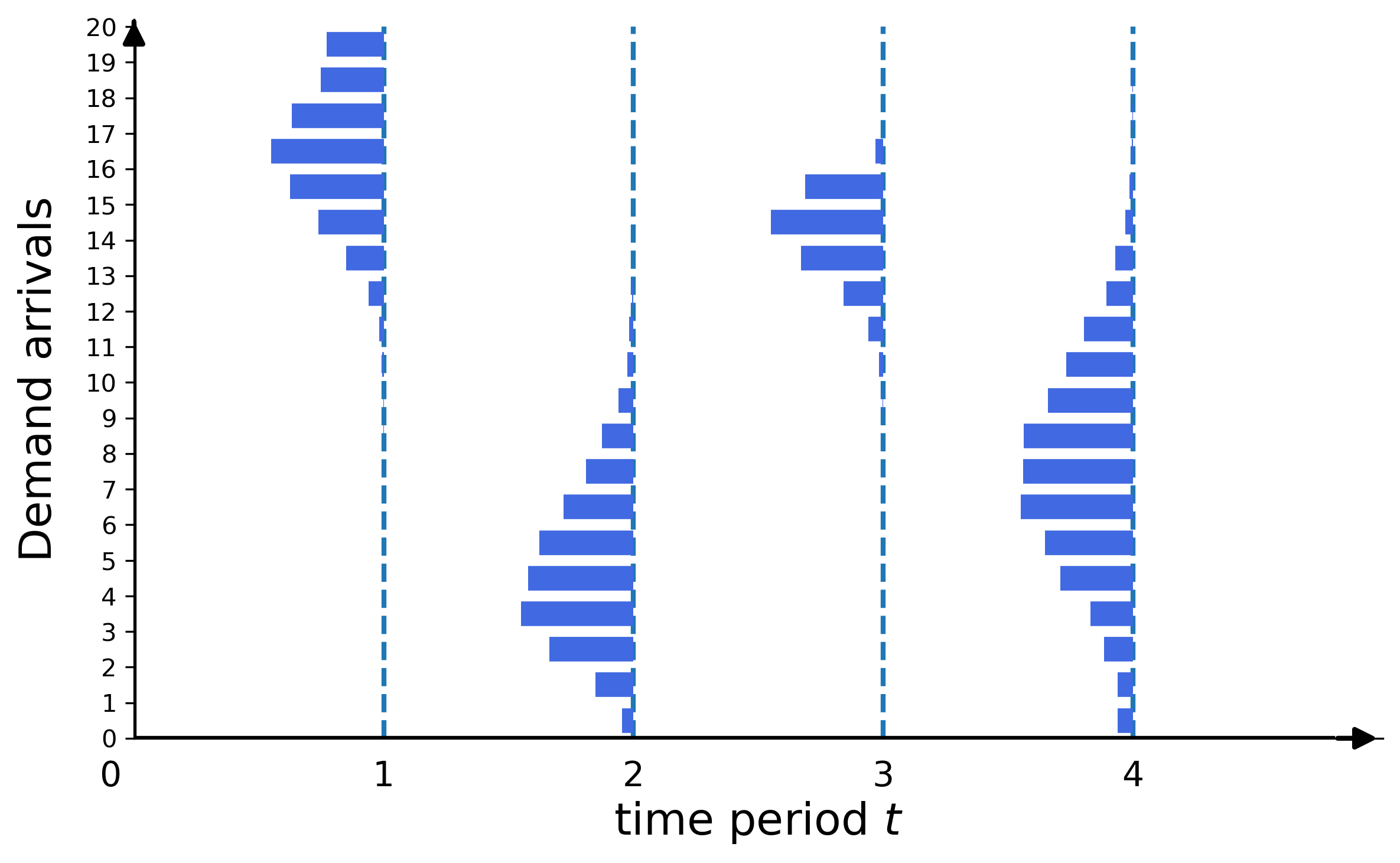}
    \subcaption{Observed frequency of demand over time.}
    \label{fig:example1:hist}
  \end{subfigure}
  \hfill
  \begin{subfigure}[t]{0.49\linewidth}
    \centering
    \includegraphics[width=\linewidth]{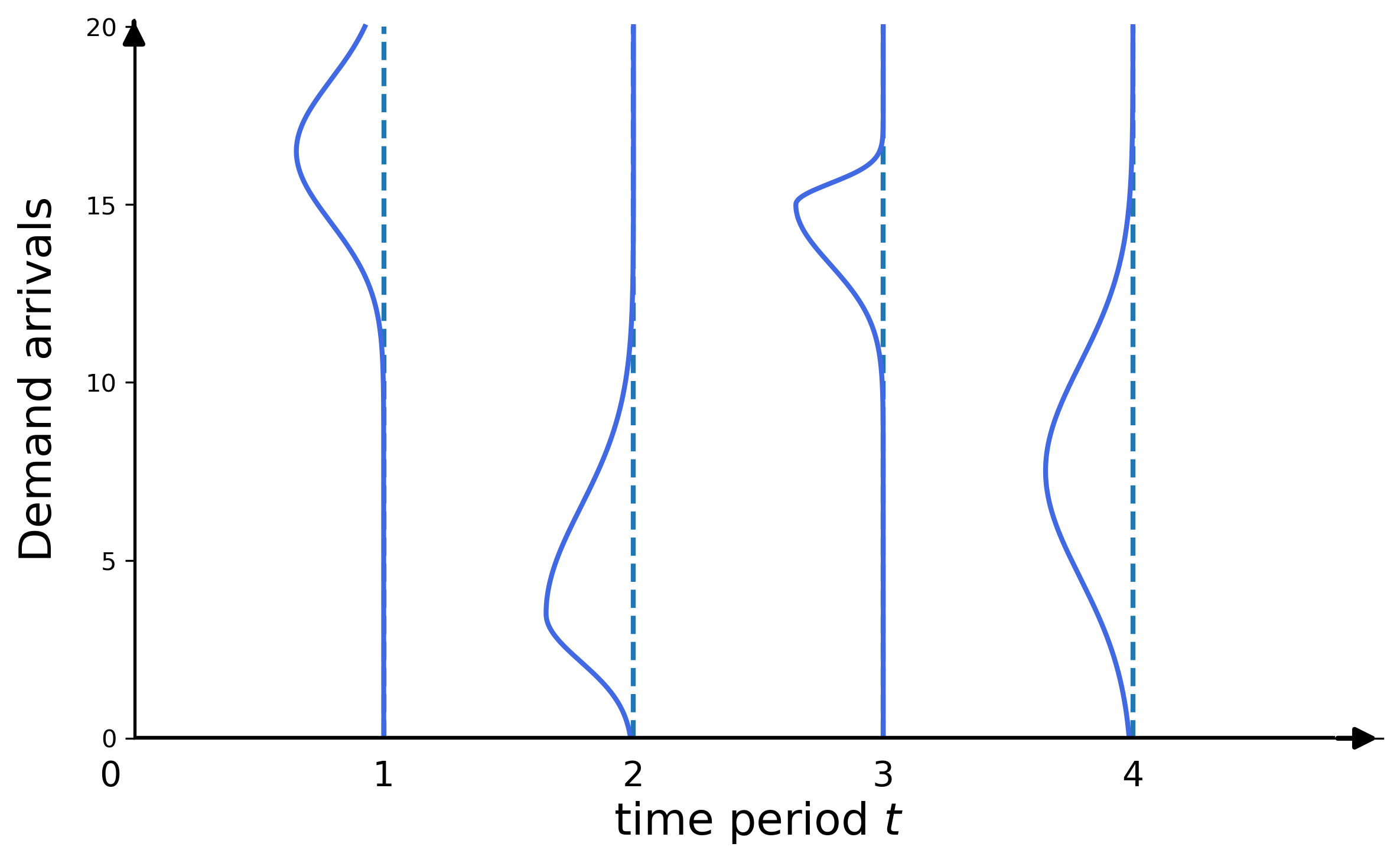}
    \subcaption{Estimated demand distribution over time.}
    \label{fig:example1:density}
  \end{subfigure}
  \caption{Example of demand arrivals for one customer type.}
  \label{fig:example1}
\end{figure}

Let $n$ denote the number of customer classes, and $m$ denote the number of distinct resource pools. We divide the planning period (e.g., a workday) into $T$ discrete time intervals, indexed by $t = 1, \dots, T$ (e.g., hourly blocks). For each time period $t$, let $D_t \in \mathbb{R}_+^n$ denote the vector of demand for $n$ customer classes, where $\mathbb{R}_+:=\{x\ge0,x\in \mathbb{R}\}$. The vector $\bbD = (D_1, \dots, D_T)$ represents the demand profile over the full planning period. In this paper, we assume the demand $D_t$ is time-varying across different periods within the planning period. The joint distribution of $\bbD$, which is denoted by $\mathcal{D}$, is fixed and bounded but unknown to the decision-maker. However, the decision-maker can collect i.i.d. observations of $\bbD$ from $\mathcal{D}$. Figure \ref{fig:example1} provides one example of the observed demand and estimated distribution $\mathcal{D}$ for one customer type.

In the example shown in Figure~\ref{fig:example1}, the planning period is one day. We display the first four hours of the day as the time periods (\(t=1,2,3,4\)). For hour~1, one might observe 17 arrivals on day~1, 20 arrivals on day~2, and so on. Aggregating these arrivals across multiple days yields the observed frequency of demand for hour~1, which is plotted vertically at \(t=1\) in Figure~\eqref{fig:example1:hist}. Similar observed demand frequencies are plotted at \(t=2,3,4\). Based on these observed frequencies, the decision-maker can then estimate the demand distribution, as illustrated in Figure~\eqref{fig:example1:density}. Note that in this paper, the joint distribution $\mathcal{D}$ considers $n$ types of customers and allows for \textit{dependence and correlation} in demand both across customer types and across time periods.

We note that the discrete-time formulation above is a modeling choice and not a substantive restriction. The framework accommodates per-period totals generated by arbitrary underlying arrival processes, including continuous-time Poisson-type processes aggregated over each epoch. In particular, our setup is not restricted to the low-variance regime tied to the classical Bernoulli-process fluid approximations; as \cite{bai2022fluid} show, discrete-time fluid models can accommodate arrival processes with arbitrary mean and variance.

Based on the estimated demand distribution $\mathcal{D}$, in the first stage, the decision-maker must choose a procurement vector $b \in \mathbb{R}_+^m$ for $m$ types of resources, where $b_h$ represents the number of resources procured of type $h$. Procurement decisions are made at the start of the planning period and remain fixed across all time intervals. Here, $b_h$ represents the service capacity of resource pool $h$, sometimes referred to simply as capacity, meaning the maximum amount of service pool $h$ can provide in each period. 
Once procurement is fixed, in the second stage, the actual demand $\bbD$ is observed. Based on the observed demand, resources are assigned to handle incoming customer demand through a set of routing decisions that can vary over time. The structure of this service system is illustrated as a bipartite graph in Figure \ref{fig:bipartite}.

\begin{figure}[H]
\centering
\begin{tikzpicture}[
    x=1cm,y=1cm,>=stealth,
    square/.style={draw,minimum size=6mm,inner sep=0pt},
    circ/.style={draw,circle,minimum size=6mm,inner sep=0pt},
    every node/.style={font=\sffamily\footnotesize} 
]

\def\ytop{0.6}   
\def\ybot{-3.0}

\node[square] (c1) at (0,\ytop) {1};
\node[square] (c2) at (2,\ytop) {2};
\node[square] (c3) at (4,\ytop) {3};
\node            at (5.5,\ytop) {$\cdots$};
\node[square] (cn) at (7,\ytop) {$n$};

\node[circ] (s1) at (0,\ybot) {1};
\node[circ] (s2) at (2,\ybot) {2};
\node            at (4,\ybot) {$\cdots$};
\node[circ] (sm) at (7,\ybot) {$m$};

\draw (c1) -- (s1);
\draw (c2) -- (s2);
\draw (c3) -- (s2);
\draw (c3) -- (s1);
\draw (cn) -- (s2);
\draw (cn) -- (sm);

\node[align=center,anchor=east] at (-1,{(\ytop+\ybot)/2}) { edges of\\feasible allocations};

\node[anchor=west] at (8.0,\ytop) {$n$ types of customers};
\node[anchor=west] at (8.0,{\ytop-1.3}) {Binary routing matrix ${R}_{n\times k}$};
\node[anchor=west] at (8.0,{\ybot+1.3}) {Capacity requirement matrix ${A}_{m\times k}$};
\node[anchor=west] at (8.0,\ybot) {$m$ types of resources};

\end{tikzpicture}
\caption{Bipartite representation of $n$ customer types and $m$ resource types. Edges indicate feasible job assignments. The binary routing matrix $\mathbf{R}$ captures possible assignments; the capacity requirement matrix $\mathbf{A}$ captures resource usage.}
\label{fig:bipartite}
\end{figure}

In the second stage, when assigning customers to different resources, the feasibility of assignments is represented by the routing matrix $R \in \{0, 1\}^{n \times k}$, where $k\le n \times m$ is the total number of possible allocations. The allocation from one type of customer to one type of resource is also referred to as one job.  In Figure \ref{fig:bipartite}, the feasible allocations are represented by the edges connecting customers and resources in a bipartite graph. In the routing matrix $R$, each entry $R_{ij}$ is a binary indicator of whether the edge corresponding to job $j$ involves customer class $i$. In other words, $R_{ij} = 1$ if customer $i$ can be served by job $j$, and $R_{ij} = 0$ otherwise.

The second-stage decision is represented by the allocation vector $x \in \mathbb{R}^k_+$, where each component $x_j$ denotes the number of customers assigned to job $j$. Accordingly, the product $R x$ represents the total number of served customers (across all $n$ customer types) under the allocation decision $x$.

We further define a capacity requirement matrix $A \in \mathbb{R}_+^{m \times k}$, where $A_{hj}$ represents the amount of capacity required from resource pool $h$ when \textit{one unit of job $j$} is used. Therefore, the product $A x$ represents the consumed resource for $m$ types of resources under the allocation decision $x$.

Since each edge in Figure~\ref{fig:bipartite} connects one customer type and one resource type, both matrices $A$ and $R$ have exactly one non-zero element in each column. We introduce index mappings to formalize these relationships. For each job (column)~$j$ in either matrix~$A$ or~$R$, let $i(j)$ denote the customer class served by job~$j$, and let $h(j)$ denote the resource pool from which job~$j$ draws capacity. Hence, each job~$j$ serves exactly one customer class~$i(j)$ and utilizes capacity from exactly one resource pool~$h(j)$.

Given the service network, the job allocation vector $x$ must satisfy the following component-wise constraints for all time periods $t \in 1,\dots, T$:
$$R x_t \leq D_t, \text{ and } A x_t \leq b.$$
The first set of inequalities, $R x_t \leq D_t$, ensures that the total service allocated to each customer class does not exceed its arrival rate. The second set of inequalities, $A x_t \leq b$, ensures that the total service assigned to each resource pool does not exceed its available capacity.

If the available capacity is not sufficient to serve all incoming demand, we assume this unsatisfied demand is lost. To capture the cost of this, we define a penalty vector $p \in \mathbb{R}_+^n$, where $p_i$ represents the cost associated with each unit of unserved demand from customer class $i$. Given a job allocation $x$, the total lost demand cost is given by $\sum_{t=1}^T p^\top (D_t - R x_t)$. Thus, the total cost incurred by the system consists of two components:
\begin{itemize}
    \item \textbf{Total Procurement cost}: $c^\top b ~ T$, where $c \in \mathbb{R}_+^m$ is the cost per unit of capacity for each resource per period;
    \item \textbf{Total Lost Demand cost}: $\sum_{t=1}^T p^\top (D_t - R x_t)$, which arises when demand is not fully met.
\end{itemize}

To determine the best capacity $b$ for all types of resources, given the distribution of demand $\mathcal{D}$, we solve the following two-stage stochastic optimization problem:

The first stage problem is to determine the capacity $b$:
\begin{equation}\label{eq:outer_problem}
\min_{b \ge 0} \quad c^\top b ~ T+ \bbE_{\bbD \sim \mathcal{D}} \left[ \pi(b, \bbD) \right].\tag{P1}
\end{equation}
The second stage problem $\pi: \bbR^m \times (\mathbb{R}^n_+)^T  \rightarrow \bbR$ determines the allocation $x_t$ for $t = 1,\dots,T$:
\begin{equation}\label{eq:inner_problem}
\pi(b, \bbD) := 
\begin{array}[t]{l}
\displaystyle \min_{\{x_t \ge 0\}} \sum_{t=1}^T p^\top (D_t - R x_t) \\
\text{s.t.} \quad  
    \begin{aligned}[t]
    & A x_t \le b, \quad \text{for all } t = 1, \dots, T, \\
    & R x_t \le D_t, \quad \text{for all } t = 1, \dots, T.
    \end{aligned}
\end{array}
\end{equation}
This model is equivalent to the popular model in \cite{harrison2005method} that has wide applications in operations. Below are three examples:

\paragraph{Application 1: Facility Location and Transshipment.}
Problem \eqref{eq:outer_problem} can be used to determine the optimal locations and capacities of facilities, such as warehouses or distribution centers; see, e.g., \cite{mieghem2002newsvendor}. In the first stage, a company makes a strategic decision about where to open facilities and how much capacity to install, represented by $b$, incurring a setup cost $c^{\top} b T$. This decision is made before future customer demand is realized. In the second stage, after demand $D_t$ from different regions is observed, the company makes operational decisions $x_t$ on how to route products from available facilities to satisfy demand. Geographic and service constraints are captured by the routing matrix $R$.

\paragraph{Application 2: Workforce Staffing in Large-Scale Service Systems.}
Another application of this framework is staffing in large-scale service systems, such as hospital emergency departments or multi-skill call centers; see e.g., \cite{bassamboo2010optimal,ibrahim2016modeling}. In these systems, managers make a first-stage capacity decision, such as the number of available nurse-hours or the number of agents in each skill pool, represented by $b$. These decisions incur a setup cost $c^\top b T$ and must be made before the actual volume of patient or customer arrivals is known. In the second stage, as requests arrive in real time, they are routed to available resources through the assignment decision $x_t$. Following \cite{harrison2005method}, customers waiting in queue may abandon the system before receiving service, resulting in a penalty $p^\top(D_t - R x_t)$ that reflects costs such as patient boarding, customer abandonment, or reduced service quality. In an asymptotic regime, the performance of this system can be well approximated by \eqref{eq:outer_problem}.

\paragraph{Application 3: Manufacturing and Energy Capacity Planning.}
Problem \eqref{eq:outer_problem} also applies to manufacturing and energy capacity planning; see e.g., \cite{van2003commissioned,devalve2020primal}. In manufacturing, a firm decides in advance how much production or resource capacity $b$ to reserve for multiple products, incurring an upfront cost $c^\top b T$. After uncertain demand $D_t$ is realized, the second-stage decision $x_t$ allocates available capacity to assemble and satisfy demand. Any unmet demand incurs a penalty $p^\top(D_t - R x_t)$, representing lost sales, backlog costs, or expedited production costs.
A similar structure arises in energy systems, where base generation capacity is committed in advance based on forecasted household and industrial demand. After actual demand is realized, generated power is dispatched across the grid. If the committed generation is insufficient, grid operators must rely on expensive peaking plants or emergency spot-market purchases, leading to significant penalty costs.

\subsection{Fluid approximation for predicting unknown demand}\label{sec:motivation_fa}

As motivated in the introduction, the main challenge in solving \eqref{eq:outer_problem} comes from the unknown joint distribution $\mathcal{D}$ of all customer types across all periods. Thus, we consider a fluid-approximation-type method. In the traditional fluid approximation, rather than estimating the entire distribution $\mathcal{D}$, the decision-maker only predicts the mean demand, $\bbE[D_1], \bbE[D_2], \ldots, \bbE[D_T]$. 
The expectation of demand is conventionally referred to as the arrival rate, which is denoted by the vector $\bblambda:=(\bbE[D_1], \bbE[D_2],...,\bbE[D_T]).$ 
Substituting the demand distribution with $\bblambda$, the two-stage problem becomes the following deterministic problem:
\begin{align}\label{eq:fa}\tag{P2}
\min_{b \ge 0} \quad & c^\top b ~ T+  \pi(b, \bblambda) ,\\
& \pi(b, \bblambda) = 
\begin{array}[t]{l}
\displaystyle \min_{\{x_t \ge 0\}} \sum_{t=1}^T p^\top (\lambda_t - R x_t) \nonumber\\
\text{s.t.} \quad  
    \begin{aligned}[t]
    & A x_t \le b, \quad \text{for all } t = 1, \dots, T, \\
    & R x_t \le \lambda_t, \quad \text{for all } t = 1, \dots, T.
    \end{aligned}
\end{array}
\end{align}

The fluid approximation, \eqref{eq:fa}, has been widely used in queueing networks, and the predict-then-optimize paradigm, e.g., \cite{harrison2005method, bassamboo2009data, hu2021prediction}. The fluid approximation is also consistent with the original two-stage problem, \eqref{eq:outer_problem}, in the asymptotic regime, i.e., when $\bbE[\bbD]$ goes to infinity and variance is bounded, then the fluid approximation yields the same capacity decision as the original problem; e.g., see \cite{whitt2002stochastic,gurvich2015dynamic}.

However, in practice, when the system size and realized demand are finite (for example, demand represented by a bounded random vector), the fluid approximation may lead to a biased capacity decision. Simply replacing the demand distribution with its mean loses a significant amount of information, although it makes the computation tractable. Below is a simple example of a one-to-one service system, where the service system reduces to a newsvendor problem, and the fluid approximation leads to a suboptimal decision.

\begin{example}[Newsvendor special case: One customer type, one resource system]\label{example:newsvendor}
Consider a simple system with a single customer class, a single resource pool, and a single period, so that all vectors reduce to scalars with $R=A=1$. To avoid triviality, suppose $c<p$; otherwise, purchasing no resource would be optimal. Given random demand $\bbD$, the second-stage problem is
\begin{align*}
\pi^*(b,\bbD)\;=\;\min_{x\ge0}\; p\,(\bbD-x)\quad\text{s.t.}\quad x\le \bbD,~ x\le b.    
\end{align*}
Thus, the first stage problem can be reduced to
$\min_{b\ge0}\;\; c\,b\;+\;\bbE\!\left[p\,(\bbD-b)^+\right]$.
Denoting $F$ as the cumulative distribution function (CDF) of the demand $\bbD$, for a continuous $F$, the objective is convex and the first-order optimality gives
$b^*=
F^{-1}\!\big(\tfrac{p-c}{p}\big)$.

This aligns with the quantile value in the newsvendor problem.
By contrast, the traditional fluid approximation replaces $\bbD$ with its mean $\bar D=\bbE[\bbD]$. The resulting minimizer of the fluid formulation is $b=\bar D$.

However, in general $b^* \neq  \bar D$ (i.e., $F^{-1}(\frac{p-c}{p})\neq \bar D$), so the fluid approximation leads to a suboptimal decision. \hfill\Halmos
\end{example}

Example~\ref{example:newsvendor} shows that, even in a simple setting, using a mean-based arrival rate can yield suboptimal capacity decisions. This suggests that, with multiple customer types and multiple resources, and with more complex routing and service matrices $R$ and $A$, a mean-based fluid approximation may be biased and thus lead to suboptimal choices. Our central question is whether there exists a corrected arrival rate, $\lambda_t \in \bbR^n_+$ a summary statistic of $D_t$, that can be used as the fluid input and still recover the optimal capacity decision.

\begin{figure}[h]
    \centering
    \includegraphics[width=0.7\linewidth]{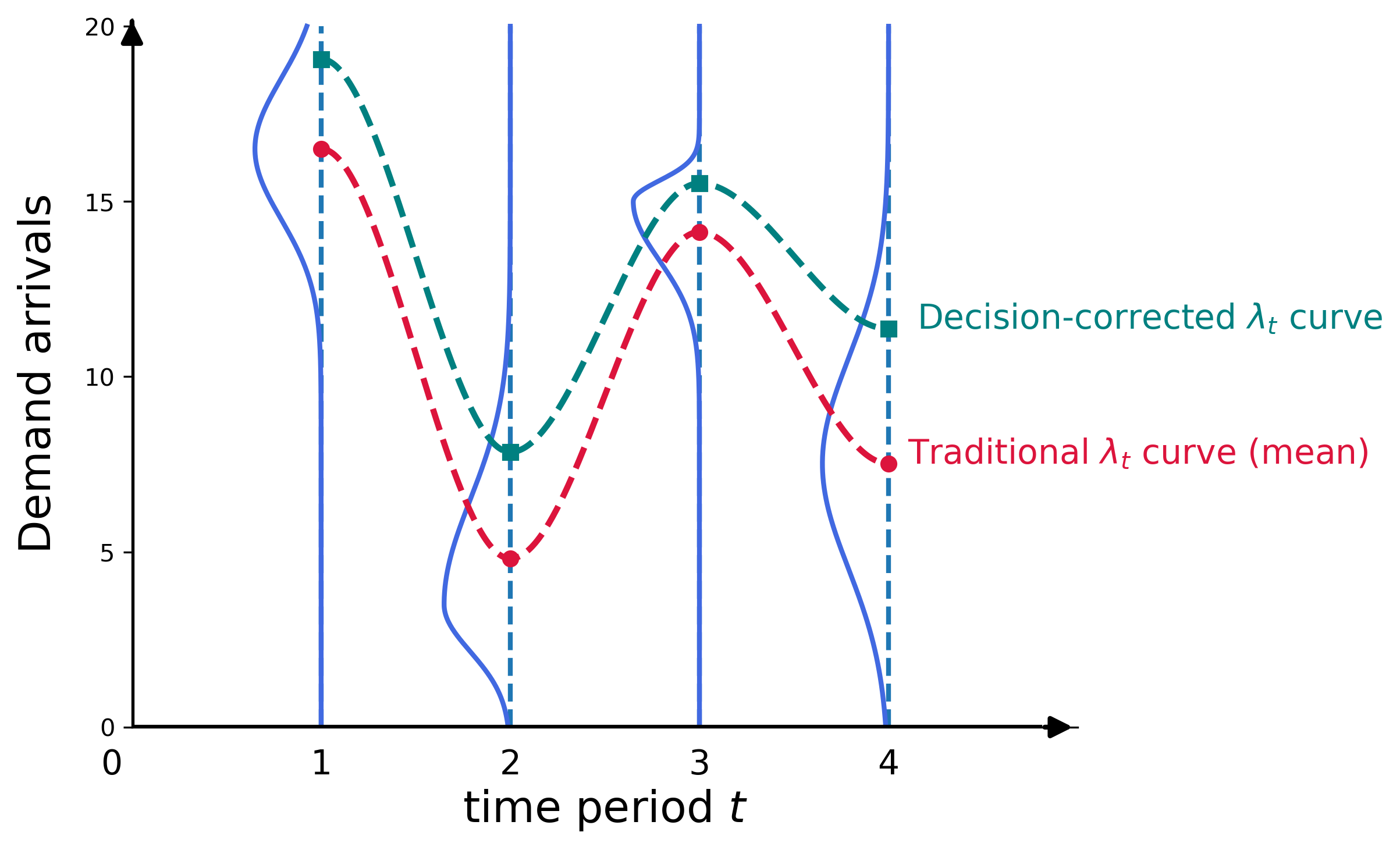}
    \caption{Comparison between the traditional fluid approximation and the decision-focused fluid approximation for one type of customer}
    \label{fig:example2}
\end{figure}

Here, the term ``arrival rate $\lambda_t$'' refers to the input to the fluid approximation; it is not necessarily the average arrival rate. In a predict-then-optimize workflow, this quantity typically corresponds to a point estimate of demand used to guide procurement decisions. Figure \ref{fig:example2} provides one illustrative comparison between the traditional mean-based arrival rate and our \textit{decision-corrected arrival rate} for the input of the fluid approximation in the setting with one type of customer.

Given one arrival rate $\bblambda=(\lambda_1,\ldots,\lambda_T) \in (\mathbb{R}^n_+)^T$, define the set of optimal procurement decisions as
\[
b^*(\bblambda)\ :=\ \arg\min_{b\ge0}\ \Big\{\,c^\top b ~ T\ +\ \pi\big(b,\bblambda\big)\,\Big\}.
\]
When $b^*(\bblambda)$ is not a singleton, we apply a predetermined tie-breaking rule to select a procurement decision from
$b^*(\bblambda)$. Thus, for notational convenience, we assume $b^*(\cdot ): (\mathbb{R}^n_+)^T\mapsto \bbR^m$ is a deterministic mapping to one optimal procurement decision. Then, the decision-corrected arrival rate is defined as follows:

\begin{definition}[Decision-corrected arrival rate]\label{def:arrival_rate}
 We call an arrival rate $\bblambda \in (\bbR_+^n)^T$ decision-corrected with respect to the demand distribution $\mathcal{D}$, if $b^*(\bblambda)\in \arg \min_{b \ge 0} \quad c^\top b ~ T + \bbE_{\bbD \sim \mathcal{D}} \left[ \pi(b, \bbD) \right]$.
\end{definition}

As outlined in Figure \ref{fig:workflow}, using historical decision-corrected $\bblambda$ to forecast future $\bblambda$ is much more efficient than using the entire historical demand distribution to forecast the future distribution. A key motivation for decision-corrected $\bblambda$ is that it provides a simple, low-dimensional label that captures the decision-relevant information in $\mathcal{D}$. Such a label can be used to train prediction models for future (or contextual) arrival rates, and it enables data aggregation or transfer across different environments while ensuring that the resulting procurement decisions remain unbiased.

Although in the simplified setting of Example~\ref{example:newsvendor}, the decision-corrected arrival rate coincides with the newsvendor quantile, this relationship does not generalize to a broader class of service networks. In fact, for general networks, such decision-unbiased arrival rates may not even exist. Therefore, before introducing our bias-correction method for fluid approximation, Section~\ref{sec:existence} investigates whether unbiased arrival rates exist under arbitrary service system structures and demand distributions, and establishes both necessary and sufficient conditions for their existence.

\section{Existence of Unbiased Fluid Approximation}\label{sec:existence}

In this section, we investigate the existence of the decision-corrected arrival rate.  To avoid triviality of the system, we consider the following structure of the service network:

\begin{enumerate}[label=\textbf{(\alph*)}., leftmargin=*]
    \item \textit{Every customer can be served by at least one resource}, i.e.,  
    for each customer class \(i \in \{1, \dots, n\}\), there exists at least one
    job column \(j\) such that \(R_{ij} > 0\).

    \item \textit{Every resource can serve at least one customer class, i.e., }  
    for each resource pool \(h \in \{1, \dots, m\}\), there exists at least one
    job column \(j\) such that \(A_{hj} > 0\).
\end{enumerate}

We use these standard assumptions to keep the model focused on meaningful decisions, as customer classes with no service options only contribute a fixed cost and resource pools that serve no one would always be assigned zero capacity in an optimal solution.

Before studying the existence of the unbiased fluid approximation, we note that the first-stage procurement decision $b$ remains constant across planning periods. Consequently, when searching for the decision-corrected arrival rate, a natural question arises: can we also set the corrected arrival rate $\lambda_t$ to be constant across periods? Intuitively, since the capacity decision is fixed for all periods, we may be able to ignore demand fluctuations and instead summarize them by a constant effective arrival rate. However, Example \ref{counterexample} below shows that in general, we should not expect a constant arrival rate to be decision-corrected. 

\begin{example}[Non-existence of constant decision-corrected arrival rate] \label{counterexample}
Consider a simple service system with two types of customers and three resources. The network is shown in Figure~\ref{fig:customer-server-small}, where resource~2 is more expensive than the other two but can serve both customer types. The penalty for lost demand is $p = (30,32)^\top$. The capacity requirement matrix $A$ has all non-zero entries equal to~1. Suppose the planning horizon is $T=2$, with deterministic demand given by $\bbD_{t=1} = (3,0)^\top$ and $\bbD_{t=2} = (0,3)^\top$. By solving Problem~\eqref{eq:fa}, we obtain the optimal capacity decision $b^* = (0,3,0)$. However, there is no constant arrival rate $\hat{\lambda}$ that yields $b^*$.

\begin{figure}[H]
\centering
\begin{tikzpicture}[
    x=1cm,y=1cm,>=stealth,
    square/.style={draw,minimum size=6mm,inner sep=0pt},
    circ/.style={draw,circle,minimum size=6mm,inner sep=0pt},
    every node/.style={font=\sffamily\footnotesize}
]

\def\ytop{0.8}    
\def\ybot{-1.2}

\node[square] (u1) at (1,\ytop) {1};
\node[square] (u2) at (4,\ytop) {2};

\node[circ] (v1) at (0,\ybot) {1};
\node[circ] (v2) at (3,\ybot) {2};
\node[circ] (v3) at (6,\ybot) {3};

\draw (u1) -- (v1);
\draw (u1) -- (v2);
\draw (u2) -- (v2);
\draw (u2) -- (v3);

\node[anchor=east] at (-0.8,\ytop)  {customer type};
\node[anchor=east] at (-0.8,\ybot)  {resource type};

\node[anchor=north] at (v1.south) {$c_1=4$};
\node[anchor=north] at (v2.south) {$c_2=6$};
\node[anchor=north] at (v3.south) {$c_3=5$};

\end{tikzpicture}
\caption{Counterexample: Non-existence of constant unbiased corrected arrival rate.}
\label{fig:customer-server-small}
\end{figure}

To see this, note that under any constant arrival rate $\hat{\lambda}$, the corresponding allocation decision $x$ must also be constant across periods. In this case, if resource~2 is used to serve type~1 customers, it is always strictly better to switch to resource~1. Likewise, if resource~2 is used to serve type~2 customers, it is always better to switch to resource~3. Hence, with a constant arrival rate, we must have $b_2=0$. This contradicts the optimal solution $b^*=(0,3,0)$. Therefore, no constant unbiased corrected arrival rate can lead to $b^*$. \hfill\Halmos
\end{example}

Example~\ref{counterexample} shows that, when searching for the decision-corrected arrival rate, we should not, in general, restrict attention to constant arrival rates. The next example further demonstrates that even a non-constant (general) decision-corrected arrival rate may fail to exist.

\begin{example}[Non-existence of the decision-corrected arrival rate]\label{ex:t1-stochastic}
Consider the same network and costs as in Figure~\ref{fig:customer-server-small}. Let \(T=1\) and suppose the demand is stochastic with finite support:
\[
D=\begin{cases}
(3,0) & \text{with probability } \tfrac{1}{2},\\
(0,3) & \text{with probability } \tfrac{1}{2}.
\end{cases}
\]
For the expectation problem \eqref{eq:outer_problem}, the capacity \(b^*=(0,3,0)\) is optimal: purchasing three units on resource pool~2 allows all demand to be served in either scenario using the flexible links (class~1 via pool~2 in the first scenario; class~2 via pool~2 in the second), yielding zero lost sale cost and procurement cost \(3c_2 = 18\). Any plan that avoids pool~2 must purchase both pool~1 and pool~3 at level~3 to cover the two scenarios, costing \(3c_1+3c_3=27>18\).

However, there does not exist a point estimator \(\hat\lambda\) for the fluid problem \eqref{eq:fa} that yields \(b^*\). Indeed, for any arrival rate, using pool~2 to serve class~1 (respectively, class~2) is strictly dominated by using the cheaper dedicated pool~1 (respectively, pool~3). Therefore, under any input \(\widehat{\lambda}_1\), every optimal solution to \eqref{eq:fa} satisfies \(b_2=0\), which contradicts \(b^*=(0,3,0)\). Hence, when \(T=1\) and demand is stochastic as above, no (time-varying or constant, here equivalent to \(T=1\)) decision-corrected arrival rate exists. \hfill\Halmos
\end{example}

Example~\ref{ex:t1-stochastic} also violates the condition for single-scenario optimality in the general stochastic programming framework of \cite{homem2024forecasting}. As a result, their analysis does not apply to the multi-resource, multi-customer newsvendor formulation. To the best of our knowledge, this is the first work to examine whether a decision-corrected arrival rate exists for this newsvendor network problem. 
Examples~\ref{counterexample} and~\ref{ex:t1-stochastic} show that one should not, in general, expect the decision-corrected arrival rate to exist. Therefore, we next examine the conditions under which a decision-corrected arrival rate does exist. To formally establish the necessary and sufficient conditions, we first present the KKT conditions for optimality.

\begin{lemma}[KKT conditions for the two-stage fluid model]\label{lem:fluid-kkt}
Fix an arrival profile $\bblambda = (\lambda_1,\dots,\lambda_T) \in (\mathbb{R}^n)^T$. For the fluid model \eqref{eq:fa}, there exist dual variables $y_t\in\mathbb{R}^m_+$ (for $A x_t\le b$) and $z_t\in\mathbb{R}^n_+$ (for $R x_t\le \lambda_t$) such that a tuple $(b,\{x_t\}_{t=1}^T;\{y_t,z_t\}_{t=1}^T)$ is optimal if and only if:

\noindent\textbf{(P) Primal feasibility:}
\[
b \ge 0,\quad x_t \ge 0,\quad A x_t \le b,\quad R x_t \le \lambda_t,\qquad \forall t= 1,\cdots,T.
\]

\noindent\textbf{(D) Dual feasibility:}
\[
\sum_{t=1}^T y_t \;\le\; c ~ T,\qquad A^\top y_t + R^\top z_t \;\ge\; R^\top p, \quad y_t\ge 0, \quad z_t \ge 0,  \qquad \forall t = 1,\cdots,T.
\]

\noindent\textbf{(CS) Complementary slackness:}
\[
\begin{aligned}
& y_t (A x_t - b) = 0, \qquad &&\forall t= 1,\cdots,T,\\
& z_t (R x_t - \lambda_t)=0, \qquad &&\forall t= 1,\cdots,T,\\
& b \Big(c ~ T - \sum_{t=1}^T y_t\Big) =0, \\
& x_t \big(A^\top y_t + R^\top z_t - R^\top p\big) =0, \qquad &&\forall t = 1,\cdots,T.
\end{aligned}
\]
\end{lemma}

\subsection{Universal existence of decision-focused fluid approximation}\label{subsec:universal-point}

In this section, we investigate the conditions on the network structure specified by parameters \((c, p, A, R)\) under which decision-corrected arrival rates exist universally for any demand distribution \(\mathcal{D}\). 
The high-level idea is to check the feasibility of the KKT conditions for Problem \eqref{eq:outer_problem}. However, due to the randomness of constraints, we first transform Problem \eqref{eq:outer_problem} to a deterministic expanded counterpart (Lemma~\ref{lem:scenario-expansion-eqprob}). 
Building on this equivalence, we establish a bridge (Theorem~\ref{thm:existence-iff-constant}) showing that universal existence of decision-corrected arrival rates is equivalent to universal existence of constant corrected arrival rates. Lastly, for the deterministic expanded formulation, we give the necessary and sufficient condition under which such a constant correction exists (Theorem~\ref{thm:const-characterization}).

When deriving a deterministic expanded formulation of the original Problem~\eqref{eq:outer_problem}, for the purpose of illustration, in this section, we consider the case where the demand distribution has finite support. Specifically, there exist scenarios (joint demand paths) \(z=1,\dots,Z\) with rational probabilities \(p_z>0\) (and \(\sum_{z=1}^Z p_z=1\)) and realizations \(D_{1:T}^{(z)}\in\mathbb{R}_+^{n\times T}\) such that
$\mathbb{P}\{\bbD = D_{1:T}^{(z)}\}=p_z$ for each $z=1,\dots,Z$.
This finite support of demand can be easily relaxed to the continuous support of demand, which is detailed in Appendix \ref{append:equivalence}. 

Given the demand support $\{D_{1:T}^{(z)}\}_{z=1}^Z$, we transform it to a scenario-expanded sequence using the following approach: We first duplicate each support path $p_z$ times so that the expanded set of support $\{D_{1:T}^{(\tilde z)}\}_{\tilde z=1}^{\tilde Z}$ has an equal weight. Next, we generate an expanded deterministic demand sequence by treating $(t,\tilde z)$ as a single
``period'' and appending all the periods into a $T\times \tilde Z$ vector as $\{D_\tau\}_{\tau=1}^{T\tilde Z}$, where $D_\tau \in \bbR_+^n$. 
Lemma~\ref{lem:scenario-expansion-eqprob} shows that this expanded sequence yields a deterministic problem that is equivalent to the original problem \eqref{eq:outer_problem}.

\begin{lemma}[Scenario–expansion equivalence]\label{lem:scenario-expansion-eqprob}
Given the demand support $\{D_{1:T}^{(z)}\}_{z=1}^Z$ and its scenario-expanded sequence $\{D_\tau\}_{\tau=1}^{T\tilde Z}$, Problem \eqref{eq:outer_problem}
is equivalent to the deterministic scenario–expanded problem
\begin{equation}\label{eq:equal-weight-expansion}
\min_{b\ge0}\; c^\top b ~ T+\sum_{z=1}^Z p_z\,\pi\!\big(b,D_{1:T}^{(z)}\big)
\;=\;
T \left\{\min_{b\ge0}\; c^\top b+\frac{1}{\tilde Z T}\pi\!\big(b,\{D_\tau\}_{\tau=1}^{T\tilde Z}\big)\right\}.
\end{equation}
In other words, Problem \eqref{eq:equal-weight-expansion} has the same objective function and set of minimizers as Problem \eqref{eq:outer_problem}.
\end{lemma}

Lemma \ref{lem:scenario-expansion-eqprob} indicates that, to determine the best capacity level, \eqref{eq:outer_problem} is equivalent to minimizing the per-period average procurement cost plus the per-period average lost demand cost. This transformation in Lemma \ref{lem:scenario-expansion-eqprob} from a stochastic setting to a deterministic expansion enables us to establish the following equivalence for the universal existence in Theorem \ref{thm:existence-iff-constant}.

\begin{theorem} {\bfseries (Universal existence $\Longleftrightarrow$ universal existence for the constant arrival for the deterministic expansion)}
\label{thm:existence-iff-constant}
Given the system parameters $ A, c, p,$ and $R$, the following statements are equivalent:
\begin{enumerate}[label=\textbf{(\Alph*)}]
    \item For any $T\ge 1$ and any joint demand distribution $\mathcal{D}$, there exists a decision-corrected arrival rate $\{\widehat\lambda_t\}_{t=1}^T$ for Problem \eqref{eq:outer_problem}. 
    \item For any $T\ge 1$, any positive integer $\tilde Z$ and any deterministic demand sequence $\{D_\tau\}_{\tau=1}^{T\tilde Z}$, there exists a {\bfseries constant} decision-corrected arrival rate $\widehat\bblambda\in(\bbR_+^n)^{T\tilde{Z}}$, where $\widehat\lambda_1=\cdots= \widehat\lambda_{T\tilde{Z}} \in \bbR^n_+$, such that $b^*(\widehat\bblambda) \in \arg\min_{b\ge 0}\Big\{c^\top b+\frac{1}{\tilde Z T}\pi\big(b,\{D_\tau\}_{\tau=1}^{T\tilde Z}\big)\Big\}$. 
\end{enumerate}
\end{theorem}

Theorem \ref{thm:existence-iff-constant} shows that to study the condition when the decision-corrected arrival rate universally exists, it is equivalent to studying the condition when the constant decision-corrected arrival rate universally exists. This transformation to a constant arrival rate enables us to obtain some interpretable insights for the conditions in Theorem \ref{thm:const-characterization}.

\begin{theorem}[Necessary and Sufficient condition for the universal existence]\label{thm:const-characterization}
Given the system parameters $c, p, A,$ and $R$, the following statements are equivalent.
\begin{enumerate}[label=\textbf{(\Alph*)}]
  \item For any integer $T\geq1$, for every demand outcome \(\bbD\in (\mathbb{R}^n_+)^T\) and every
        optimal capacity vector \(b^*(\bbD)\), there exists a
        \emph{constant} per-period arrival vector
        \(\widehat\bblambda\in(\mathbb{R}^n_+)^{T}\) (the same each period) such that
        \(b^*(\widehat\bblambda)\) is optimal for the constant-profile fluid program
        \begin{equation}\label{eq:det-const-primal}
           \left\{\min_{\substack{b\ge0,\,x\ge0}} \;
          c^\top b T \;+\; \, T p^\top\!\big(\widehat\lambda - R x\big)\right\}
          \quad\text{s.t.}\quad
          A x \le b,\;\; R x \le \widehat\lambda.
        \end{equation}
  \item For every resource pool \(h\) that is positively procured in an
        optimal solution for  \emph{at least} one demand outcome
        \(\bbD\) (i.e., there exists \(\bbD\) such that
        \(b^*_h(\bbD)>0\)), this resource is also the most cost–effective option for at
        least one customer class it can serve, i.e., there exists a column \(j\) with
        \(h(j)=h\) such that
        \[
          (A^\top c)_j \;=\; \min_{k:\, i(k)=i(j)} (A^\top c)_k
          \qquad\text{and}\qquad
          (A^\top c)_j \;\le\; \,p_{i(j)}.
        \]
\end{enumerate}
\end{theorem}

Theorem \ref{thm:const-characterization} provides an intuitive way to check whether the decision-corrected arrival rate always exists, given the system parameters. Intuitively, according to \textbf{(B)} in Theorem \ref{thm:const-characterization}, if each resource is the most cost-effective resource for at least one customer it can serve, and the procurement cost is not too large (i.e., less than the penalty for lost demand), then it guarantees the existence of a decision-corrected constant arrival rate. Note that this is both a necessary and sufficient condition for the universal existence. Combining this with Theorem \ref{thm:existence-iff-constant}, this condition is also a necessary and sufficient condition for the universal existence of the potential time-varying decision-corrected arrival rate under stochastic demand $\mathcal{D}$.

\begin{remark}[Quick check examples for Condition \textbf{(B)} in Theorem~\ref{thm:const-characterization}]
We present several typical cases that allow for a quick verification of Condition~\textbf{(B)} in Theorem~\ref{thm:const-characterization}, thereby determining whether the universal existence of the decision-corrected arrival rate holds. For the simplicity of expression, we focus on the resource that is positively procured in an optimal solution for at least one
demand outcome. 
\begin{enumerate}[label=\textbf{(\arabic*)}, leftmargin=2.2em]
\item \textbf{More resources than classes ($m>n$).}
If the number of resources $m$ is larger than the number of customer classes $n$, then the universal existence typically fails. Intuitively, we cannot expect a vector in $\bbR^n$ to determine a capacity vector in $\bbR^m$ uniformly when $n<m$.

\item \textbf{Buffer / multi-function resource.}
If a resource is not the most cost-effective option for serving \emph{any} customer class it can serve, but it is still procured as a buffer (multi-function) resource in an optimal solution, then the universal existence fails. Intuitively, the inefficiency of this buffer resource can ``hide'' demand fluctuations from being reflected in its capacity.

\item \textbf{Extremely powerful resource.}
If one resource is the most cost-effective option for \emph{all} customer classes, then the universal existence typically holds. Intuitively, in this case, all other resources will be zero.
\end{enumerate}
\end{remark}

\subsection{Distribution-dependent existence of decision-focused fluid approximation}\label{subsec:dist-dependent}

In this section, we examine the existence of a decision-corrected arrival rate from a distribution-dependent (data-dependent) perspective: Given the system parameters and the (empirical) demand distribution $\mathcal{D}$, we ask whether there exists a (possibly time-varying) arrival-rate profile that is decision-unbiased. 
 
Checking the existence is equivalent to asking whether there exists $\widehat{\lambda}$ such that $b^*(\widehat{\lambda})$ belongs to the set of minimizers of Problem~\eqref{eq:outer_problem} under $\mathcal{D}$.
To formalize which minimizers of Problem~\eqref{eq:outer_problem} under $\mathcal{D}$ admit such a corrected arrival rate, we introduce the fluid-correctable decision set $\mathcal{B}$ (Definition~\ref{def:B}). Our main
result in this section shows that $b^*$ admits a distribution-dependent decision-corrected arrival-rate profile if and only if $b^* \in \mathcal{B}$
(Theorem~\ref{thm:dist-dependent-existence}). 

\begin{definition}[Fluid-correctable decision set]\label{def:B}
Let $T \in \mathbb{N}$ and parameters $(A,R,c,p)$ be given.  
We define $\mathcal{B} \subseteq \mathbb{R}_+^m$ as the set of all $b$ for which there exist vectors $y_1,\dots,y_T \in \mathbb{R}_+^{m}$ satisfying:

\begin{enumerate}[label=(B\arabic*), ref=B\arabic*]
    \item \label{eq:B1}
    $\sum_{t=1}^T y_t \le c ~T$.

    \item \label{eq:B2}
    For every resource pool $h$ with $b_h>0$, $\sum_{t=1}^T (y_t)_h = c _h ~ T$.

    \item \label{eq:B3}
    For every time period $t$ and every resource pool $h$ with $(y_t)_h>0$, there exists a job $j$ with $h(j)=h$ such that  
    $(A^\top y_t)_j = \min_{k':\, i(k') = i(j)} (A^\top y_t)_{k'}$ and $(A^\top y_t)_j \le p_{i(j)}$.
\end{enumerate}

We call $\mathcal{B}$ the \emph{fluid-correctable decision set}.
\end{definition}

Note that the fluid-correctable decision set is independent of the distribution $\mathcal{D}$. With Definition~\ref{def:B}, verifying the existence of a decision-corrected $\bblambda$ is equivalent to checking whether $b^*(\mathcal{D})$ belongs to the set $\mathcal{B}$, as shown in the following theorem.

\begin{theorem}[Distribution-dependent existence for decision-corrected arrival rate]\label{thm:dist-dependent-existence}
Let $b^*$ be an optimal capacity decision for the expectation problem \eqref{eq:outer_problem} under distribution $\mathcal{D}$.
The following statements are equivalent.
\begin{enumerate}[label=\textbf{(\Alph*)}]
  \item There exists an arrival-rate profile
        $\widehat\lambda=(\widehat\lambda_1,\dots,\widehat\lambda_T)\in(\mathbb{R}_+^{\,n})^{T}$
        such that $b^*\in\arg\min_{b\ge0}\big\{\,c^\top T b+\pi(b,\widehat\lambda)\,\big\}$, i.e.,
        for the deterministic fluid problem
        \[
        \min_{\substack{b\ge0,\,x_t\ge0}}\ c^\top Tb+\sum_{t=1}^T p^\top(\widehat\lambda_t-Rx_t)
        \ \ \text{s.t.}\ \ A x_t\le b,\ \ R x_t\le \widehat\lambda_t\quad(\forall t=1,\dots,T),
        \]
        the vector $b^*$ is optimal.
      \item $b^*$ is within the fluid-correctable decision set $\mathcal{B}$, i.e., $b^* \in \mathcal{B}$.
    
\end{enumerate}
\end{theorem}

Theorem~\ref{thm:dist-dependent-existence} shows that, given the distribution information $\mathcal{D}$, we first compute its optimal procurement decision $b^*$ and then determine existence based on this $b^*$. To verify whether $b^*$ lies in the set $\mathcal{B}$, we show that this is equivalent to checking the feasibility of a mixed-integer optimization model in Proposition~\ref{prop:support-feasible} in Appendix~\ref{subsec:mip}. To further draw insights, we establish structural properties of the set $\mathcal{B}$ in Proposition~\ref{prop:support-B}.

\begin{proposition}[Some properties of fluid-correctable decision sets]\label{prop:support-B}
Given parameters $(T, A,R,c,p)$ and a procurement decision $b\in\mathbb{R}_+^m$, define the set
\[
H(b)\ :=\ \{\,h\in\{1,\dots,m\}:\ b_h>0\,\}.
\]
Then, for any two procurement decisions $b^0$ and $\tilde{b}$:
\begin{enumerate}
  \item[(1)] If $H(b^0)=H(\tilde b)$, then $b^0\in\mathcal{B}$ if and only if $\tilde b\in\mathcal{B}$.
  \item[(2)] If $b^{0}\in\mathcal{B}$ and $H(\tilde b)\subseteq H(b^{0})$, then $\tilde b\in\mathcal{B}$.
\end{enumerate}
\end{proposition}

In this proposition, $H(b)$ denotes the set of purchased resources. Proposition~\ref{prop:support-B}.(1) implies that for any two procurement vectors with the same procured-resource set, the existence result is identical, regardless of the specific values of $b$. Proposition~\ref{prop:support-B}.(2) further shows that if the decision-corrected arrival rate exists for a given set of procured resources, then not procuring any subset of these resources preserves existence; that is, the decision-corrected arrival rate continues to exist for the reduced procurement decision.

\section{Correction Method for Fluid Approximation}\label{sec:correction}

Section \ref{sec:existence} establishes when a decision-corrected arrival rate exists. We now turn to constructing one from the empirical demand distribution when such an arrival rate exists. 

Algorithm~\ref{alg:1} gives a data-guided construction that applies to a general service network. When the network is decomposable, so that a resource pool is specialized to a single customer class, Remark~\ref{rmk:quantile} shows that this construction reduces to a classical newsvendor quantile and can be obtained in closed form without solving for the dual variables in Algorithm~\ref{alg:1}. The derivation and the extensions to multiple resources per class and to a shared resource serving several classes are deferred to Appendix~\ref{append:quantile}.

Algorithm~\ref{alg:1} outlines the steps for checking and constructing the decision-corrected arrival rate given the observed demand profile. At the beginning of Algorithm~\ref{alg:1}, we solve for the optimal procurement decision under the empirical demand distribution. This optimization problem
$\min_{b \ge 0} \{c^\top T b + \frac{1}{Z}\sum_{z=1}^Z \pi(b,\bbD^{(z)})\}$
can be solved efficiently because the total cost is convex in the procurement decision $b$ (see Proposition~1 in \cite{harrison2005method}). Step~1 and Step~2 then verify the existence of an unbiased correction. If such a correction exists, Step~3 provides the procedure for constructing it based on $b^*$.

\begin{algorithm}[p]
\caption{Data-guided Computation of the Decision-Corrected Arrival Rate}
\begin{algorithmic}[1]
\State \textbf{Input:} Length of planning period $T\in\mathbb{N}$; System parameters $(A,R,c,p)$. Observed $Z$ demand outcomes $\{\bbD^{(1)}, \bbD^{(2)} ,\cdots,\bbD^{(Z)}\}$ across $T$ periods for all types of customers. A predetermined tie-breaking rule $\Gamma$.

\State Solve the following problem to find the optimal procurement decision $b^*$ under the collected demand profile: 
\[\min_{b \ge 0}  \quad c^\top Tb + \frac{1}{Z}\sum_{z = 1}^Z \left[ \pi(b, \bbD^{(z)}) \right].\]\label{eq:stochastic-expectation}
\State \textbf{Step 1: Check universal existence of unbiased correction.}
\If { system parameters $(A,R,c,p)$ satisfy the conditions in Theorem \ref{thm:const-characterization}}
\State For each pool \(h\) with \(b_h^*>0\), use tie-breaking rule $\Gamma$ to select a class \(i(h)\) and a column \(j(h)\)
that: (i) uses pool \(h\), (ii) attains the classwise minimum
\(\min_{k:\,i(k)=i(h)}(A^\top c)_k\), and (iii) satisfies the penalty cap
\((A^\top c)_{j(h)}\le \,p_{i(h)}\).
\State Set \(x_{j(h)}\gets b_h^*/A_{h,\,j(h)}\); Set other components of \(x\) to zero. 
\State \(\widehat\lambda \gets R x\).
\State Set \(\widehat\lambda_t \gets \widehat\lambda\) for all \(t = 1, \dots, T\).
\State \Return $\ \text{the constant decision corrected arrival rate } \{\widehat\lambda_t\}_{t=1}^T$.
\Else
    \State  Go to Step 2.
\EndIf

\State \textbf{Step 2: Check distribution-dependent existence of unbiased correction.}
\If { $b^*$ is within the fluid-correctable decision set $\mathcal{B}$}
\State Go to Step 3.
\Else
    \State \Return Decision-corrected arrival rate does not exist.
\EndIf

\State \textbf{Step 3: Calculate the decision-corrected demand arrival rate.}

\State Let $\{y_t\}_{t=1}^T$ be the vectors satisfying all conditions in Definition \ref{def:B}, selected according to the tie-breaking rule $\Gamma$.
\State
For each $t$ and each pool $h$ with $(y_t)_h>0$,
pick $j_t(h)$ with $h(j_t(h))=h$ and $(A^\top y_t)_j =\min_{k:\, i(k)=i} (A^\top y_t)_k$. Set $i_t(h)\gets i(j_t(h))$, and set
$x_{t,\,j_t(h)} \gets\ b^*_h / A_{h,\,j_t(h)}$,
with all other components of $x_t$ zero.
\State Set \(\widehat\lambda_t \gets R x_t\) for all \(t = 1, \dots, T\).
\State \Return $\{\widehat\lambda_t\}_{t=1}^T$.

\end{algorithmic}
\label{alg:1}
\end{algorithm}

In Step~3 of Algorithm~\ref{alg:1}, given the optimal procurement decision $b^*$, we construct the corresponding arrival rate as follows. We first form the auxiliary dual variables $y_t$, as specified in the definition of the set $\mathcal{B}$. Using these variables, we then construct the allocation decisions $x_t$, from which we obtain the decision-corrected arrival rate $\lambda_t$.
Theorem \ref{thm:alg} below further shows that the output of Algorithm \ref{alg:1} converges to the unbiased corrected arrival rate as the sample size goes to positive infinity.

When multiple decision-corrected arrival rates exist, a predetermined tie-breaking rule is used. Such a rule provides flexibility in tailoring the arrival rate for different purposes. Examples include penalizing changes in the arrival rate to improve smoothness, penalizing high values on specific days to enhance interpretability, or selecting the arrival-rate profile that is closest to the mean. The choice of the tie-breaking rule only affects the forecasting stability, not in-sample optimality, since every such profile is decision-corrected by construction. See the numerical section for a detailed discussion.

\begin{theorem}[Consistency of Algorithm~\ref{alg:1}]\label{thm:alg}
Suppose that the decision-corrected arrival rate exists and the demand is at most $\bar D$ for all periods and customer classes. 
Assume that the predetermined tie-breaking rule $\Gamma$ in Algorithm~\ref{alg:1} is a bounded measurable selection: there exists a compact set $\mathscr K_\lambda\subseteq(\mathbb R^n_+)^T$ such that every arrival-rate profile returned by Algorithm~\ref{alg:1} belongs to $\mathscr K_\lambda$.
Let $\widehat \lambda^{(Z)}$ be the output of Algorithm~\ref{alg:1} with $Z$ i.i.d. observations of demand profile sampled from $\mathcal{D}$.
Define $\Lambda^\star$ as the set of decision-corrected arrival rates with respect to the underlying distribution $\mathcal{D}$.
Then, under the fixed tie-breaking rule $\Gamma$, as $Z\to\infty$, we have
$
\inf_{\tilde{\lambda}\in\Lambda^\star}
\|\widehat \lambda^{(Z)}-\tilde{\lambda}\|_2
\to 0$ almost surely.
\end{theorem}

\begin{remark}[Quantile-guided correction in decomposable networks]\label{rmk:quantile}

For certain network structures, the decision-corrected arrival rate admits a far simpler description than the general construction of Algorithm~\ref{alg:1}. A network is decomposable when it splits into disjoint components that share no resources, so the first-stage objective separates and each component can be solved in isolation; the simplest such component is a single resource pool dedicated to one customer class. When a resource pool $h$ is specialized, dedicated to a single customer class $i$ through a job $j$ with capacity requirement $A_{hj} > 0$, its procurement decision does not interact with the rest of the network, so the component can be solved in isolation and reduces to a classical newsvendor problem. Writing $F_{{mix},i}(x) := \tfrac{1}{T}\sum_{t=1}^{T} F_{it}(x)$ for the pooled (time-averaged) CDF of demand for class $i$, the optimal capacity and decision-corrected rate are
\begin{equation*}
b^*_h = A_{hj}\, F^{-1}_{{mix},i}\!\Big(1 - \tfrac{c_h A_{hj}}{p_i}\Big),
\qquad
\hat\lambda_{it} = F^{-1}_{{mix},i}\!\Big(1 - \tfrac{c_h A_{hj}}{p_i}\Big),
\end{equation*}
whenever $c_h A_{hj} < p_i$ (and $b^*_h = 0$ otherwise); in practice $F_{{mix},i}$ is replaced by its empirical counterpart. A more detailed derivation and the case of multiple resources serving a single class are given in Appendix~\ref{append:quantile}.
\end{remark}

\section{Numerical Experiments}\label{sec:numerical}

In this section, we conduct numerical experiments on semi-synthetic demand data calibrated to real patient arrivals to compare the performance of our proposed decision-corrected method with three benchmark approaches: the traditional fluid approximation using forecasted mean arrival rates, a stationary SAA that directly optimizes the capacity decision from the pooled historical data, and a nonstationary SAA that leverages time-series forecasting. We aim to demonstrate the potential benefits of the decision-corrected method in a setting designed to mimic service system characteristics like daily patterns and longer-term trends.

\subsection{Experimental setup}

\begin{figure}[H]
    \centering
    \includegraphics[width=\linewidth]{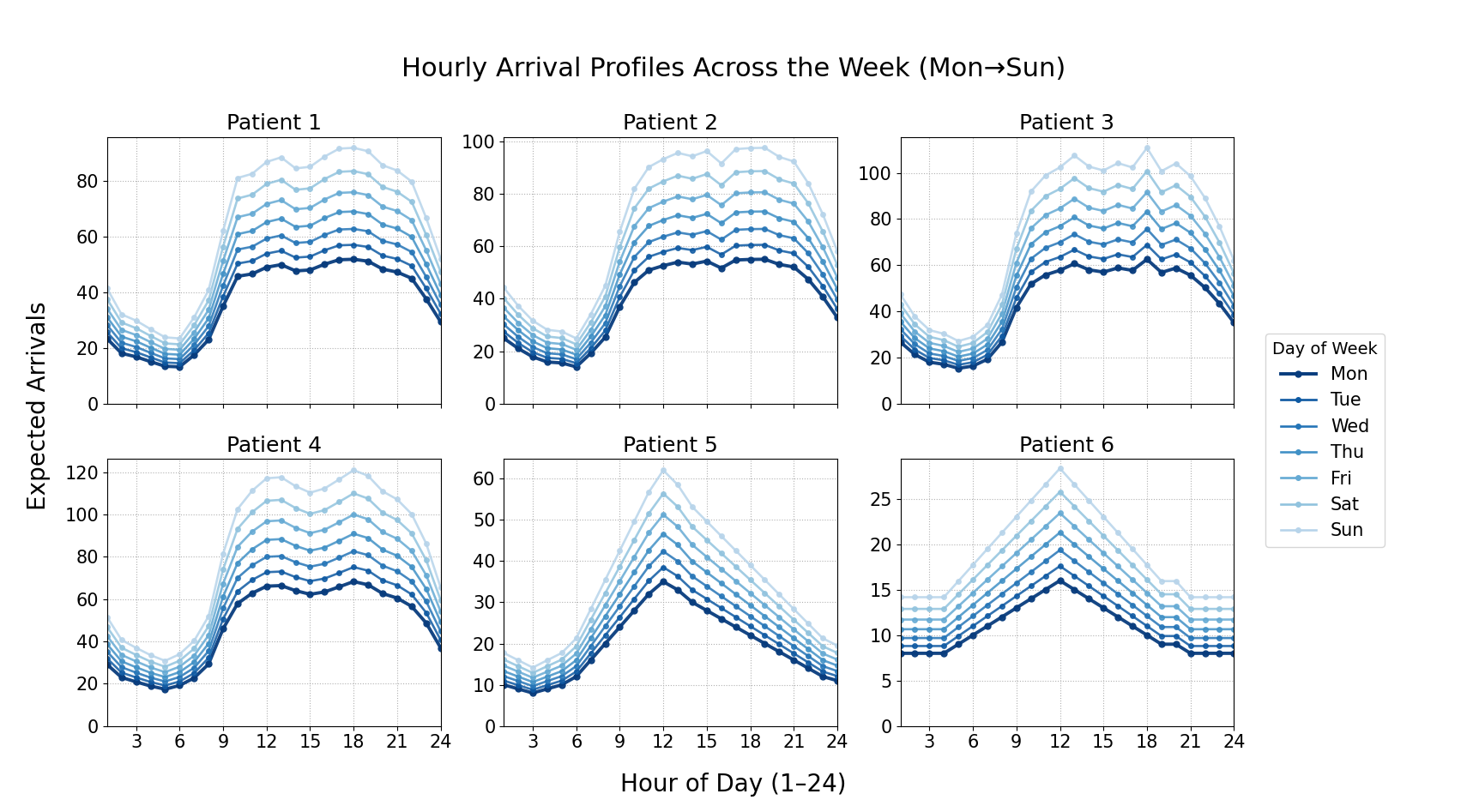}
    \caption{Average Hourly Arrival Profiles in the Training Set}
    \label{fig:arrival_real}
\end{figure}

We use hourly arrival data from four years (2014--2017), sourced from an emergency department in a hospital in Iowa, US \citep{ChoudhuryUrena2020ForecastingHourlyED}.  Based on this dataset, we construct six different patient arrival profiles. The planning period is 24 hours, and each training observation corresponds to a full 7-day week (Monday through Sunday). Let $\lambda_{d,i,t}$ denote the Poisson arrival rate for patient type~$i$ on day~$d$ of the week and hour~$t$. For $i=1,2,3,4$, we use the average hourly arrival rates across all days in 2014, 2015, 2016, and 2017, respectively. For $i=5$ and $6$, we design synthetic time-varying arrival rates that peak around noon. The values of $\lambda_{d,i,t}$ across days, hours, and patient types are shown in Figure~\ref{fig:arrival_real}. As illustrated in the figure, there is a clear increasing trend from Monday to Sunday, indicating that predictions for the following Monday should follow this weekly pattern.

\begin{figure}[H]
\centering
\begin{tikzpicture}[
    x=1cm,y=1cm,>=stealth,
    square/.style={draw,minimum size=6mm,inner sep=0pt},
    circ/.style={draw,circle,minimum size=6mm,inner sep=0pt},
    every node/.style={font=\sffamily\footnotesize}
]

\def\ytop{0.8}    
\def\ybot{-1.2}

\node[square] (u1) at (1,\ytop) {1};
\node[square] (u2) at (3,\ytop) {2};
\node[square] (u3) at (5,\ytop) {3};
\node[square] (u4) at (7,\ytop) {4};
\node[square] (u5) at (9,\ytop) {5};
\node[square] (u6) at (11,\ytop) {6};

\node[circ] (v1) at (1,\ybot) {1};
\node[circ] (v2) at (2,\ybot) {2};
\node[circ] (v3) at (4,\ybot) {3};
\node[circ] (v4) at (6,\ybot) {4};
\node[circ] (v5) at (8,\ybot) {5};
\node[circ] (v6) at (10,\ybot) {6};
\node[circ] (v7) at (12,\ybot) {7};

\draw (u1) -- (v1);
\draw (u2) -- (v2);
\draw (u2) -- (v3);
\draw (u3) -- (v4);
\draw (u4) -- (v4);
\draw (u5) -- (v5);
\draw (u5) -- (v6);
\draw (u6) -- (v6);
\draw (u6) -- (v7);

\node[anchor=east] at (-0.5,\ytop)  {patient type};
\node[anchor=east] at (-0.5,\ybot)  {resource type};

\end{tikzpicture}
\caption{The service network considered in the numerical experiments. }
\label{fig:dec_network}
\end{figure}

We consider a service network composed of $n=6$ patient classes and $m=7$ resource pools. To include various configurations, such as dedicated resources for a single patient type, shared resources, and resources handling multiple patient types, we use the service network shown in Figure \ref{fig:dec_network}, resulting in $k=9$ defined jobs. The universal existence does not hold for this network as it's an extension of the network from Example~\ref{counterexample}. The planning horizon is one day, divided into $T=24$ hourly periods.  The routing matrix $R$ specifies which job serves which patient class. We assume each job consumes one unit of capacity from its designated resource. Thus the routing, and capacity requirement matrices $R$ and $A$, as well as the per-unit per-day procurement costs $c$, and per-unit per-hour penalty cost for lost demand $p$ are defined as:

{\small 
\[ 
R =
\begin{pmatrix}
1 & 0 & 0 & 0 & 0 & 0 & 0 & 0 & 0 \\
0 & 1 & 1 & 0 & 0 & 0 & 0 & 0 & 0 \\
0 & 0 & 0 & 1 & 0 & 0 & 0 & 0 & 0 \\
0 & 0 & 0 & 0 & 1 & 0 & 0 & 0 & 0 \\
0 & 0 & 0 & 0 & 0 & 1 & 1 & 0 & 0 \\
0 & 0 & 0 & 0 & 0 & 0 & 0 & 1 & 1
\end{pmatrix},
\qquad \qquad 
A =
\begin{pmatrix}
1 & 0 & 0 & 0 & 0 & 0 & 0 & 0 & 0 \\
0 & 1 & 0 & 0 & 0 & 0 & 0 & 0 & 0 \\
0 & 0 & 1 & 0 & 0 & 0 & 0 & 0 & 0 \\
0 & 0 & 0 & 1 & 1 & 0 & 0 & 0 & 0 \\
0 & 0 & 0 & 0 & 0 & 1 & 0 & 0 & 0 \\
0 & 0 & 0 & 0 & 0 & 0 & 1 & 1 & 0 \\
0 & 0 & 0 & 0 & 0 & 0 & 0 & 0 & 1
\end{pmatrix},
\qquad \qquad
cT = 
\begin{pmatrix}
50 \\
60 \\
40 \\
70 \\
60 \\
80 \\
70
\end{pmatrix},
\qquad \qquad
p =
\begin{pmatrix}
140 \\
135 \\
130 \\
120 \\
150 \\
175
\end{pmatrix}.
\]
}

We compare four approaches for determining the first-stage capacity vector $b^*$ for the next planning period (next Monday). The first is the mean plug-in method. This approach first constructs an estimated 7-day arrival rate profile by averaging the training data. For a training set of size $N$ weeks, let $D^{(w)}_{d,i,t}$ be the observed arrivals in week $w$ on day $d$ for patient $i$ at hour $t$. The benchmark's 7-day profile is computed as the element-wise average across the $N$ sample weeks:
$$ \lambda^{\text{avg}}_{d,i,t} = \frac{1}{N} \sum_{w=1}^{N} D^{(w)}_{d,i,t} \quad \text{for } d=0..6, i=1..6, t=1..24 .$$
This results in a single average week profile, $\lambda^{\text{avg}}_{\text{week}}$, with shape $(7, 6, 24)$. Next, this 7-day profile is used to forecast the arrival rates for the next day (Monday) using a SARIMA time series forecasting model. For each patient type $i$, the 168 average hourly rates ($\lambda^{\text{avg}}_{d,i,t}$ for $d=0..6, t=1..24$) are treated as a single time series. The fitted model is then used to predict the next 24 hourly rates, yielding the Monday forecast $\lambda^{\text{forecast}}_{\text{bench}}$. Finally, the deterministic fluid model \eqref{eq:fa} is solved using $\lambda^{\text{forecast}}_{\text{bench}}$ as the input arrival rate profile $\lambda_t$ to obtain the capacity vector $b^*_{\text{bench}}$.

The second approach is the stationary SAA method. Unlike the mean-based benchmark, which constructs an arrival-rate profile and feeds it into the fluid model, the SAA benchmark bypasses the fluid approximation entirely and directly solves the two-stage stochastic optimization problem \eqref{eq:outer_problem} using all available training data. Specifically, all $N \times 7$ daily demand profiles from the training set are pooled into a single collection of $N \times 7$ demand scenarios, regardless of day of the week. The SAA problem
\[
  \min_{b \geq 0} \; c^\top Tb 
    + \frac{1}{N \times 7} \sum_{z=1}^{N \times 7} \pi(b, \bm{D}^{(z)})
\]
is then solved to obtain the capacity vector $b^*_{\text{SAA}}$, which is used directly for the next planning period without any forecasting step. This approach is a natural baseline when the demand process is approximately stationary, since in that case all historical days are exchangeable and pooling them yields the most data-efficient estimate of the optimal capacity. However, because the SAA benchmark does not produce an arrival-rate statistic, it cannot leverage time-series forecasting to account for trends or seasonality in the demand process. As a result, it is expected to perform well under stationarity but may be limited in settings with non-stationary demand patterns, such as the weekly trend present in our experimental design.

The third approach is the nonstationary SAA method. Unlike the stationary SAA, this method incorporates time-series forecasting to account for demand trends. Specifically, it uses the same SARIMA model as the mean plug-in method to forecast the expected arrival rates for the target day (Monday). Instead of feeding these mean forecasts directly into the fluid model, it generates 500 independent synthetic demand scenarios. Crucially, this step assumes that the underlying arrival process follows a known parametric distribution (in this case, Poisson) parameterized by the forecasted rates. A sample average approximation of the true two-stage stochastic problem \eqref{eq:outer_problem} is then formulated and solved over these 500 synthetic scenarios to determine the capacity vector. This allows the nonstationary SAA to capture both the weekly demand trend and the variance of the arrival process, although with a strict distributional assumption.

The fourth approach is our proposed decision-corrected method. By the decomposable-network analysis in Remark \ref{rmk:quantile}, this method determines capacity using a hybrid approach: quantile calculations for the decomposable components and the correction derived from Algorithm \ref{alg:1} for the coupled component. For the two decomposable sub-networks (resources $s_1$ serving $c_1$, and $s_3$ serving $c_2$; $b_2^*$ is set to 0 due to cost dominance), capacities $b_1^*, b_3^*$ are calculated directly using the quantile formula $b_h^* = A_{hj} \hat{F}_{mix,i}^{-1}(1 - \frac{c_h A_{hj}}{p_i})$ derived in Remark \ref{rmk:quantile}. The empirical CDF $\hat{F}_{mix,i}$ is constructed by pooling all $N \times 7 \times 24$ raw observations from the full training dataset for the relevant patients. For the two coupled sub-networks ($s_4$ serving $c_3, c_4$, and $s_5, s_6, s_7$ serving $c_5, c_6$), 7-day corrected arrival rate profiles, $\lambda^{\text{corr}}_{\text{week}}$ (both with shapes $(7, 2, 24)$), are constructed. This is done by applying the core logic of Algorithm \ref{alg:1} independently for each day of the week $d \in \{0, \dots, 6\}$ (Monday to Sunday). For a given day $d$, we first solve the SAA problem restricted to the sub-network, using all $N$ training samples of demand for patients $c_3,c_4$, and , $c_5, c_6$ on day $d$. Let the resulting optimal capacity vector for the sub-network resources be $b^*_{\text{day}}$. We then solve the MIP certificate in Appendix \ref{subsec:mip} (restricted to each sub-network), using the smoothness-penalizing objective that minimizes successive hour-to-hour differences in $y_{t,h}$ to obtain the corresponding dual vectors $y_t$ for $t=1..24$. Namely, among all dual certificates that yield a valid decision-corrected arrival rate for $b^*_{\text{day}}$, the tie-breaking rule $\Gamma$ returns the one minimizing $\sum_{h} \max_{t} |y_{t,h} - y_{t+1,h}|$, the sum over resources of the largest hour-to-hour change in the dual profile. Finally, using $b^*_{\text{day}}$ and the obtained $y_t$ vectors, we construct the corrected arrival rate profile for day $d$, $\lambda^{\text{corr}}_{\text{day}}$ (shapes $(2, 24)$), following the procedure outlined in Step 3 of Algorithm \ref{alg:1}. These 7 daily profiles ($\lambda^{\text{corr}}_{\text{day}}$ for $d=0..6$) are combined to form $\lambda^{\text{corr}}_{\text{week}}$. This 7-day corrected profile is then used to forecast the next Monday for the sub-network ($\lambda^{\text{forecast}}_{\text{sub,ours}}$) using the same time series prediction model as the benchmark. Deterministic fluid models \eqref{eq:fa}, restricted to the sub-network parameters, are solved using $\lambda^{\text{forecast}}_{\text{sub,ours}}$ to find $b_4^*$, and $b_5^*, b_6^*, b_7^*$. The final capacity $b^*_{\text{ours}}$ combines the quantile-based $b_1^*, \dots, b_3^*$ and the fluid-based $b_4^*, \dots, b_7^*$.

Finally, we describe the evaluation procedure. For each method and each training set size $N$, the resulting capacity vector $b^*$ is evaluated against the 30 samples in the test dataset. The total cost is calculated as the sum of the first-stage procurement cost ($c^\top T b^*$) and the average second-stage lost demand cost over the 30 test samples. The lost demand cost for each individual test sample is found by solving the second-stage linear program \eqref{eq:inner_problem}. The final results presented represent the average costs obtained over the 10 independent trials for each training size $N$.

\subsection{Results}

\begin{figure}[H]
    \centering
    \includegraphics[width=0.95\linewidth]{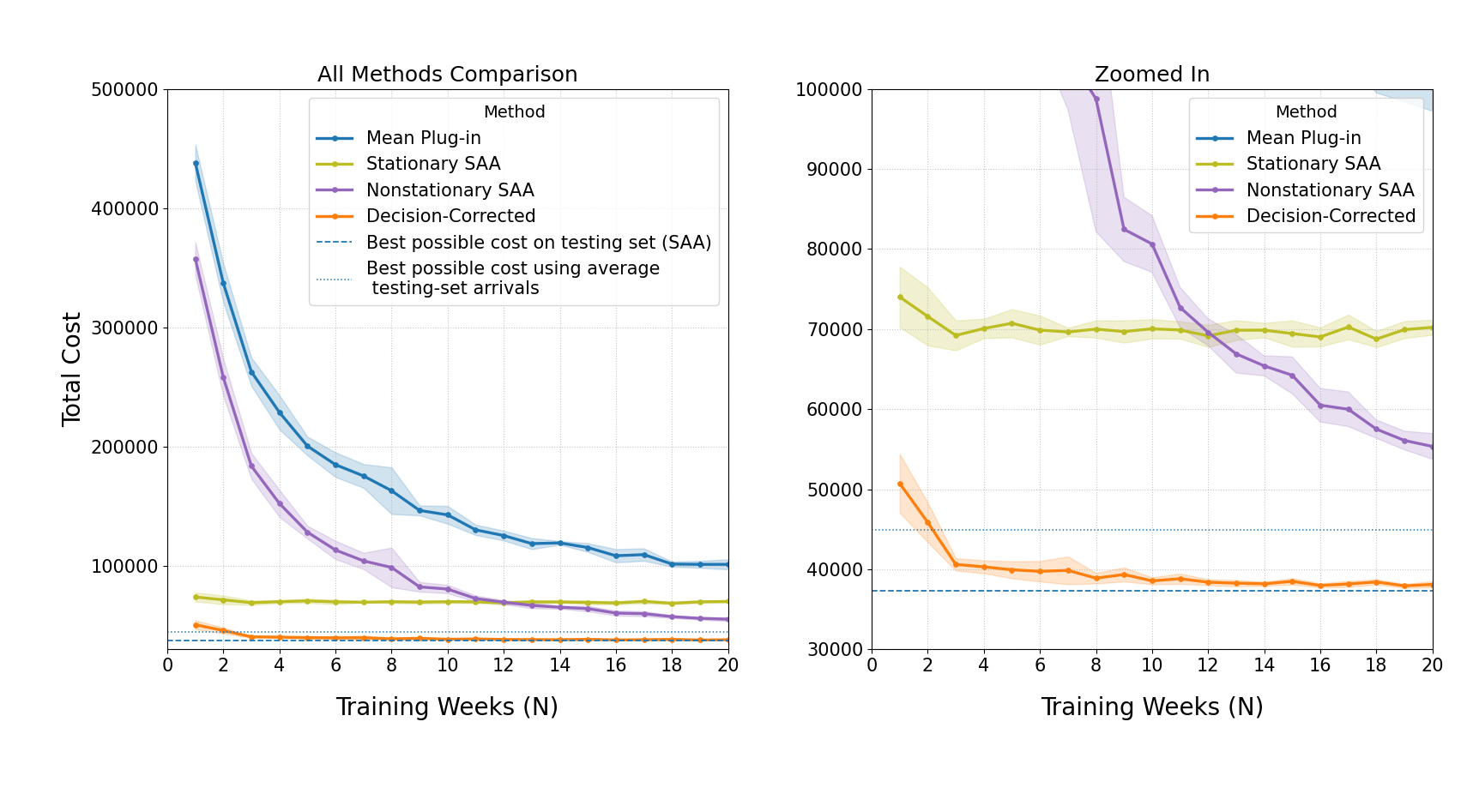}
    \caption{Average total cost on the test set as a function of training sample size for the mean plug-in benchmark, the stationary SAA benchmark, the nonstationary SAA benchmark, and the decision-corrected method.}
    \label{fig:total}
\end{figure}

Figure~\ref{fig:total} presents the average total cost during the testing period for the four approaches: the mean plug-in benchmark, the stationary SAA benchmark, the nonstationary SAA benchmark, and the decision-corrected method.

Figure~\ref{fig:total} shows that as more observations become available, all four methods achieve lower costs during the testing period. However, the traditional fluid approximation is biased, causing its total cost to converge to a higher value than that of our decision-corrected method. At our largest training size of 20 weeks, the decision-corrected approach achieves an average test cost that is approximately 62.5\% lower than that of the mean plug-in benchmark.

The stationary SAA benchmark, which directly optimizes the capacity decision from the pooled historical data, outperforms the mean-based fluid approximation. However, the decision-corrected approach still achieves an average test cost approximately 45.8\% lower than the stationary SAA benchmark at $N=20$ weeks. This gap arises because the stationary SAA benchmark pools all training days regardless of day of the week and therefore cannot capture the weekly trend present in our demand process: it optimizes capacity for an average historical day rather than for the specific future Monday being planned. In contrast, the decision-corrected method produces arrival-rates that are passed through a time-series forecasting model, enabling it to extrapolate the weekly trend to the target day when making capacity decisions.

The nonstationary SAA benchmark, described above, is the strongest of the three benchmarks. It captures both the weekly demand trend and the variance of the arrival process, although it requires a parametric Poisson assumption. Even so, the decision-corrected method achieves an average test cost approximately 31.1\% lower than the nonstationary SAA at $N=20$ weeks. Figure~\ref{fig:total} shows that the nonstationary SAA method has lower costs as training data grows and it trends toward the decision-corrected level. However, that convergence is slow, and the gap remains substantial across all training sizes considered.

The two horizontal reference lines in Figure~\ref{fig:total} are computed with full knowledge of the test set. The lower (dashed) line is the cost of solving the two-stage stochastic program~\eqref{eq:outer_problem} directly on the realized test-set arrivals, which is the lowest total cost any capacity decision can achieve in hindsight. The upper (dotted) line is the cost of plugging the test-period average arrival profile, computed as the per-hour mean across the test Mondays, into the deterministic fluid model~\eqref{eq:fa}. The gap between the two lines reflects the bias that remains when a mean profile is plugged into the fluid model.

Taken together, the comparison against the three benchmarks highlights three complementary properties of the decision-corrected method:

\begin{itemize}

    \item Relative to the mean plug-in benchmark, the decision-corrected method is unbiased. Plugging the mean arrival rate into the fluid model is biased, and the bias does not vanish as more data become available. The mean plug-in benchmark stabilizes at a strictly suboptimal cost, while the decision-corrected statistic removes this bias by construction.

    \item Relative to the stationary SAA benchmark, the decision-corrected method can account for trend and seasonality in demand. The stationary SAA pools all historical days and cannot extrapolate these patterns to the target planning day. The decision-corrected method instead produces an arrival-rate label that can be combined with time-series forecasting, so the weekly trend in our experimental design is reflected in the capacity decision.

    \item Relative to the nonstationary SAA benchmark, the decision-corrected method learns faster and does not require a distributional assumption on demand. The nonstationary SAA assumes Poisson arrivals and solves a stochastic program over many sampled scenarios, which converges slowly in the training horizon. The decision-corrected method extracts the point statistic directly from observed data and attains competitive cost even with very few training weeks.

\end{itemize}

Appendix~\ref{append:numerical_experiments} further presents a detailed breakdown of costs across the decomposable and non-decomposable components of the network.

\begin{figure}[H]
    \centering
    \includegraphics[width=0.9\linewidth]{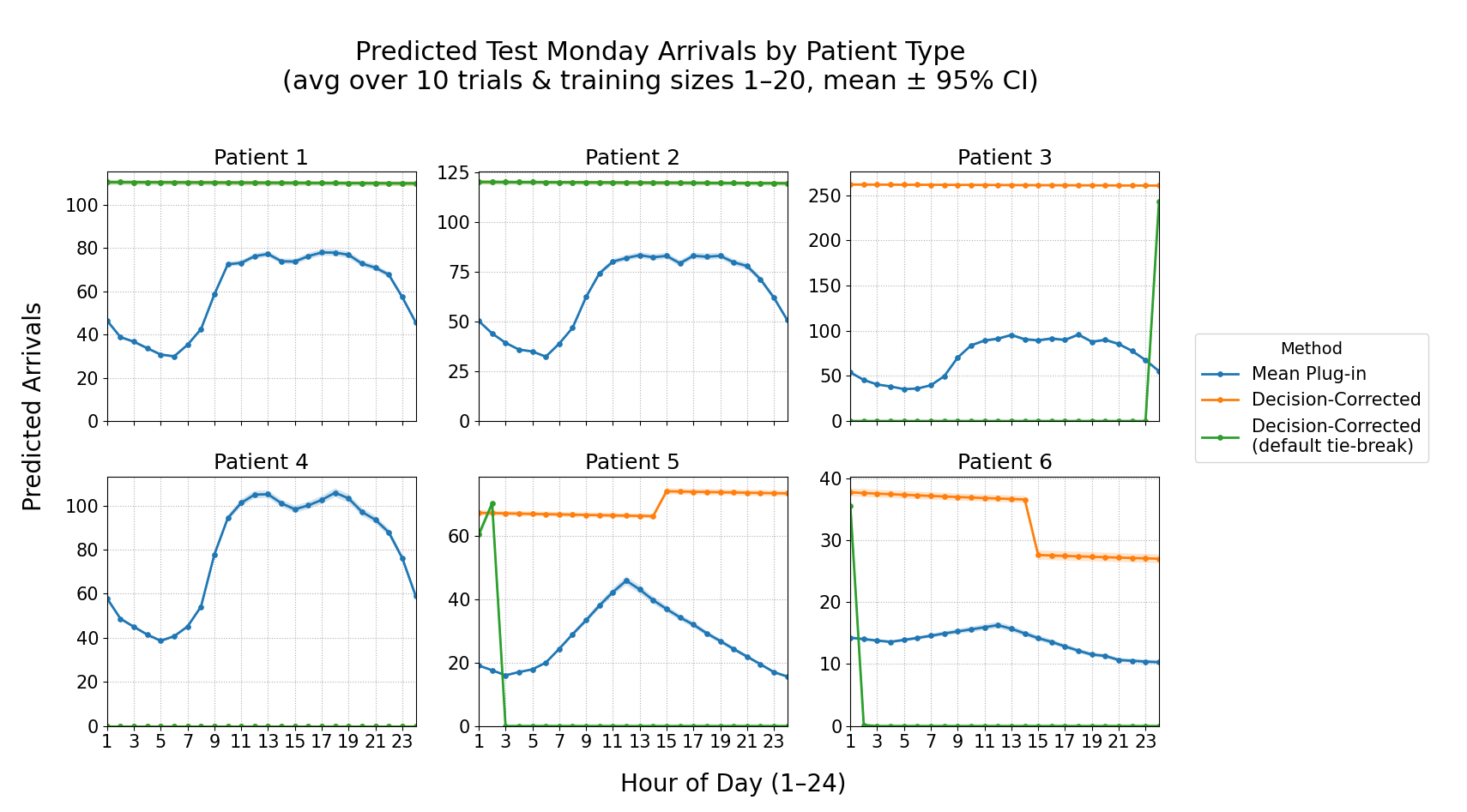}
    \caption{Predicted 24-hour arrival-rate profiles for each patient class under the mean plug-in and under the decision-corrected approach, averaged over ten trials trained on twenty weeks of data.}
    \label{fig:forecast}
\end{figure}

Figure~\ref{fig:forecast} further illustrates the predicted arrival rates~$\boldsymbol{\lambda}$ obtained from either the mean plug-in method or our decision-corrected approach over a 24-hour testing period. 
The results represent the average of ten trials, each trained on 20 weeks of data. 
The blue curve corresponds to the traditional method, which closely follows the demand trend observed in the training set and can be viewed as a conventional forecast aimed at matching the empirical mean of demand. The two SAA benchmarks do not appear in this figure: the stationary SAA produces a capacity decision directly from pooled data without constructing an arrival-rate profile, and the nonstationary SAA uses the same underlying SARIMA forecast as the mean plug-in method but samples scenarios from it rather than feeding the rate directly into the fluid LP. In contrast, the orange and green curves represent two alternative predicted arrival rates produced by the decision-corrected approach. Although these two arrival-rate profiles differ in structure and interpretability, both are decision-corrected for the training distribution, so each induces an optimal procurement decision in-sample. The green curve uses the default tie-breaking rule, where the MIP certificate carries a zero objective, and the solver returns a feasible set of dual variables $y_{t,h}$ depending on its internal rules. The orange curve uses the rule that penalizes period-to-period variability in $y_{t,h}$, yielding a smoother profile. (see Appendix~\ref{subsec:mip} for both formulations). The orange curve is the one used for the analysis in this section.

\begin{itemize}
    \item The decision-corrected arrival rate~1 (green curve) remains constant for patient types~1, 2, and~4, and exhibits a sharp jump for patient types~3, 5, and~6.
    \item The decision-corrected arrival rate~2 (orange curve) is obtained by applying tie-breaking rules that penalize excessive hour-to-hour variability. For patient types~1, 2, and~4, it matches the green curve, while for patient types~3, 5, and~6, it is noticeably smoother and closer to the real demand.
\end{itemize}

\begin{figure}[ht]
    \centering
    \includegraphics[width=0.9\linewidth]{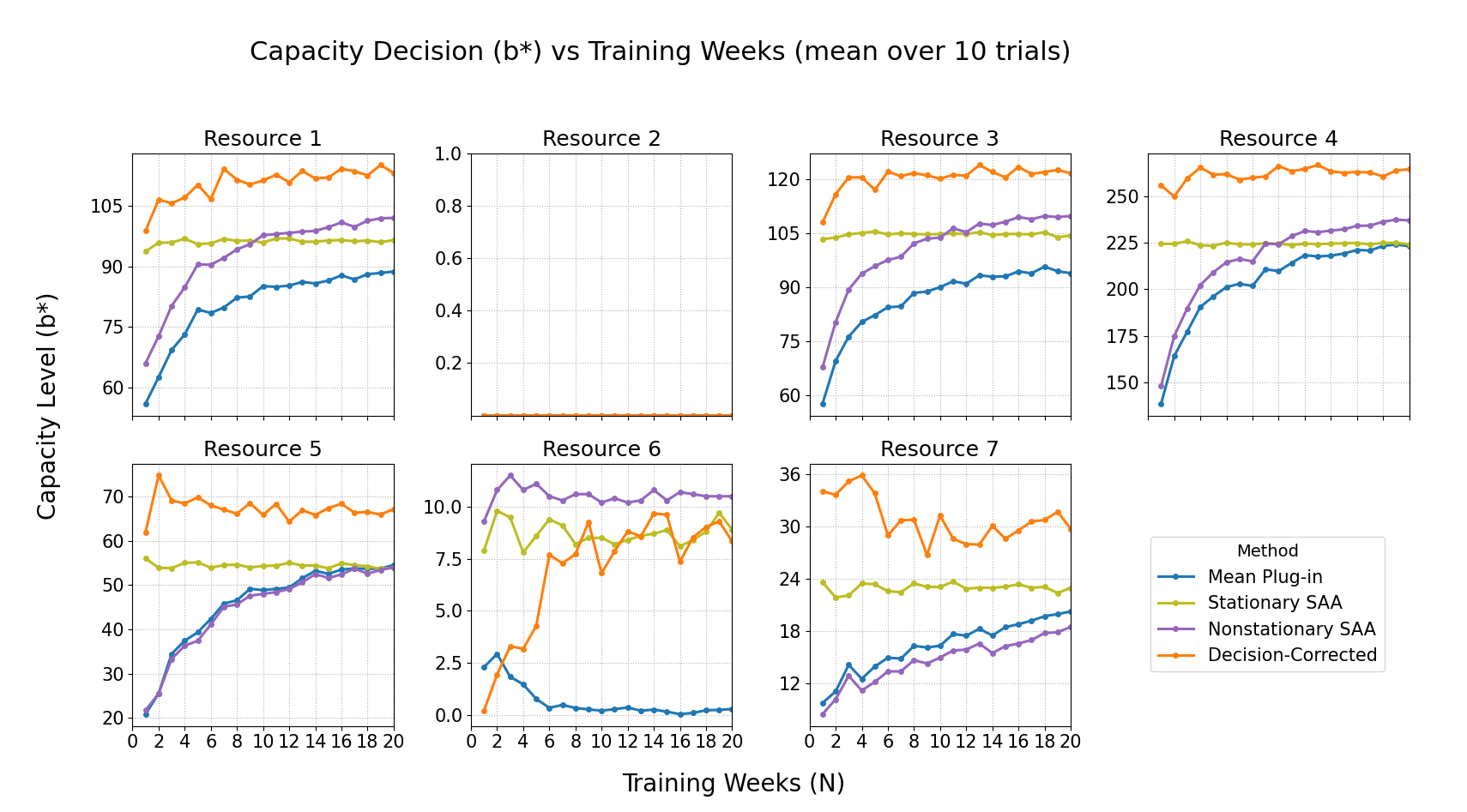}
    \caption{Average capacity levels by resource type as a function of training sample size, comparing policies derived from the mean plug-in benchmark, the stationary SAA benchmark, the nonstationary SAA benchmark, and the decision-corrected arrival rate.}
    \label{fig:bstar}
\end{figure}

Figure~\ref{fig:bstar} presents the average capacity decisions derived from the decision-corrected approach, the mean plug-in method, the stationary SAA, and the nonstationary SAA as more weeks of training data become available.

Recall that the network in our experimental setup has three independent components. The first component contains resources 1, 2, and 3, with customers 1 and 2; after pruning the dominated resource 2, this component reduces to two independent single-resource, single-customer problems, so the optimal capacity has a closed-form newsvendor quantile (see Remark~\ref{rmk:quantile}). On resources 1 and 3, the decision-corrected method computes the quantile directly and reaches its capacity level essentially from $N=1$, as shown in Figure \ref{fig:bstar}. The other three methods' capacity levels start lower and either climb toward the decision-corrected level as training data grows (mean plug-in and nonstationary SAA) or remain flat at an intermediate level (stationary SAA).

The second and third components have a richer structure. The second component is a single resource (resource 4) that serves two customer classes, 3 and 4. Because one resource is shared by two demand streams, the optimal capacity benefits from a variance-pooling effect, and procuring the sum of the two single-customer quantiles would over-procure. The decision-corrected method captures this pooling effect by construction and reaches its capacity level on resource 4 from the start. Figure \ref{fig:bstar} shows that among the other methods, nonstationary SAA approaches the decision-corrected level as $N$ grows and nearly matches it at $N=20$, while stationary SAA stays flat below it and the mean plug-in method climbs but remains the lowest.

The third component consists of resources 5, 6, and 7 and customers 5 and 6. Resource 5 serves only customer 5, resource 7 serves only customer 6, and resource 6 serves both customers. Therefore, the three capacity decisions in this component are coupled: changing the capacity of resource 6 affects how much capacity is needed on resources 5 and 7. On resource 6, the decision-corrected method and the two SAA benchmarks procure similar capacity levels, whereas the mean plug-in method collapses toward zero. On resources 5 and 7, however, the decision-corrected method procures meaningfully higher capacities than all three benchmarks.

This component further demonstrates the advantage of the decision-corrected method beyond decomposable networks. As shown in Appendix~\ref{append:numerical_experiments}, the decision-corrected method also achieves the lowest cost on the non-decomposable subsystem, where the correction is computed through Algorithm~\ref{alg:1} rather than through a closed-form quantile. Thus, its overall savings come from both decomposable and non-decomposable parts of the network. Taken together, these patterns indicate that the advantage of the decision-corrected method comes from accurately capturing the joint capacity decision across the three resources, rather than from simply inflating any single capacity level.

\section{Conclusion}\label{sec:Conclusion}

We consider a two-stage stochastic programming problem in the form of a multi-resource, multi-customer, multi-period newsvendor problem. Given an unknown demand distribution, we study whether one can identify a decision-corrected arrival rate as the input to the fluid approximation. Using such a decision-corrected arrival rate can significantly reduce estimation complexity. We show that the mean of the demand process is generally biased for decision-making, and that such a decision-corrected arrival rate may not even exist.
We establish both necessary and sufficient conditions for the existence of decision-corrected arrival rates. Numerical studies based on semi-synthetic demand data calibrated to real patient arrivals demonstrate that using decision-corrected arrival rates for demand forecasting yields lower total cost than the traditional mean-based fluid approximation, a stationary SAA benchmark, and a nonstationary SAA benchmark that incorporates time-series forecasting.

Future research directions include extending the proposed bias-correction method to more general simulation systems, developing non-asymptotic performance guarantees, and implementing the framework in facility location, healthcare, and other application domains. More broadly, our results highlight that, in capacity sizing problems, using the naive per-period average demand as the Poisson arrival rate can introduce systematic bias in fluid approximation. This observation underscores the need to align arrival-rate estimation with the downstream decision problem, rather than treating demand estimation as a purely predictive task.

\ACKNOWLEDGMENT {The authors thank Vidyadhar G. Kulkarni, Serhan Ziya, and Nilay Tanik Argon for early discussions that helped shape this work. We also thank J. Michael Harrison, Assaf Zeevi, Itai Gurvich, Michael R. Wagner, Huseyin Topaloglu, Tito Homem-de-Mello, Ali Aouad, George Shanthikumar, Jing Dong, and Ningyuan Chen for their insightful comments and feedback. We thank our colleague Xiangying Huang for discussions regarding the proof of Theorem 4.}

\bibliographystyle{informs2014} 
\bibliography{references}

\newpage
\begin{APPENDICES}
\section{Details about Numerical Experiments}\label{append:numerical_experiments}

We generate semi-synthetic demand data calibrated to real arrival rates, with daily patterns, incorporating a weekly trend. The fundamental unit of training data is a 7-day week (Monday to Sunday). Hourly base arrival rates $\lambda_{0,i,t}$ for the first day (Monday) are established using average hourly arrival data from four years (2014-2017) for patients of type 1-4, sourced from a dataset of emergency department arrivals, and specified profiles for patient types 5-6. A multiplicative day-over-day trend is applied within each 7-day week. The base rate for day $d$ ($d=0$ for Monday) is $(trend\_factor)^d$ times the base Monday rate. We set $trend\_factor = 1.1$, representing a 10\% daily increase in mean arrivals from Monday to Sunday. Let $\lambda_{d,i,t}$ denote the resulting mean rate for day $d$, patient $i$, hour $t$. The training dataset consists of $N$ independent sample weeks, where $N$ ranges from 1 to 20. For a given training set size $N$, we generate $N$ i.i.d. weeks. The demand $D_{w,d,i,t}$ for week $w \in \{1, \dots, N\}$, day $d \in \{0, \dots, 6\}$, patient $i$, and hour $t$ is drawn independently from a Poisson distribution with rate $10 \times \lambda_{d,i,t}$. The scaling factor 10 increases arrival volumes. The resulting training dataset for a specific run has shape $(N, 7, 6, 24)$. The test dataset comprises 30 independent samples of a single day's demand profile, representing 30 Mondays that follow the established trend. The mean rates for these test Mondays, $\lambda^{\text{test}}_{i,t}$, are calculated by applying the daily trend factor 7 times to the base Monday rates: $\lambda^{\text{test}}_{i,t} = \lambda_{0,i,t} \times (trend\_factor)^7$. The 30 test samples $D^{\text{test}}_{s,i,t}$ are drawn independently: $D^{\text{test}}_{s,i,t} \sim \text{Poisson}(10 \times \lambda^{\text{test}}_{i,t})$. The test dataset has shape $(30, 6, 24)$. 

For the mean plug-in benchmark, the nonstationary SAA, and the decision-corrected method, a seasonal ARIMA (SARIMA) model is fitted to this series. Specifically, we use a SARIMA model with order $(p=1, d=0, q=0)$ for the non-seasonal part and $(P=1, D=0, Q=0, s=24)$ for the seasonal part, including a constant trend term (denoted SARIMA(1,0,0)x(1,0,0,24)). This model structure assumes the current hour's rate depends on the previous hour's rate (lag 1) and the rate from the same hour on the previous day (lag 24). 

The stationary SAA benchmark does not use any time-series forecasting. Instead, for each training set of $N$ weeks, all $N \times 7$ daily demand profiles are pooled into a single collection of demand scenarios, treating each day identically regardless of its position in the week. The SAA optimization problem is then solved over this pooled collection to obtain a capacity vector $b^*_{\text{SAA}}$ which is applied directly to the test period. Because the stationary SAA benchmark does not distinguish between days of the week, it effectively optimizes for the average historical demand level rather than for the specific future day (Monday) being planned. This makes the SAA benchmark a natural comparison that isolates the value of the forecasting component.

When solving the second-stage allocation problem, the SAA benchmark optimization, and MIP certificate, all linear programs (including the SAA problems, the deterministic fluid models, and the second-stage evaluation LPs) were formulated and solved using the \texttt{cvxpy} modeling library with the ECOS solver. The mixed-integer programs for the MIP certificate were also formulated with \texttt{cvxpy} and solved using the \texttt{GLPK\_MI} optimizer. The SARIMA forecasting models were implemented using the \texttt{statsmodels} library. The entire process was repeated 10 times using independently seeded random number streams for the training data generation.

Figure~\ref{fig:subnetwork} decomposes the total cost into contributions from the decomposable sub-network (resources 1 to 3 serving patients 1 and 2) and the non-decomposable sub-network (resources 4 to 7 serving patients 3 to 6). For both subsystems, the decision-corrected method achieves the lowest cost, followed by the nonstationary SAA benchmark, the stationary SAA benchmark, and finally the mean plug-in benchmark resulting in the highest cost. The decision-corrected method maintains a consistent cost advantage regardless of the network structure. It outperforms all three benchmarks on both the simple decomposable resources and the complex shared resources.

\begin{figure}[H]
    \centering
    \includegraphics[width=0.9\linewidth]{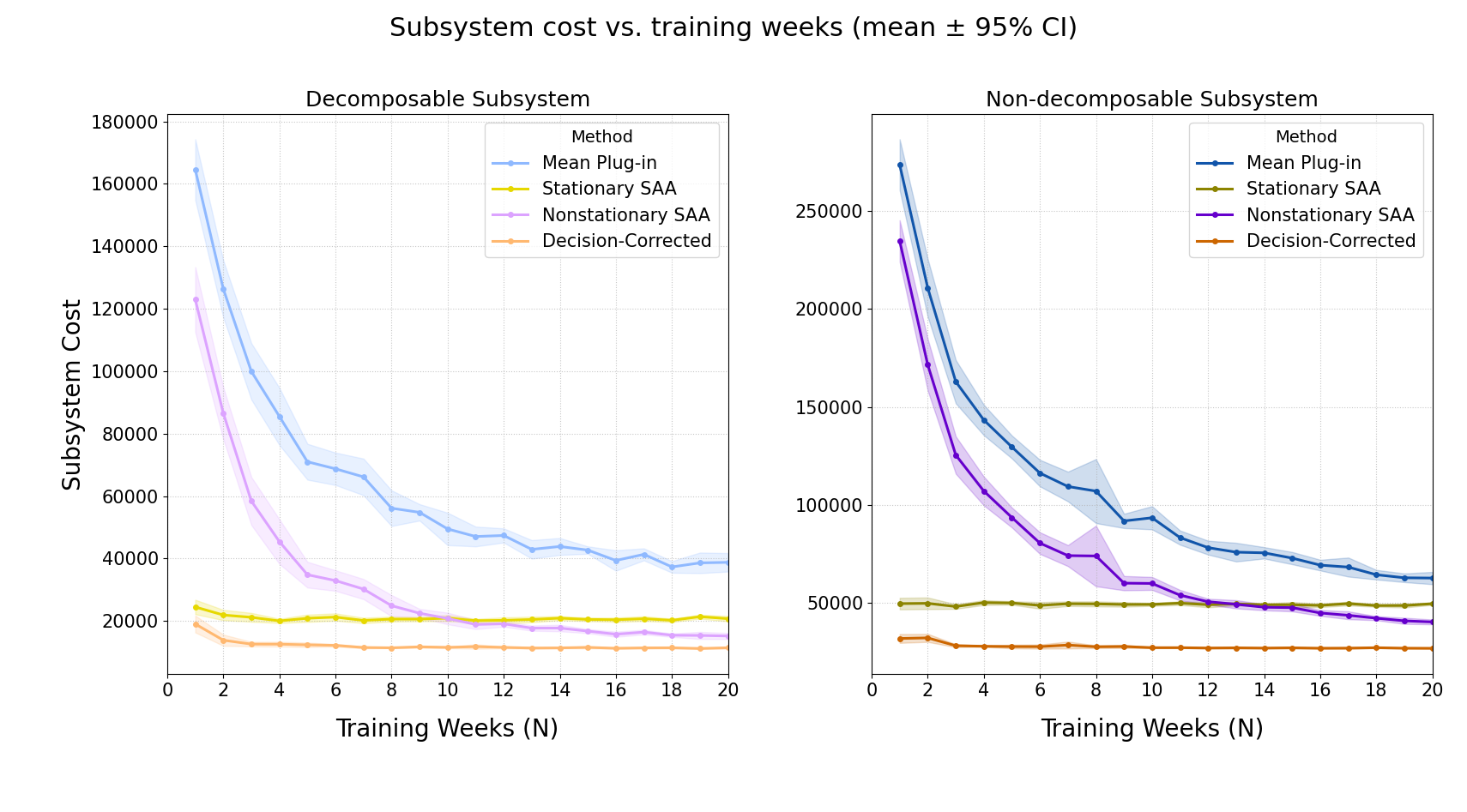}
    \caption{Subsystem-level decomposition of total cost under the mean plug-in benchmark, the stationary SAA benchmark, the nonstationary SAA benchmark, and the decision-corrected approach, separating decomposable and non-decomposable components of the network.}
    \label{fig:subnetwork}
\end{figure}

\section{Details about Quantile-guided correction in decomposable networks}\label{append:quantile}

While Algorithm~\ref{alg:1} provides a data-guided method for constructing the decision-corrected arrival rate $\hat\lambda$, it requires first solving for the optimal decision $b^*$ and then constructing the associated dual variables, both of which may be computationally intensive. For certain network structures, however, the correction admits a far simpler, closed-form description.

The quantile characterizations developed below specialize the classical newsvendor-network analysis, in \cite{mieghem2002newsvendor} and \cite{HARRISON199917}, in which the optimal capacities solve generalized critical-fractile conditions of the multivariate demand distribution. For decomposable structures these conditions reduce to ordinary newsvendor quantiles, whereas when a resource is shared across classes the pooling effect can make the optimal capacity smaller than, equal to, or larger than the capacity obtained by managing each class separately.

A network is decomposable if it can be partitioned into disjoint connected components $G = \sqcup_l G_l$, where each component $G_l$ involves its own subset of resources $H_l$, customers $I_l$, and the jobs connecting them. Because there are no links across components, the first-stage objective separates into a sum over components, and each component can be analyzed in isolation:
\begin{equation*}
\min_{b \ge 0} \; c^\top T b + \mathbb{E}[\pi(b, \bm{D})]
= \sum_l \min_{b_{H_l} \ge 0} \Big\{ c_{H_l}^\top T\, b_{H_l} + \mathbb{E}[\pi_l(b_{H_l}, D_{I_l})] \Big\}.
\end{equation*}
This separation lets us study the optimal capacity $b^*$ and look for quantile-based corrections one component at a time.

The simplest such case is a specialized resource: a single resource pool $h$ dedicated to serving a single customer class $i$ through a job $j$ with capacity requirement $A_{hj} > 0$. Since this component's procurement decision does not affect the rest of the network, it can be solved on its own, and it reduces exactly to a newsvendor problem. Writing $F_{{mix},i}(x) := \tfrac{1}{T}\sum_{t=1}^{T} F_{it}(x)$ for the pooled (time-averaged) CDF of demand for class $i$, the optimal capacity and the decision-corrected arrival rate are
\begin{equation}
b^*_h = A_{hj}\, F^{-1}_{{mix},i}\!\Big(1 - \tfrac{c_h A_{hj}}{p_i}\Big),
\qquad
\hat\lambda_{it} = F^{-1}_{{mix},i}\!\Big(1 - \tfrac{c_h A_{hj}}{p_i}\Big),
\label{eq:quantile-correction}
\end{equation}
whenever $c_h A_{hj} < p_i$ (and $b^*_h = 0$ otherwise). In practice $F_{{mix},i}$ is replaced by its empirical counterpart $\hat F_{{mix},i}$ from pooled data. The correction is thus a single newsvendor quantile, recovering the classical solution of Example~\ref{example:newsvendor} and requiring no dual construction.

When several resources can serve the same single customer class, the same quantile rule applies after pruning: all but the most cost-effective resource are dominated and procured at zero, reducing the component to the specialized-resource case above. 

Besides the above two sub-network structures with a single customer, we can also find a closed-form optimal capacity when one resource serves several customers, but this case is more complex than the single-customer quantile. For a single customer, the optimal capacity is a newsvendor quantile that depends only on the cost parameters and the demand distribution of that customer, and the decision-corrected arrival rate can be found by the empirical quantile. When one resource serves several customers, the optimal capacity instead solves a condition that links the served classes through their joint demand distribution. Even under independence, computing it requires integrating over the individual distributions, as the following remark shows, and dependence across classes makes it more difficult. In addition, even when the optimal capacity is known, the decision-corrected arrival rate is not given by a single quantile, because the capacity can be served by one class, the other, or a mix of the two. We therefore construct the corrected rate from the optimal capacity using Algorithm~\ref{alg:1} instead of writing it in closed form. The general multi-resource multi-customer system is even harder and usually has no simple quantile solution unless dominance-based pruning reduces it, in which case the procedure in Algorithm~\ref{alg:1} can again be used.

\begin{remark}[Shared resource serving two customer classes]
Consider a single resource pool, with unit procurement cost $c > 0$, serving two customer classes $i = 1, 2$ in a single period. Class $i$ has demand $D_i \sim F_i$ and a per-unit lost-demand penalty $p_i$, with $p_1 \ge p_2 > 0$. Because class~1 carries the higher penalty, it is served first; any capacity left after serving class~1 is then used for class~2.

For a capacity level $Q$, the realized lost demand is $(D_1 - Q)^+$ for class~1 and 
$\big(D_2 - (Q - D_1)^+\big)^+$ for class~2, so the expected total cost is
\[
C(Q) \;=\; cQ \;+\; p_1\,\mathbb{E}\big[(D_1 - Q)^+\big] 
\;+\; p_2\,\mathbb{E}\Big[\big(D_2 - (Q - D_1)^+\big)^+\Big].
\]
The marginal cost of one additional unit of capacity is
\[
\frac{dC}{dQ} \;=\; c \;-\; p_1\,\mathbb{P}(D_1 \ge Q) 
\;-\; p_2\,\mathbb{P}\big(D_1 < Q \le D_1 + D_2\big),
\]
so the optimal capacity $Q^*$ satisfies the first-order condition
\[
c \;=\; p_1\,\mathbb{P}(D_1 \ge Q^*) \;+\; p_2\,\mathbb{P}\big(D_1 < Q^* \le D_1 + D_2\big).
\]
When $D_1$ and $D_2$ are independent, the second probability has the closed form $\int_0^{Q^*}\big(1 - F_2(Q^* - x)\big)\,dF_1(x)$, and therefore $Q^*$ solves
\[
c \;=\; p_1\big(1 - F_1(Q^*)\big) 
\;+\; p_2\int_0^{Q^*}\big(1 - F_2(Q^* - x)\big)\,dF_1(x).
\]
This equates the unit procurement cost with the expected penalty avoided by the marginal unit: the first term is the probability that the unit prevents a lost sale from the high-penalty class~1, and the second is the probability that the unit is not needed by class~1 but prevents a lost sale from class~2. Once $Q^*$ is obtained, it plays the role of $b^*_h$ for the shared pool, and the decision-corrected arrival rate follows from Step~3 of Algorithm~\ref{alg:1}.

As shown by \cite{HARRISON199917} and \cite{mieghem2002newsvendor}, pooling demand onto a shared resource lowers the optimal cost, but the resulting capacity can be smaller than, equal to, or larger than the total capacity of dedicated resources, depending on the cost and demand parameters. In our setting this means $Q^*$ is generally not the sum $F_1^{-1}(1 - c/p_1) + F_2^{-1}(1 - c/p_2)$ of the two dedicated-resource quantiles. The following two cases show that $Q^*$ can fall on either side of this sum.  If $D_1, D_2$ are independent and uniform on $[0,1]$ with $c = 1$, $p_1 = 4$, and $p_2 = 2$, the dedicated quantiles are $Q_1^* = \frac{3}{4}$ and $Q_2^* = \frac{1}{2}$, summing to $\frac{5}{4}$, whereas $Q^* = 1 < \tfrac{5}{4}$. If instead $D_1, D_2$ are independent $\text{Beta}(1,3)$ with $c = \frac{1}{2}$ and $p_1 = p_2 = 1$, the dedicated quantiles are 
$Q_1^* = Q_2^* = 1 - 2^{-1/3} \approx 0.206$, summing to $2 - 2^{2/3} \approx 0.413$, whereas $Q^* \approx 0.467 > 0.413$.
\end{remark}

\section{Equivalence between finite support demand and continuous support demand}\label{append:equivalence}

The derivation of the necessary and sufficient conditions for the eThe derivation of the necessary and sufficient conditions for the existence of a decision-corrected arrival rate in the main body is based on demand distributions with finite support. This finite-support assumption allows us to reformulate the stochastic problem as a deterministic scenario-expanded problem (Lemma \ref{lem:scenario-expansion-eqprob}), which in turn enables the definition of dual variables and the analysis of the KKT conditions. However, this assumption can be relaxed.

In this appendix, we show that, for the purpose of characterizing the universal existence of a decision-corrected arrival rate, finite-support demand and continuously distributed demand with possibly infinite support are equivalent. In particular, we prove the following theorem, which extends the existence conditions for the decision-corrected arrival rate to the continuous infinite-support setting. Moreover, by the scenario-expansion equivalence in Lemma \ref{lem:scenario-expansion-eqprob}, Theorem \ref{thm:finite_to_inf} also implies that the continuous-time demand setting, obtained as the period length tends to zero, is equivalent to the finite discrete-time setting analyzed in the main body.

\begin{theorem}[Finite Support $\Leftrightarrow$ Continuous Support]\label{thm:finite_to_inf}
Given fixed system parameters $A, c, p,$ and $R$, consider the two-stage stochastic optimization problem
\[
\min_{b \ge 0} \Big\{\, c^\top b~T + \bbE_{\bbD \sim \mathcal{D}} \big[\pi(b, \bbD)\big] \Big\},
\]
where $\pi(b,\bbD)$ is defined by the inner problem
\[
\pi(b, \bbD) = 
\begin{array}[t]{l}
\displaystyle \min_{\{x_t \ge 0\}} \sum_{t=1}^T p^\top (D_t - R x_t) \\
\text{s.t.} \quad  
    \begin{aligned}[t]
    & A x_t \le b, \quad t = 1, \dots, T, \\
    & R x_t \le D_t, \quad t = 1, \dots, T.
    \end{aligned}
\end{array}
\]
Then, the following two statements are equivalent:
\begin{enumerate}
    \item For any distribution $\mathcal{D}$ with \textbf{finite support}, and any $T \ge 1$, there exist parameters $\{\widehat\lambda_i\}_{i=1}^T$ such that
    \[
        b^*(\widehat\lambda)\ \in\ \arg\min_{b\ge 0}\Big\{\,c^\top b~T+\bbE_{D_Z\sim\mathcal{D}}[\pi(b,D_z)]\,\Big\}.
    \]
    \item For any bounded distribution $\mathcal{D}$, whether it has finite or continuous support, and any $T \ge 1$, there exist parameters $\{\widehat\lambda_i\}_{i=1}^T$ such that
    \[
    b^*(\widehat\lambda)\ \in\ \arg\min_{b\ge 0}\Big\{\,c^\top b~T+\bbE_{D_Z\sim\mathcal{D}}[\pi(b,D_z)]\,\Big\}.
    \]
\end{enumerate}
\end{theorem}
Before showing the proof of Theorem \ref{thm:finite_to_inf}, we first prove the following lemma.
\begin{lemma}[Uniform Discrete Approximation of  $\pi(\cdot,\cdot)$]\label{lem:uniform_discrete}
Let $\mathscr{D}\subset(\mathbb R_+^{\,n})^T$ be bounded and let $B\subset\mathbb R_+^{\,m}$ be compact. 
\begin{enumerate}
    \item For each $\bbD\in \mathscr{D}$, the mapping $b\mapsto \pi(b,\bbD)$ is convex and continuous.
    \item For each $b\in B$, the mapping $\bbD\mapsto \pi(b,\bbD)$ is Lipschitz continuous and bounded uniformly over $b\in B$.
\end{enumerate}
As a consequence, for any probability measure $\mu$ supported on $\mathscr{D}$ and any $\varepsilon>0$, there exists a finite-support measure $\mu_\varepsilon$ supported on $\{D_1,\dots,D_Z\}\subset \mathscr{D}$ such that
\begin{equation*}
\sup_{b\in B}\Big|\bbE_\mu[\pi(b,D)]-\bbE_{\mu_\varepsilon}[\pi(b,D)]\Big|\le \varepsilon.
\end{equation*}
\end{lemma}

\begin{proof}{Proof of Lemma \ref{lem:uniform_discrete}}

By strong duality, each second-stage problem admits the dual representation
\[
\pi(b,\bbD)=\sum_{t=1}^T\max_{(u_t,v_t)\in\mathcal U}\big(-u_t^\top b-v_t^\top D_t+p^\top D_t\big),
\]
where
\[
\mathcal U:=\bigl\{(u,v)\in\mathbb R_+^{m}\times\mathbb R_+^{n}:A^\top u+R^\top v\ge R^\top p\bigr\}
\]
is a fixed polyhedron independent of $(b,\bbD)$. Because $b \ge 0$ and $D_t \ge 0$, the objective function is bounded from above and cannot diverge along any extreme ray of $\mathcal U$. Therefore, the maximum is finite and must be attained at one of the extreme points of $\mathcal U$. Let $\{(u^k,v^k)\}_{k=1}^K$ be the finite set of these extreme points.

Then for each $t$,
\[
\max_{(u_t,v_t)\in\mathcal U}(-u_t^\top b-v_t^\top D_t+p^\top D_t)
=\max_{1\le k\le K}\big(-(u^k)^\top b-(v^k)^\top D_t+p^\top D_t\big).
\]
Hence,
\[
\pi(b,\bbD)=\sum_{t=1}^T\max_{1\le k\le K}\big(-(u^k)^\top b-(v^k)^\top D_t+p^\top D_t\big).
\]
Because $\pi(b,\bbD)$ is a finite maximum of affine functions in $(b,\bbD)$, it is jointly continuous and convex in $b$. 
For any compact $B$ and bounded $\mathscr{D}$, each affine term is Lipschitz in $\bbD$ with constant 
\[
L:=\sum_{t=1}^T \max_{k\le K} \|p - v^k\|_2,
\]
uniform over $b\in B$.

Next, we construct the discrete finite support by empirical measure.
Given i.i.d.\ samples $D^{(1)},\dots,D^{(N)}\sim\mu$, define the empirical measure
\[
\mu_N:=\frac{1}{N}\sum_{i=1}^{N}\delta_{D^{(i)}},
\]
which is supported on the finite set $\{D^{(1)},\dots,D^{(N)}\}$.
By Varadarajan’s theorem, $\mu_N\Rightarrow \mu$ almost surely as $N\to\infty$ (weak convergence on separable metric spaces).
On bounded metric spaces the bounded–Lipschitz metric $d_{\text{BL}}$ metrizes weak convergence
(e.g., \cite{dudley2018real} Thm. 11.3.3), hence
\[
d_{\text{BL}}(\mu_N,\mu)\ \xrightarrow[N\to\infty]{\text{a.s.}}\ 0.
\]
Let $B\subset\bbR_+^m$ be compact. Since
$\mathcal F:=\{\pi(b,\cdot):b\in B\}$ is uniformly bounded and Lipschitz on $\mathcal{D}$ with envelope $M$ and
Lipschitz constant $L$ (both independent of $b\in B$).
By definition of $d_{\text{BL}}$,
\[
\sup_{b\in B}\Big|\bbE_\mu[\pi(b,D)]-\bbE_{\mu_N}[\pi(b,D)]\Big|
\ \le\ \max\{M,L\}\, d_{\text{BL}}(\mu_N,\mu)\ \xrightarrow{\text{a.s.}}\ 0.
\]
Therefore, for any $\varepsilon>0$, with probability $1$ there exists $N_\varepsilon$ such that
for all $N\ge N_\varepsilon$,
\[
\sup_{b\in B}\Big|\bbE_\mu[\pi(b,D)]-\bbE_{\mu_N}[\pi(b,D)]\Big|\le \varepsilon.
\]
Since each $\mu_N$ is finitely supported, this yields the desired discrete approximation.
\hfill\Halmos
\end{proof}

\begin{proof}{Proof of Theorem \ref{thm:finite_to_inf}}
The difference between (1) and (2) lies only in the class of distributions: finite-support versus general (possibly continuous) distributions. 
Since $b^*(\cdot)$ is the same mapping, it suffices to show that the set of minimizers
\[
\mathcal B_{\text{fin}}=\bigcup_{\text{finite }\mathcal D}\arg\min_{b\ge 0}\big(c^\top b~T+\bbE_{\mathcal D}[\pi(b,D)]\big)
\quad\text{and}\quad
\mathcal B_{\text{cont}}=\bigcup_{\text{general }\mathcal D}\arg\min_{b\ge 0}\big(c^\top b~T+\bbE_{\mathcal D}[\pi(b,D)]\big)
\]
coincide.

\paragraph{Contrapositive argument.}
Assume statement (2) fails. Then there exists a bounded-support (possibly continuous) distribution $\mathcal D$ such that
\[
b^*(\hat\lambda)\notin\arg\min_{b\ge 0} J_{\mathcal D}(b), \quad 
J_{\mathcal D}(b):=c^\top b~T+\bbE_{\mathcal D}[\pi(b,D)],
\]
for all $\hat\lambda$. 

Let $B$ be a compact set containing all possible $b^*(\hat\lambda)$ and the minimizers of $J_{\mathcal D}$ (such a set exists by coercivity of $J_{\mathcal D}$). 
Since $b^*(\hat\lambda)$ never minimizes $J_{\mathcal D}$, there exists a constant $\eta>0$ such that
\begin{equation}\label{eq:gap}
J_{\mathcal D}(b^*(\hat\lambda)) \ge \min_{b\ge 0} J_{\mathcal D}(b) + 2\eta, \qquad \forall\,\hat\lambda.
\end{equation}

By Lemma~\ref{lem:uniform_discrete}, for any $\varepsilon>0$, we can find a finite-support law $\mathcal D_\varepsilon$ such that
\begin{equation}\label{eq:uniform}
\sup_{b\in B}\big|J_{\mathcal D}(b)-J_{\mathcal D_\varepsilon}(b)\big|\le\eta.
\end{equation}
Combining \eqref{eq:gap} and \eqref{eq:uniform} gives, for all $\hat\lambda$,
\[
J_{\mathcal D_\varepsilon}(b^*(\hat\lambda))
\ge J_{\mathcal D}(b^*(\hat\lambda))-\eta
\ge \min_{b} J_{\mathcal D}(b)+\eta
\ge \min_{b} J_{\mathcal D_\varepsilon}(b),
\]
where the last inequality follows from \eqref{eq:uniform}.  
Therefore, $b^*(\hat\lambda)$ also fails to minimize $J_{\mathcal D_\varepsilon}$ for any $\hat\lambda$.  
This contradicts statement (1), which asserts the existence of some $\hat\lambda$ making $b^*(\hat\lambda)$ optimal for \emph{every} finite-support distribution.

Hence, statement (2) must hold.  
The reverse direction $(2)\Rightarrow(1)$ is immediate because finite-support distributions are a special case of general distributions.\hfill\Halmos

\end{proof}

\section{Pitfall of Stochastic Right-Hand Side}\label{append:pitfall}

The two-stage stochastic problem studied in this paper is closely related to stochastic linear programming. Specifically, our paper addresses the predict-then-optimize problem in which uncertainty appears in the right-hand side (RHS) of the constraints, and provides theoretically sound guarantees. This RHS uncertainty distinguishes our work from much of the broader predict-then-optimize literature, where uncertainty typically appears in the objective function. Indeed, RHS uncertainty is more challenging than linear objective uncertainty. This appendix provides further insight into the reasons behind this difficulty.

Particularly, we consider the following general two-stage optimization problem. We use $x^{\text{1st}}$ to denote the first stage decision, and $x^{\text{2nd}} $ to denote the second stage decision when the uncertainty is observed. The cost incurred from the first stage is a deterministic function $f(x^{\text{1st}})$, while $g(x^{\text{1st}}, x^{\text{2nd}})$ is the cost function of the second stage problem, which depends on some distribution of uncertainties. The overall problem is formulated as:

\begin{align*}
    \min_{x^{\text{1st}} }:~\left\{ f(x^{\text{1st}}) + \bbE\left[\min_{x^{\text{2nd}} }: g(x^{\text{1st}}, x^{\text{2nd}})\right]\right\}.
\end{align*}

We define the general extended-value function as follows:
\[
I_C(x) =
\begin{cases}
0, & \text{if } x \in C, \\
+\infty, & \text{if } x \notin C.
\end{cases}
\]

By the extended-value function, one example of the second-stage cost can be written as
\begin{equation}\label{append:g}
g(x^{\text{1st}}, x^{\text{2nd}}) =
I_{\{A_1 x^{\text{1st}} + A_2 x^{\text{2nd}} \le b\}}\left(
c^\top x^{\text{2nd}}\right)
\;\Longleftrightarrow\;
\begin{array}{ll}
\displaystyle & c^\top x^{\text{2nd}} \\
\text{s.t.} & A_1x^{\text{1st}} + A_2x^{\text{2nd}} \le b.
\end{array}
\end{equation}

In \eqref{append:g}, $x^{\text{2nd}}\in \bbR^m$ is the decision variable, $c \in \bbR^m$ is the coefficient in the objective, $b \in \bbR^k$ is the RHS in constraints and $A_1$ and  $A_2\in \bbR^{m\times k}$ are the known constraint matrix.

We compare two stochastic two-stage problems:
\begin{align}
\min_{x^{\text{1st}}}\;& f(x^{\text{1st}}) + \mathbb{E}_{c\sim \mathcal{D}_c}\Big[\min_{x^{\text{2nd}}} c^\top x^{\text{1st}} \;\text{s.t.}\; A_1 x^{\text{1st}} + A_2 x^{\text{2nd}} \le b \Big],
\tag{SP1}\label{append:p1}\\[6pt]
\min_{x^{\text{1st}}}\;& f(x^{\text{1st}}) + \mathbb{E}_{b\sim \mathcal{D}_b}\Big[\min_{x^{\text{2nd}}} c^\top x^{\text{2nd}} \;\text{s.t.}\; A_1 x^{\text{1st}} + A_2 x^{\text{2nd}} \le b \Big].\label{append:p2}\tag{SP2}
\end{align}
In \eqref{append:p1}, the uncertainty lies in the objective coefficients $c$, while in \eqref{append:p2} it lies in the RHS vector $b$. We suppress contextual features for notational simplicity. All statements extend to the contextual setting where one conditions on side information.

\paragraph{Are \eqref{append:p1} and \eqref{append:p2} equally difficult?}
At first glance, duality seems to suggest equivalence. That conclusion is incorrect for two reasons: the min–expectation interchange is not innocuous when the feasible region itself is random, and dual variables depend on the realized scenario. These two effects break the naive equivalence.

\paragraph{Why objective uncertainty is easier.}
When $c$ is random but the feasible region is deterministic, the second-stage optimal value is linear in $c$:
\[
\min_{x^{\text{2nd}}\in X(x^{\text{1st}})} c^\top x^{\text{1st}}.
\]
Therefore,
\[
\mathbb{E}\Big[\min_{x^{\text{2nd}}\in X(x^{\text{1st}})} c^\top x^{\text{1st}}\Big]
=
\min_{x^{\text{2nd}}\in X(x^{\text{1st}})} \mathbb{E}[c]^\top x^{\text{1st}},
\]
so a \emph{point prediction} $\hat c \approx \mathbb{E}[c]$ suffices. This reduction collapses \eqref{append:p1} to a single deterministic optimization with objective $\mathbb{E}[c]^\top x^{\text{1st}}$. Decision-focused learning can target $\mathbb{E}[c\mid \text{features}]$ and inherits known guarantees; see, for example, \cite{elmachtoub2022smart,liu2023active,hu2022fast}.

\paragraph{Why RHS uncertainty is harder.}
Let $v(b;x^{\text{1st}})$ denote the second-stage optimal value in \eqref{append:g}. For fixed $x^{\text{1st}}$, the map $b\mapsto v(b;x^{\text{1st}})$ is a \emph{convex} piecewise-linear function of $b$ (the support function of a polyhedron). Consequently,
\begin{equation}\label{eq:jensen}
\mathbb{E}\big[v(b;x^{\text{1st}})\big] \;\ge\; v\big(\mathbb{E}[b];x^{\text{1st}}\big),
\end{equation}
with strict inequality in general. Equality requires $v(\cdot;x^{\text{1st}})$ to be affine on the support of $\mathcal{D}_b$, which is rare unless the recourse is trivial. Therefore a point prediction $\hat b\approx \mathbb{E}[b]$ is \emph{insufficient}: minimizing with $\hat b$ solves $\min f(x^{\text{1st}})+v(\mathbb{E}[b];x^{\text{1st}})$, which is a lower bound to the true objective $f(x^{\text{1st}})+\mathbb{E}[v(b;x^{\text{1st}})]$ by \eqref{eq:jensen}. It implies that the optimal $x^{\text{1st}}$ can be different. 

Moreover, the feasible set for $x^{\text{2nd}}$ depends on the realization of $b$. Thus, the dual of the second-stage problem introduces dual variables $y(b;x^{\text{1st}})$ that depend on the realized scenario $b$. Some scenarios may be infeasible. The expected dual objective is $\mathbb{E}\big[b^\top y(b;x^{\text{1st}})\big]$, which again requires distributional knowledge of $b$ to evaluate or learn, so the symmetry suggested by looking only at algebraic positions of $c$ and $b$ is misleading.

In the existing literature that addresses the RHS uncertainty, e.g., \cite{iceo}, they need to estimate the entire distributional information, instead of just a single plug-in prediction. We are the first work to check if a single plug-in prediction exists for the stochastic RHS, and provide algorithms to find such a decision-unbiased prediction.

\section{Characterization of set $\mathcal{B}$}\label{subsec:mip}

We recall that in Section \ref{subsec:dist-dependent}, when studying the distribution-dependent existence,  we have that $b\in\mathcal{B}$ iff there exist $y_1,\dots,y_T\in\mathbb{R}_+^m$ satisfying (B1)–(B3).
Let $h(j)\in\{1,\dots,m\}$ denote the pool served by job $j$, and $i(j)\in\{1,\dots,n\}$ the class of job $j$.

\paragraph{Feasibility certificate via MIP.}
Given horizon $T$, matrix $A\in\mathbb{R}_+^{m\times k}$, costs $c\in\mathbb{R}_+^m$, penalties $p\in\mathbb{R}_+^n$, and a capacity vector $b\in\mathbb{R}_+^m$.

Choose Big-$M$ constants
\[
M^y_h := c_h ~ T,\qquad 
M^r_j := (A^\top c~T)_j ,\qquad
M^p_j := (A^\top c ~T)_j .
\]

\emph{Decision variables} (for $t=1,\dots,T$):
\[
y_{t,h}\ge 0\ (h=1,\dots,m),\quad
r_{t,j}\ge 0\ (j=1,\dots,k),\quad
w_{t,h}\in\{0,1\},\quad
v_{t,j}\in\{0,1\}.
\]
\begin{align*}
\min\quad & 0\\[2mm]
\text{s.t.}\quad
&\sum_{t=1}^T y_{t,h} \;\le\; c_h~T && \forall\,h=1,\dots,m \tag{C1}\\
&\sum_{t=1}^T y_{t,h} \;=\; c_h~T && \forall\,h \ \text{with}\ b_h>0 \tag{C2}\\
& r_{t,j} \;=\; \sum_{h=1}^m A_{h j}\, y_{t,h} && \forall t=1,\dots,T, \forall j=1,\dots,k \tag{C3}\\
& 0 \;\le\; y_{t,h} \;\le\; M^y_h\, w_{t,h} && \forall t=1,\dots,T, \forall h=1,\dots,m \tag{C4}\\
& \sum_{\{j:\, h(j)=h\}} v_{t,j} \;\ge\; w_{t,h} && \forall t=1,\dots,T, \forall h=1,\dots,m \tag{C5}\\
& v_{t,j} \;\le\; w_{t,\,h(j)} && \forall t=1,\dots,T, \forall j=1,\dots,k \tag{C6}\\
& r_{t,j} \;\le\; r_{t,k} + M^r_j\,(1 - v_{t,j}) && \forall t=1,\dots,T, \forall j=1,\dots,k \\ &&& \forall\,j'\ \text{with}\ i(j')=i(j) \tag{C7}\\
& r_{t,j} \;\le\; p_{\,i(j)} + M^p_j\,(1 - v_{t,j}) && \forall t=1,\dots,T, \forall j=1,\dots,k \tag{C8}\\[1mm]
& y_{t,h}\ge0,\ \ r_{t,j}\ge0,\ \ w_{t,h}\in\{0,1\},\ \ v_{t,j}\in\{0,1\} && \forall t=1,\dots,T, \forall j=1,\dots,k \\ &&&\forall h=1,\dots,m \, .
\end{align*}

\begin{proposition}\label{prop:support-feasible}
The MIP above is feasible if and only if $b\in\mathcal{B}$.
\end{proposition}

\begin{proof}{Proof of Proposition \ref{prop:support-feasible}}
($\Rightarrow$) Suppose the MIP is feasible and take any feasible solution.
Let $y_t$ be the vectors with components $y_{t,h}$. 
Then (C1) gives $\sum_t y_t\le c~T$ (i.e., (B1)); (C2) gives $\sum_t (y_t)_h=c_h~T$ for every $h$ with $b_h>0$ (i.e., (B2)).
By (C3), $r_t=A^\top y_t$.
Fix $t$ and $h$ with $(y_t)_h>0$. 
From (C4) we must have $w_{t,h}=1$, and then (C5) ensures there exists at least one $j$ with $h(j)=h$ and $v_{t,j}=1$.
For such a witness $j$, (C7) yields $r_{t,j}\le r_{t,k}$ for all $k$ in the same class as $j$, i.e., $j$ attains the classwise minimum of $A^\top y_t$; and (C8) yields $r_{t,j}\le p_{i(j)}$.
Thus (B3) holds.
Hence $y_1,\dots,y_T$ satisfy (B1)–(B3), so $b\in\mathcal{B}$.

($\Leftarrow$) Conversely, assume $b\in\mathcal{B}$. 
Then there exist $y_1,\dots,y_T\in\mathbb{R}_+^m$ satisfying (B1)–(B3).
Set $y_{t,h}$ accordingly, let $r_{t,j}:=(A^\top y_t)_j$, and define 
\[
w_{t,h}:=\begin{cases}1,& (y_t)_h>0,\\ 0,& \text{otherwise.}\end{cases}
\]
For each $t$ and each $h$ with $(y_t)_h>0$, pick a witness $j$ guaranteed by (B3) (so $h(j)=h$, $r_{t,j}\le r_{t,k}$ for all $k$ with $i(k)=i(j)$, and $r_{t,j}\le p_{i(j)}$) and set $v_{t,j}=1$; set $v_{t,j}=0$ for all other $j$.
Then (C1) holds by (B1); (C2) holds by (B2); (C3) holds by construction; (C4)–(C6) hold by the definitions of $w$ and $v$; (C7)–(C8) hold by the witness properties from (B3) (the Big-$M$ terms are inactive when $v_{t,j}=1$, and safely relax the inequalities when $v_{t,j}=0$ since $r_{t,j}\le (A^\top c~T)_j=M^p_j$).
Therefore, the MIP is feasible.
\hfill\Halmos
\end{proof}

Note that the above MIP uses a zero objective function. In practice, when multiple feasible solutions for $y_{t,h}$ exist, we may modify the objective to impose a tie-breaking rule. The decision-corrected method in our experiments uses the following objective in place of the zero objective above, which penalizes successive differences in $y_{t,h}$ across consecutive hours for each resource and thereby produces a smoother, more interpretable arrival-rate profile. 

\begin{align*}
\min\quad & \sum_{h=1}^m\delta_h\\[2mm]
\text{s.t.}\quad
& |y_{t,h} - y_{t+1,h}| \;\le\; \delta_h, && \forall\,h=1,\dots,m, \forall t = 1,\cdots, T-1,\\
& (C1) - (C8),\\
& y_{t,h}\ge0,\ \ r_{t,j}\ge0,\ \ w_{t,h}\in\{0,1\},\ \ v_{t,j}\in\{0,1\} && \forall t=1,\dots,T, \forall j=1,\dots,k \\ &&&\forall h=1,\dots,m \, .
\end{align*}

The above formulation produces a more stable sequence of $y_{t,h}$ values, which in turn yields a smoother decision-corrected arrival rate, as shown in Figure~\ref{fig:forecast}.

\section{Omitted Proofs}

\begin{proof}{Proof of Lemma \ref{lem:fluid-kkt}}
Fix $b\ge 0$. The inner problem in \eqref{eq:fa} is
\[
\min_{\{x_t\ge 0\}} \ \sum_{t=1}^T \Big(p^\top\lambda_t - (R^\top p)^\top x_t\Big)
\quad\text{s.t.}\quad
A x_t \le b,\ \ R x_t \le \lambda_t \qquad \forall t=1,\dots,T.
\]
Introduce Lagrange multipliers $y_t\ge 0$ for $A x_t\le b$, $z_t\ge 0$ for $R x_t\le \lambda_t$, and $s_t\ge 0$ for $x_t\ge 0$ (equivalently $-x_t\le 0$). The Lagrangian is
\[
\mathcal{L}(\{x_t\};\{y_t,z_t,s_t\})
=\sum_{t=1}^T\!\Big[
p^\top\lambda_t - (R^\top p)^\top x_t
+ y_t^\top(A x_t - b) + z_t^\top(R x_t - \lambda_t) - s_t^\top x_t
\Big].
\]
Collecting the $x_t$ terms gives
\[
\mathcal{L}
=\sum_{t=1}^T\!\Big[
p^\top\lambda_t - y_t^\top b - z_t^\top\lambda_t
+ \big(A^\top y_t + R^\top z_t - R^\top p - s_t\big)^\top x_t
\Big].
\]
The infimum of $\mathcal{L}$ over $x_t\in\mathbb{R}^k$ is finite if and only if there exists $s_t\ge 0$ such that
\[
A^\top y_t + R^\top z_t - R^\top p - s_t \geq 0,
\]
which is equivalent (after eliminating $s_t$) to
\[
A^\top y_t + R^\top z_t \ \ge\ R^\top p \qquad \forall t=1,\dots,T.
\]
Under this condition the dual function equals $\sum_{t=1}^T\big(p^\top\lambda_t - y_t^\top b - z_t^\top\lambda_t\big)$; otherwise it is $-\infty$. Hence, for fixed $b$ the inner dual is
\[
\max_{\substack{y_t\ge 0,\,z_t\ge 0}}
\ \sum_{t=1}^T \lambda_t^\top(p - z_t) \ -\ b^\top\sum_{t=1}^T y_t
\quad\text{s.t.}\quad
A^\top y_t + R^\top z_t \ \ge\ R^\top p \qquad \forall t=1,\dots,T .
\]
The inner primal is feasible (take $x_t\equiv 0$) and, because $R x_t\le \lambda_t$ with $p\ge 0$, has a finite optimal value; thus strong duality identifies $\pi(b,\bblambda)$ with the optimal value of this dual problem.

Substituting this representation into the outer minimization yields
\[
\min_{b\ge 0}\ \max_{\substack{y_t\ge 0,\,z_t\ge 0\\ A^\top y_t + R^\top z_t \ge R^\top p}}
\Big\{\, \sum_{t=1}^T \lambda_t^\top(p - z_t) \ +\ (c ~ T-\textstyle\sum_{t=1}^T y_t)^\top b \,\Big\}.
\]
For fixed $\{y_t,z_t\}$, consider the inner minimization
\[
\min_{b\ge 0}\ \big(c ~ T-\textstyle\sum_{t=1}^T y_t\big)^\top b .
\]
If some coordinate satisfies $\big(c ~ T-\sum_{t=1}^T y_t\big)_h<0$, sending $b_h\to+\infty$ drives the value to $-\infty$.
Thus finiteness forces
\[
\sum_{t=1}^T y_t \ \le\ c ~ T \qquad \text{(componentwise)}.
\]
Under this condition, all coefficients $\big(c ~ T-\sum_{t=1}^T y_t\big)_h$ are nonnegative, so the minimum over $b\ge 0$ equals $0$.
Moreover, the minimizers are characterized componentwise:
\[
\big(c ~ T-\textstyle\sum_{t=1}^T y_t\big)_h>0 \ \Longrightarrow\ b_h=0,
\qquad
\big(c ~ T-\textstyle\sum_{t=1}^T y_t\big)_h=0 \ \Longrightarrow\ b_h\ \ge 0.
\]
Equivalently,
\[
b_h>0 \ \Longrightarrow\ \big(c ~ T-\textstyle\sum_{t=1}^T y_t\big)_h=0,
\]
which is the vector identity
\[
b^\top\!\Big(c ~ T-\sum_{t=1}^T y_t\Big)=0.
\]
This is the complementary slackness condition for $b$.

With $\sum_{t=1}^T y_t \le c ~ T$ enforced, the dual of the full two-stage problem is
\[
\begin{aligned}
\max_{\{y_t\ge 0,\,z_t\ge 0\}} \quad & \sum_{t=1}^T \lambda_t^\top(p - z_t) \\
\text{s.t.}\quad & \sum_{t=1}^T y_t \le c ~ T,\qquad A^\top y_t + R^\top z_t \ge R^\top p \quad \forall t=1,\dots,T.
\end{aligned}
\]
Primal feasibility is the set of constraints in \eqref{eq:fa}; dual feasibility is the system just displayed. Complementary slackness pairs each multiplier with the slack of its constraint: $y_t^\top(Ax_t-b)=0$ and $z_t^\top(Rx_t-\lambda_t)=0$ from the inner inequalities; $b^\top\!\big(c ~ T-\sum_{t=1}^T y_t\big)=0$ from the outer minimization in $b$; and, using the stationarity identity $s_t=A^\top y_t + R^\top z_t - R^\top p\ge 0$, the relation $x_t^\top s_t=0$ yields $x_t^\top\!\big(A^\top y_t + R^\top z_t - R^\top p\big)=0$. Because the overall problem is a feasible linear program and strong duality holds, these (P), (D), and (CS) conditions are necessary and sufficient for optimality. \hfill \Halmos
\end{proof}

\begin{proof}{Proof of Lemma \ref{lem:scenario-expansion-eqprob}}
Because the support of $\bbD$ is $\{D_{1:T}^{(z)}\}_{z=1}^Z$, the linearity of expectation gives
\(
\mathbb E_{\bbD\sim\mathcal D}\big[\pi(b,\bbD)\big]
=\sum_{z=1}^Z p_z\,\pi\big(b,D_{1:T}^{(z)}\big)
\)
for every fixed $b$. This establishes that Problem \ref{eq:outer_problem} is equivalent to $\min_{b\ge0}\; c^\top b ~ T+\sum_{z=1}^Z p_z\,\pi\!\big(b,D_{1:T}^{(z)}\big)$, the first result in \eqref{eq:equal-weight-expansion}. 

For the equal–weight form, construct integers $N_z\ge 0$ and $\tilde Z:=\sum_z N_z$ such that $N_z/\tilde Z=p_z$ (e.g., by clearing denominators when the $p_z$ are rational; more generally, duplicate paths in proportions approaching $p_z$ and use that the objective is linear in the distribution). Duplicating each path $D_{1:T}^{(z)}$ exactly $N_z$ times produces the multiset $\{D_{1:T}^{(\tilde z)}\}_{\tilde z=1}^{\tilde Z}$ with
\(
\sum_{z=1}^Z p_z\,\pi\!\big(b,D_{1:T}^{(z)}\big)
=\frac{1}{\tilde Z}\sum_{\tilde z=1}^{\tilde Z}\pi\!\big(b,D_{1:T}^{(\tilde z)}\big)
\)
for every $b$.
Reindexing $(t,\tilde z)\mapsto\tau$ is purely notational and does not alter the objective or minimizers.

Finally, to move the summation  $\sum_{\tilde z=1}^{\tilde Z}$ inside the function $\pi(\cdot,\cdot)$, we observe that in the definition of $\pi(\cdot,\cdot)$ in \eqref{eq:inner_problem}, both variables and constraints are decomposable by time period $t$. It implies that the $\min$ and summation over $t$ are interchangeable in \eqref{eq:inner_problem}. Thus, we have that $\frac{1}{\tilde Z}\sum_{\tilde z=1}^{\tilde Z}\pi\!\big(b,D_{1:T}^{(\tilde z)}\big) = \frac{1}{\tilde Z}\pi\!\big(b,\{D_\tau\}_{\tau=1}^{T\tilde Z}\big)$. To obtain  Lemma \ref{lem:scenario-expansion-eqprob}, we divide the entire objective function by $T$, and move the constant $T$ to the front of the objective as a multiplier. 
\hfill \Halmos
\end{proof}

\begin{proof}{Proof of Theorem \ref{thm:existence-iff-constant}}
We prove the necessary and sufficient conditions separately. First, we prove A $\Rightarrow$ B via contrapositive.

\emph{(A $\Rightarrow$ B) via contrapositive.}
Assume $\lnot\textbf{(B)}$: there exists an expanded collection $\{D_\tau\}_{\tau=1}^{T\tilde Z}$ such that
\[
\forall\lambda\in\mathbb{R}_+^{T\tilde Z}:\quad
b^*\!\big(\lambda^{(TZ)}\big)\ \notin\ \arg\min_{b\ge0}\Big\{c^\top b+\tfrac{1}{Z T}\pi\big(b,\{D_\tau\}_{\tau=1}^{T\tilde Z}\big)\Big\}.
\]

To show $\lnot\textbf{(A)}$, we consider the planning period $T$ to be 1, and construct the demand $\widetilde{\mathcal D}$ according to the equal weight scenarios $\{D_\tau\}_{\tau=1}^{T\tilde Z}$ for the single period. Then, by Lemma \ref{lem:scenario-expansion-eqprob}, Problem \eqref{eq:outer_problem} is equivalent to $\min_{b\ge0}\Big\{c^\top b+\tfrac{1}{Z T}\pi\big(b,\{D_\tau\}_{\tau=1}^{T\tilde Z}\big)\Big\}$. By $\lnot\textbf{(B)}$, there is no constant decision-corrected arrival rate $\lambda$. In other words, we cannot find a single period arrival rate that is decision-unbiased. Consequently, for demand distribution $\widetilde{\mathcal D}$ and $T=1$, there is no single period decision-corrected arrival rate $\lambda$, i.e., we obtain $\lnot\textbf{(A)}$.

\emph{(B $\Rightarrow$ A) via contrapositive.}
Assume $\lnot\textbf{(A)}$: there exists a planning period $T$ and a finite-support distribution $\mathcal D$ such that
\[
\forall\,\bblambda \in \bbR^T_+:\quad
b^*(\bblambda)\ \notin\ \arg\min_{b\ge0}\Big\{c^\top b ~ T+\mathbb{E}_{\bbD\sim\mathcal D}\big[\pi(b,\bbD)\big]\Big\}.
\]
Using Lemma \ref{lem:scenario-expansion-eqprob}, the objective is equivalent to the following problem with $Z T$ periods
$$\min_{b\ge0}\Big\{c^\top b+\tfrac{1}{Z T}\pi\big(b,\{D_\tau\}_{\tau=1}^{T\tilde Z}\big)\Big\}.$$
If \textbf{(B)} holds, then we can find a constant arrival rate $\widehat\lambda \in \bbR^{TZ}_+$ that is decision-unbiased. Then, we can set $\bblambda$ in \textbf{(A)} to be the same value as the constant entry in $\widehat\lambda$. Then, this $\bblambda$ will be decision-unbiased for the original distribution $\mathcal D$ in \textbf{(A)}. It contradicts $\lnot\textbf{(A)}$. Hence, no such $\widehat\lambda$ exists, and $\lnot\textbf{(B)}$ holds.
\hfill \Halmos
\end{proof}

\begin{proof}{Proof of Theorem \ref{thm:const-characterization}}
By Lemma~\ref{lem:fluid-kkt}, when \(\widehat\lambda_t\equiv\widehat\lambda\) for all \(t\),
the \(T\) identical inner problems coincide. The objective reduces to \(c^\top b + \,p^\top(\widehat\lambda - R x)\)
with constraints \(A x\le b\), \(R x\le \widehat\lambda\), i.e., \eqref{eq:det-const-primal}.
The KKT conditions for \eqref{eq:det-const-primal} are (componentwise):
\[
\begin{aligned}
&\textbf{P: } b\!\ge0,\ x\!\ge0,\ A x\!\le b,\ R x\!\le \widehat\lambda.\\
&\textbf{D: } y\!\ge0,\ z\!\ge0,\ \ c-y\!\ge0,\ A^\top y + R^\top z \;\ge\; \,R^\top p.\\
&\textbf{CS: } y(b-Ax)=0,\ z(\widehat\lambda-Rx)=0,\ b(c-y)=0,\
   x\big(A^\top y + R^\top z - \,R^\top p\big)=0.
\end{aligned}
\]

\noindent\textbf{(B) $\Rightarrow$ (A).}
Fix any \(\mathcal{D}\) and let \(b^*=b^*(\mathcal{D})\) be expectation–optimal.
For each pool \(h\) with \(b_h^*>0\), pick a class \(i(h)\) and a column \(j(h)\)
that (i) uses pool \(h\), (ii) attains the classwise minimum
\(\min_{k:\,i(k)=i(h)}(A^\top c)_k\), and (iii) satisfies the penalty cap
\((A^\top c)_{j(h)}\le \,p_{i(h)} \le \,p_i\).
Define \(x\) by tightening capacity at each resource pool with strictly positive capacity:
\(A_{h,\,j(h)}\,x_{j(h)} = b_h^*\) (equivalently \(x_{j(h)}=b_h^*/A_{h,\,j(h)}\)); set other
components of \(x\) to zero. For pools $h$ with $b_h^* =0$ no job is routed through $h$, so $(Ax)_h=0=b_h^*$ and the slackness term $y_h(b_h^*-(Ax)_h)$ vanishes regardless of $y_h$. The choice $y=c$ therefore remains dual-feasible and complementary on these pools. For any customer class served only by such zero-capacity pools, set $z_i$ large enough to maintain $A^\top y+ R^\top z \geq R^\top p$ on its columns (this is always possible and these columns carry $x_j=0$, so the slackness term $x\big(A^\top y + R^\top z - \,R^\top p\big)=0$. Let \(\widehat\lambda := R x\).

We verify the KKT system of \eqref{eq:det-const-primal} with primal \((b^*,x)\).
Choose \(y:=c\) and define \(z\in\mathbb{R}_+^{n}\) by
\(z_i := \,p_i - (A^\top c)_{j_i}\) for each served class \(i\), where \(j_i\) denotes the
selected column serving class \(i\); for unserved classes, take \(z_i\) large enough
for dual feasibility. Since each \(j_i\) minimizes \((A^\top c)\) over its class, for any
\(k\) with \(i(k)=i\), \((A^\top c)_k \ge (A^\top c)_{j_i}\). Hence, for any column \(k\)
serving class \(i(k)\),
\[
(A^\top y)_k + (R^\top z)_k
= (A^\top c)_k + z_{i(k)}
\;\ge\; (A^\top c)_{j_{i(k)}} + \big(\,p_{i(k)} - (A^\top c)_{j_{i(k)}}\big)
= \,p_{i(k)},
\]
so stationarity \(A^\top y + R^\top z \ge \,R^\top p\) holds with equality on each
chosen job \(j_i\) (hence on every \(j\) with \(x_j>0\)).
We also have \(c-y=0\), so \(b^{*\top}(c-y)=0\); \(A x=b^*\) on pools with positive capacity gives
\(y^\top(b-Ax)=0\); finally \(R x=\widehat\lambda\) yields \(z^\top(\widehat\lambda-Rx)=0\).
Therefore \((b^*,x;y,z)\) satisfies the KKT conditions, proving optimality of \(b^*\)
for \eqref{eq:det-const-primal}. This establishes (A).

\noindent\textbf{(A) $\Rightarrow$ (B).}
Assume (A) is true. This states that for any integer $T \ge 1$ and any distribution $\mathcal{D}$, if $b^*=b^*(\mathcal{D})$ is an optimal solution, then $b^*$ must also be an optimal solution to the constant-profile fluid program \eqref{eq:det-const-primal} for some $\widehat\lambda$.
Now, let's establish (B). Pick any resource pool $h$ that has positive capacity in some expectation-optimal solution, i.e., there exists a $\mathcal{D}$ such that $b_h^* = b_h^*(\mathcal{D}) > 0$.
Fix any integer $T\ge1$. Let $(x^*, y, z)$ be the corresponding optimal primal-dual solution for this problem.
By Complementary Slackness (CS) on $b$, $b_h^* > 0 \implies (c-y)_h = 0$, so $y_h = c_h$.
Also by CS on $b$, $b_h^* > 0 \implies (b^* - Ax^*)_h = 0$. This means at least one job $j$ using pool $h$ must be active, i.e., $h(j)=h$ and $x_j^* > 0$.
By CS on $x$, this active $j$ must satisfy the reduced-cost equality. That is $(A^\top y)_j + z_i =  p_i$, where $i=i(j)$.
Since $y_h=c_h$ and $j$ uses pool $h$, $(A^\top y)_j = (A^\top c)_j$.
Substituting this into the equality gives $(A^\top c)_j + z_i = p_i$.
By Dual Feasibility (D), $z_i \ge 0$, which requires ${(A^\top c)_j \le p_{i(j)}}$. Since it must hold for any integer $T\ge1$, we need ${(A^\top c)_j \le p_{i(j)}}$. This establishes the penalty cap condition in (B).
Now consider any other job $k$ in the same class as $j$ (i.e., $i(k)=i(j)$). The dual feasibility constraint for $k$ is $(A^\top y)_k + z_i \ge p_i$.
Substitute $(A^\top y)_k = (A^\top c)_k$ and $z_i = p_i - (A^\top c)_j$ (from the equality) into this inequality:
\[
(A^\top c)_k + \big(p_i - (A^\top c)_j\big) \ge p_i.
\]
This simplifies to $(A^\top c)_k \ge (A^\top c)_j$. This establishes that $j$ must be a classwise-minimum job.
Since $j$ was an arbitrary active job for a resource $h$ that has positive capacity, we have shown that any such resource $h$ must have at least one job $j$ that satisfies both conditions in (B). This concludes the proof of (B). \hfill \Halmos
\end{proof}

\begin{proof}{Proof of Theorem \ref{thm:dist-dependent-existence}}
By Lemma~\ref{lem:fluid-kkt}, the KKT system for the above fluid program is:
\[
\begin{aligned}
&\textbf{P: } b\!\ge0,\ x_t\!\ge0,\ A x_t\!\le b,\ R x_t\!\le \widehat\lambda_t\quad (\forall t=1,\dots,T).\\
&\textbf{D: } y_t\!\ge0,\ z_t\!\ge0,\ \ \sum_{t=1}^T y_t \le c ~T,\ \ A^\top y_t + R^\top z_t \ge R^\top p\quad(\forall t=1,\dots,T).\\
&\textbf{CS: } y_t(b-Ax_t)=0,\ z_t(\widehat\lambda_t-Rx_t)=0,\ b\Big(c ~T-\sum_{t=1}^T y_t\Big)=0,\ 
               x_t\big(A^\top y_t + R^\top z_t - R^\top p\big)=0\quad(\forall t=1,\dots,T).
\end{aligned}
\]

\noindent\textbf{(B) $\Rightarrow$ (A).}
Given $\{y_t\}$ satisfying \eqref{eq:B1},\eqref{eq:B2}, and \eqref{eq:B3}, define
for each time period $t$ and customer class $i$ the duals
\[
z_{t,i}\ :=\ (p_i - \min_{k:\,i(k)=i}(A^\top y_t)_k)^+  \ \ge 0.
\]
Then $A^\top y_t+R^\top z_t\ge R^\top p$ with equality on every job column $j$ that attains the minimum for its class. For each period $t$ and each resource pool $h$ with $(y_t)_h>0$,
pick $j_t(h)$ with $h(j_t(h))=h$ as in \eqref{eq:B3}, let $i_t(h)=i(j_t(h))$, and set
\[
A_{h,\,j_t(h)}\,x_{t,\,j_t(h)}\ =\ b^*_h,
\]
with all other components of $x_t$ zero. For any pool $h$ with $b^*_h = 0$, no job is routed through $h$ (the corresponding $x_{t,j}=0$ for all $j$ with $h(j)=h$), so the capacity constraint $(Ax_t)_h \leq b^*_h$ holds with equality at $0$ and the complementary slackness term $y_{t,h}(b_h^*-(Ax_t)_h)=0$ is satisfied trivially. Hence the KKT conditions hold across all pools. Define $\widehat\lambda_t:=R x_t$ so that
$R x_t\le \widehat\lambda_t$ is tight and hence $z_t^\top(\widehat\lambda_t-Rx_t)=0$.
Because $(y_t)_h>0$ only where we tightened $(b-Ax_t)_h=0$, we have $y_t(b-Ax_t)=0$.
Finally, \eqref{eq:B2} gives $b^{*\top}\big(c ~T-\sum_t y_t\big)=0$.
All KKT conditions hold with primal $b^*$, thus $b^*$ is optimal and (A) follows.

\noindent\textbf{(A) $\Rightarrow$ (B).}
Assume (A) and let $(b^*,\{x_t^*\};\{y_t,z_t\})$ be primal–dual optimal.
KKT gives $\sum_t y_t\le c ~T$ and $b^{*\top}\!\big(c ~T-\sum_t y_t\big)=0$, yielding
\eqref{eq:B1}, and \eqref{eq:B2}.  
If $(y_t)_h>0$, then $y_t^\top(b-Ax_t^*)=0$ implies $(A x_t^*)_h=b^*_h$, so some
$j$ with $h(j)=h$ has $x_{t,j}^*>0$; let $i=i(j)$. Reduced-cost equality on $j$ gives
\((A^\top y_t)_j + z_{t,i} = p_i\) with $z_{t,i}\ge0$, hence \((A^\top y_t)_j \le p_i\).
For any $k$ with $i(k)=i$, the inequality \(A^\top y_t + R^\top z_t \ge R^\top p\) yields
\((A^\top y_t)_k \ge (A^\top y_t)_j\), so $j$ attains the classwise minimum.
Thus \eqref{eq:B3} holds and (B) is established. \hfill \Halmos
\end{proof}

\begin{proof}{Proof of Proposition \ref{prop:support-B}}
By Definition~\ref{def:B}, $b^0\in\mathcal{B}$ iff there exist $y_1,\dots,y_T\in\mathbb{R}_+^{m}$ such that \eqref{eq:B1}, \eqref{eq:B2}, and \eqref{eq:B3} hold. Among these, only \eqref{eq:B2} involves $b^0$, and it does so solely through the support $H(b^0)=\{h:\,b_h>0\}$: it requires $\sum_{t=1}^T (y_t)_h=c_h ~ T$ for every $h\in H(b^0)$.

(1) If $H(b^0)=H(\tilde b)$, the set of equalities required by \eqref{eq:B2} is identical for $b^0$ and $\tilde b$; therefore, any sequence $\{y_t\}$ that certifies one certifies the other, proving equivalence.

(2) If $H(\tilde b)\subseteq H(b^{0})$ and $b^{0}\in\mathcal{B}$, take the sequence $\{y_t\}$ for $b^{0}$. They satisfy $\sum_t (y_t)_h=c_h ~T$ for all $h\in H(b^{0})$ and hence for all $h\in H(\tilde b)$. With (B1)–(B3) otherwise unchanged, the same $\{y_t\}$ witness $\tilde b\in\mathcal{B}$. \hfill\Halmos
\end{proof}

\begin{proof}{Proof of Theorem~\ref{thm:alg}.}
Let $\bbD^{(1)},\ldots,\bbD^{(Z)}$ be i.i.d.\ demand profiles with common distribution $\mathcal D$.
By assumption, there exists $\bar D\in\mathbb R^n_+$ such that $0\le D_t\le \bar D$ componentwise for all $t=1,\ldots,T$ almost surely. Define
\begin{align*}
\mathscr D
:=
\Bigl\{
\bbD=(D_1,\ldots,D_T)\in(\mathbb R^n_+)^T:
0\le D_t\le \bar D,\ t=1,\ldots,T
\Bigr\}.
\end{align*}
For a set $\Lambda\subseteq(\mathbb R^n_+)^T$, write
\begin{align*}
\operatorname{dist}(\bblambda,\Lambda)
:=
\inf_{\tilde\bblambda\in\Lambda}
\|\bblambda-\tilde\bblambda\|_2 .
\end{align*}

Since $x_t=0$ is feasible for the second-stage problem for every $b\ge0$ and every $\bbD\in\mathscr D$, and since $p\ge0$, we have
\begin{align*}
0
\le
\pi(b,\bbD)
\le
\sum_{t=1}^T p^\top D_t
\le
\sum_{t=1}^T p^\top \bar D
=:
\bar\pi,
\qquad
b\ge0,\ \bbD\in\mathscr D .
\end{align*}
Let $\widehat b^{(Z)}$ be the empirical optimizer selected in Step 2 of Algorithm~\ref{alg:1}. Define the population and empirical objectives
\begin{align*}
L(b)
&:=
c^\top T b+\mathbb E_{\bbD\sim\mathcal D}[\pi(b,\bbD)], \\
L_Z(b)
&:=
c^\top Tb+
\frac1Z\sum_{z=1}^Z
\pi\bigl(b,\bbD^{(z)}\bigr),
\qquad b\ge0 .
\end{align*}
Because $c\in\mathbb R^m_{++}$, both $L$ and $L_Z$ are coercive on $\mathbb R^m_+$. Moreover,
\begin{align*}
L(0)\le\bar\pi,
\qquad
L_Z(0)\le\bar\pi .
\end{align*}
Therefore, every population minimizer and every empirical minimizer belongs to the compact set
\begin{align*}
B_0
:=
\Bigl\{
b\in\mathbb R^m_+:
 c^\top T b\le\bar\pi
\Bigr\}.
\end{align*}
Hence, all first-stage minimizations may be restricted to $B_0$.

By Lemma~\ref{lem:uniform_discrete}(i), for every fixed $\bbD\in\mathscr D$, the map $b\mapsto\pi(b,\bbD)$ is continuous on $B_0$. By Lemma~\ref{lem:uniform_discrete}(ii), for every fixed $b\in B_0$, the map $\bbD\mapsto\pi(b,\bbD)$ is Lipschitz continuous and bounded uniformly over $b\in B_0$. Hence $\pi(\cdot,\cdot)$ is jointly continuous on $B_0\times\mathscr D$, and therefore uniformly continuous there. Thus, for every $\varepsilon>0$, there exists $\delta>0$ such that
\begin{align*}
\|b-b'\|_2\le\delta
\quad\Longrightarrow\quad
\sup_{\bbD\in\mathscr D}
|\pi(b,\bbD)-\pi(b',\bbD)|
\le
\varepsilon .
\end{align*}
Let $\{b^1,\ldots,b^N\}\subset B_0$ be a finite $\delta$-net of $B_0$. For every $b\in B_0$, choose $j(b)$ such that $\|b-b^{j(b)}\|_2\le\delta$. Then, for every $Z$,
\begin{align*}
&
\left|
\frac1Z\sum_{z=1}^Z
\pi\bigl(b,\bbD^{(z)}\bigr)
-
\mathbb E[\pi(b,\bbD)]
\right| \\
&\qquad\le
\left|
\frac1Z\sum_{z=1}^Z
\pi\bigl(b^{j(b)},\bbD^{(z)}\bigr)
-
\mathbb E[\pi(b^{j(b)},\bbD)]
\right|
+
2\varepsilon .
\end{align*}
Taking the supremum over $b\in B_0$ yields
\begin{align*}
&
\sup_{b\in B_0}
\left|
\frac1Z\sum_{z=1}^Z
\pi\bigl(b,\bbD^{(z)}\bigr)
-
\mathbb E[\pi(b,\bbD)]
\right| \\
&\qquad\le
\max_{j=1,\ldots,N}
\left|
\frac1Z\sum_{z=1}^Z
\pi\bigl(b^j,\bbD^{(z)}\bigr)
-
\mathbb E[\pi(b^j,\bbD)]
\right|
+
2\varepsilon .
\end{align*}
For each fixed $j$, the random variables $\pi(b^j,\bbD^{(z)})$ are i.i.d.\ and bounded by $\bar\pi$, so the strong law of large numbers gives
\begin{align*}
\frac1Z\sum_{z=1}^Z
\pi\bigl(b^j,\bbD^{(z)}\bigr)
\to
\mathbb E[\pi(b^j,\bbD)]
\qquad
\text{almost surely}.
\end{align*}
Since the maximum is over finitely many $j$, and since $\varepsilon>0$ is arbitrary, it follows that
\begin{align*}
\Delta_Z
:=
\sup_{b\in B_0}
|L_Z(b)-L(b)|
\to 0
\qquad
\text{almost surely}.
\end{align*}

Let
\begin{align*}
B^\star
:=
\arg\min_{b\in B_0} L(b).
\end{align*}
The set $B^\star$ is nonempty and compact because $L$ is continuous on the compact set $B_0$. Fix any $b^\star\in B^\star$. By optimality of $\widehat b^{(Z)}$ for $L_Z$ over $B_0$,
\begin{align*}
L\bigl(\widehat b^{(Z)}\bigr)
&\le
L_Z\bigl(\widehat b^{(Z)}\bigr)+\Delta_Z \\
&\le
L_Z(b^\star)+\Delta_Z \\
&\le
L(b^\star)+2\Delta_Z \\
&=
\min_{b\in B_0}L(b)+2\Delta_Z .
\end{align*}
Thus,
\begin{align*}
0
\le
L\bigl(\widehat b^{(Z)}\bigr)
-
\min_{b\in B_0}L(b)
\le
2\Delta_Z
\to0
\qquad
\text{almost surely}.
\end{align*}
Fix $\epsilon>0$ and define
\begin{align*}
S_\epsilon
:=
\Bigl\{
b\in B_0:
\operatorname{dist}(b,B^\star)\ge\epsilon
\Bigr\}.
\end{align*}
If $S_\epsilon=\emptyset$, then $\operatorname{dist}(\widehat b^{(Z)},B^\star)<\epsilon$ for all $Z$. Otherwise, $S_\epsilon$ is compact and disjoint from $B^\star$, so continuity of $L$ gives
\begin{align*}
\eta_\epsilon
:=
\min_{b\in S_\epsilon}
\left(
L(b)-\min_{u\in B_0}L(u)
\right)
>0 .
\end{align*}
For all sufficiently large $Z$ on the almost-sure event $\Delta_Z\to0$, we have
\begin{align*}
L\bigl(\widehat b^{(Z)}\bigr)
-
\min_{u\in B_0}L(u)
<
\eta_\epsilon ,
\end{align*}
which implies $\widehat b^{(Z)}\notin S_\epsilon$. Hence
\begin{align*}
\operatorname{dist}\bigl(\widehat b^{(Z)},B^\star\bigr)
\to0
\qquad
\text{almost surely}.
\end{align*}

It remains to pass from consistency of $\widehat b^{(Z)}$ to consistency of the corrected arrival-rate output. The corrected arrival rate returned by Algorithm~\ref{alg:1} need not belong to the physical demand support $\mathscr D$. We therefore use the boundedness induced by the construction of Algorithm~\ref{alg:1}.

Let
\begin{align*}
 c_{\min}T
:=
\min_{h=1,\ldots,m} c_hT
>0,
\qquad
a_{\min}
:=
\min\{A_{hj}:A_{hj}>0\}
>0 .
\end{align*}
For every empirical optimizer $\widehat b^{(Z)}\in B_0$,
\begin{align*}
0
\le
\widehat b_h^{(Z)}
\le
\bar b
:=
\frac{\bar\pi}{ c_{\min}T},
\qquad
h=1,\ldots,m .
\end{align*}
In Algorithm~\ref{alg:1}, every nonzero allocation component used to construct $\widehat\bblambda^{(Z)}$ is of the form
\begin{align*}
x_{t,j}^{(Z)}
=
\frac{\widehat b_{h(j)}^{(Z)}}{A_{h(j),j}}
\end{align*}
for some selected job $j$. Therefore,
\begin{align*}
0
\le
x_{t,j}^{(Z)}
\le
\frac{\bar b}{a_{\min}} .
\end{align*}
Since there are finitely many job columns, say $k$, and since $\widehat\lambda_{t,i}^{(Z)}$ is obtained from $R x_t^{(Z)}$, we have
\begin{align*}
0
\le
\widehat\lambda_{t,i}^{(Z)}
\le
\bar\lambda
:=
\frac{k\bar b}{a_{\min}},
\qquad
i=1,\ldots,n,\quad t=1,\ldots,T .
\end{align*}
Thus all corrected arrival-rate profiles returned by Algorithm~\ref{alg:1} belong to the compact set
\begin{align*}
\mathscr K_\lambda
:=
[0,\bar\lambda]^{nT}
\subseteq
(\mathbb R^n_+)^T .
\end{align*}

Define
\begin{align*}
\bar\pi_\lambda
:=
T\bar\lambda\sum_{i=1}^n p_i
\end{align*}
and
\begin{align*}
B_1
:=
\Bigl\{
b\in\mathbb R^m_+:
c^\top T~ b\le \max\{\bar\pi,\bar\pi_\lambda\}
\Bigr\}.
\end{align*}
Then $B_1$ is compact and contains $B_0$. Moreover, for every $\bblambda\in\mathscr K_\lambda$, every minimizer of the deterministic fluid objective lies in $B_1$, because
\begin{align*}
c^\top T~b+\pi(b,\bblambda)
\ge
 c^\top T~ b
\end{align*}
and
\begin{align*}
 c^\top T~0+\pi(0,\bblambda)
=
\pi(0,\bblambda)
\le
\bar\pi_\lambda .
\end{align*}

For $\bblambda\in\mathscr K_\lambda$, define
\begin{align*}
F(b,\bblambda)
:=
c^\top T~ b+\pi(b,\bblambda),
\qquad
b\in B_1,
\end{align*}
and
\begin{align*}
\mathcal S(\bblambda)
:=
\arg\min_{b\in B_1}F(b,\bblambda).
\end{align*}
Because minimizing over $B_1$ is equivalent to minimizing over $b\ge0$ for every $\bblambda\in\mathscr K_\lambda$, $\mathcal S(\bblambda)$ is the deterministic fluid optimal-capacity set for $\bblambda$. By Lemma~\ref{lem:uniform_discrete}, applied with $B=B_1$ and the bounded set $\mathscr K_\lambda$ in place of $\mathscr D$, the function $F$ is continuous on $B_1\times\mathscr K_\lambda$.

By the definition of the population decision-corrected arrival-rate set, for $\bblambda\in\mathscr K_\lambda$,
\begin{align*}
\bblambda\in\Lambda^\star
\quad\Longleftrightarrow\quad
\mathcal S(\bblambda)\cap B^\star\neq\emptyset .
\end{align*}

Fix an outcome in the probability-one event on which $
\operatorname{dist}\bigl(\widehat b^{(Z)},B^\star\bigr)\to0 $.
Take any subsequence $\{Z_\ell\}$. Since $\widehat\bblambda^{(Z_\ell)}\in\mathscr K_\lambda$ and $\mathscr K_\lambda$ is compact, there is a further subsequence, not relabeled, and some $\bar\bblambda\in\mathscr K_\lambda$ such that
\begin{align*}
\widehat\bblambda^{(Z_\ell)}
\to
\bar\bblambda .
\end{align*}
Since $\widehat b^{(Z_\ell)}\in B_0$ and $B_0$ is compact, we may again pass to a further subsequence, not relabeled, such that
\begin{align*}
\widehat b^{(Z_\ell)}
\to
\bar b
\end{align*}
for some $\bar b\in B_0$. The convergence
$\operatorname{dist}(\widehat b^{(Z)},B^\star)\to0$ and the closedness of $B^\star$ imply
\begin{align*}
\bar b\in B^\star .
\end{align*}

Because Algorithm~\ref{alg:1} returns $\widehat\bblambda^{(Z)}$ as a decision-corrected arrival rate for $\widehat b^{(Z)}$, we have
\begin{align*}
\widehat b^{(Z)}
\in
\mathcal S\bigl(\widehat\bblambda^{(Z)}\bigr).
\end{align*}
Equivalently, for every $b'\in B_1$,
\begin{align*}
F\bigl(\widehat b^{(Z_\ell)},\widehat\bblambda^{(Z_\ell)}\bigr)
\le
F\bigl(b',\widehat\bblambda^{(Z_\ell)}\bigr).
\end{align*}
Taking $\ell\to\infty$ and using continuity of $F$ on $B_1\times\mathscr K_\lambda$, we obtain
\begin{align*}
F(\bar b,\bar\bblambda)
\le
F(b',\bar\bblambda),
\qquad
\forall b'\in B_1 .
\end{align*}
Hence
\begin{align*}
\bar b\in\mathcal S(\bar\bblambda).
\end{align*}
Together with $\bar b\in B^\star$, this gives
\begin{align*}
\bar\bblambda\in\Lambda^\star .
\end{align*}
Thus every accumulation point of $\{\widehat\bblambda^{(Z)}\}$ belongs to $\Lambda^\star$ almost surely.

Finally, suppose for contradiction that
\begin{align*}
\operatorname{dist}\bigl(\widehat\bblambda^{(Z)},\Lambda^\star\bigr)
\nrightarrow
0
\end{align*}
on this probability-one event. Then there exist $\epsilon>0$ and a subsequence $\{Z_\ell\}$ such that
\begin{align*}
\operatorname{dist}\bigl(\widehat\bblambda^{(Z_\ell)},\Lambda^\star\bigr)
\ge
\epsilon,
\qquad
\forall \ell .
\end{align*}
By compactness of $\mathscr K_\lambda$, this subsequence has a further subsequence, not relabeled, such that
\begin{align*}
\widehat\bblambda^{(Z_\ell)}
\to
\bar\bblambda
\in
\mathscr K_\lambda .
\end{align*}
The accumulation-point argument above gives $\bar\bblambda\in\Lambda^\star$. Therefore,
\begin{align*}
\operatorname{dist}\bigl(\widehat\bblambda^{(Z_\ell)},\Lambda^\star\bigr)
\le
\|\widehat\bblambda^{(Z_\ell)}-\bar\bblambda\|_2
\to
0,
\end{align*}
which contradicts the lower bound by $\epsilon$. Hence
\begin{align*}
\operatorname{dist}\bigl(\widehat\bblambda^{(Z)},\Lambda^\star\bigr)
=
\inf_{\tilde\bblambda\in\Lambda^\star}
\|\widehat\bblambda^{(Z)}-\tilde\bblambda\|_2
\to
0
\qquad
\text{almost surely}.
\end{align*}
This proves the theorem.
\hfill\Halmos
\end{proof}

\section{Sensitivity Analysis on System Parameters}

In this appendix, we evaluate the robustness of the decision-corrected method by changing key system parameters from their baseline values. To isolate the impact of specific operational factors, we vary one parameter set at a time while maintaining all others at the baseline levels established in Section \ref{sec:numerical}. All experiments utilize the service network structure defined in Figure \ref{fig:dec_network}. Similar to Section \ref{sec:numerical}, we present the results obtained over 10 independent trials.

For the nonstationary SAA method, we first tested how many synthetic scenarios are needed by varying the generated sample size from 1 to 500. Figure~\ref{fig:sample_size} shows that the test-set costs do not change significantly across these sizes. This observation is consistent with recent findings by \cite{Besbes2023-ez}, who show that in data-driven newsvendor settings, a very small number of samples (often just one) is sufficient to achieve near-optimal performance. Although we use 500 synthetic scenarios in the main experiments of Section \ref{sec:numerical}, we use 20 scenarios for the sensitivity analysis in this appendix to reduce computational time.

\begin{figure}[H]
    \centering
    \includegraphics[width=0.9\linewidth]{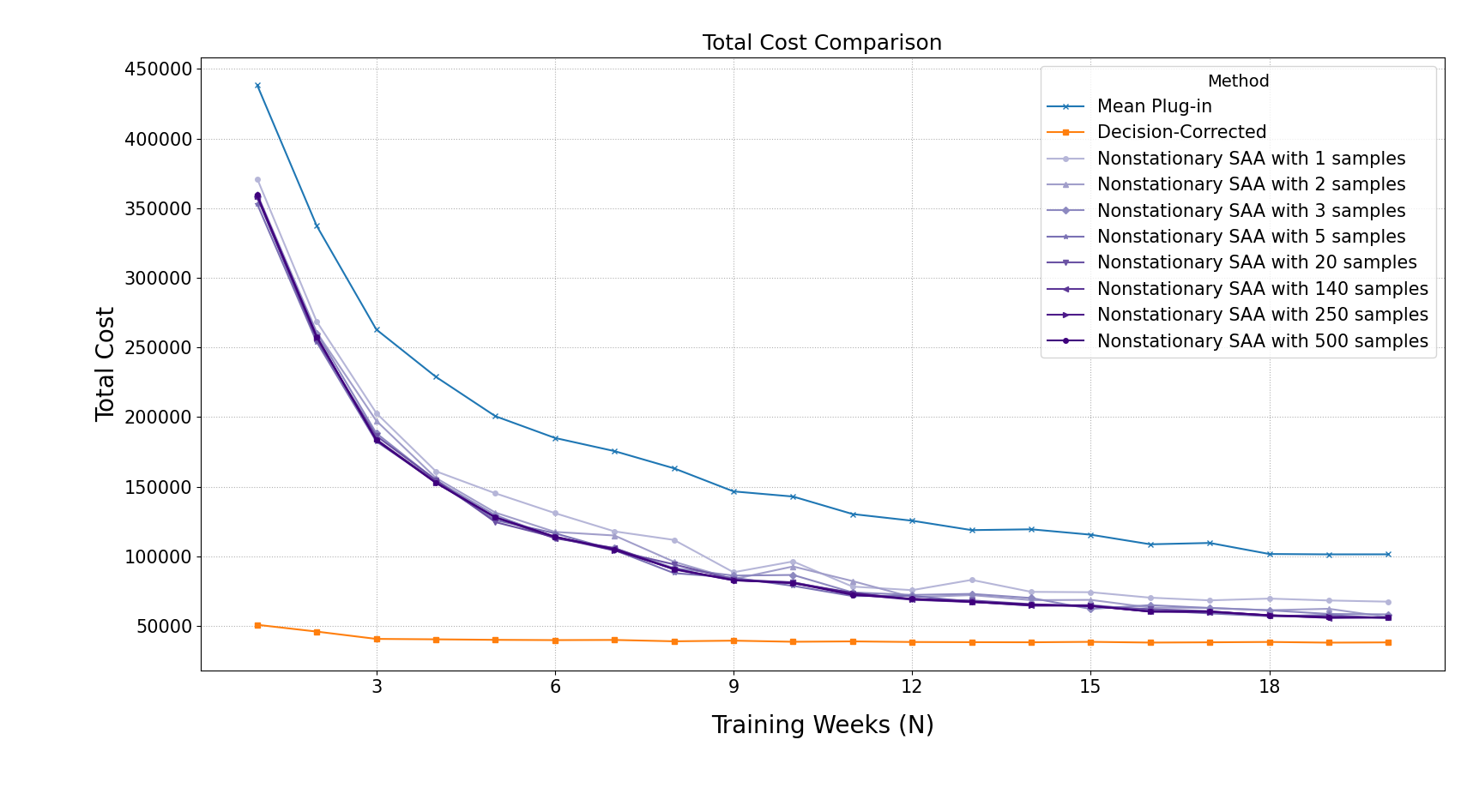}
    \caption{Total cost comparison for the nonstationary SAA benchmark under varying synthetic sample sizes (1-500).}
    \label{fig:sample_size}
\end{figure}

Table~\ref{tab:param_variations} summarizes the notation and the specific values tested for each regime.

\begin{table}[h]
\centering
\caption{Summary of Parameter Variations}
\label{tab:param_variations}
\begin{tabular}{lll}
\hline
Category & Label & Tested Values / Vector Profiles \\ \hline
Procurement Costs & $C(k)$ & $c~T \leftarrow k \cdot c_{base}~T$ for $k \in \{1, 2, 4, 8, 16\}$ \\
Lost Demand Penalties & $P(1/k)$ & $p \leftarrow p_{base} / k$ for $k \in \{1, 2, 4, 8, 16\}$ \\
Trend Factor & $F(f)$ & $f \in \{1.00, 1.05, 1.10, 1.15, 1.20\}$ \\
Demand Scale & $S(s)$ & $s \in \{3, 5, 10, 15, 20\}$ \\
Resource Usage Profiles & $A(a)$ & 6th row of $A$ varied across five profiles \\ \hline
\end{tabular}
\end{table}

\subsection{Cost Variations ($C$ and $P$)}
We examine the method's performance across different economic regimes by scaling the baseline procurement cost vector $c_{base}~T = (50, 60, 40, 70, 60, 80, 70)^\top$ and the unit lost demand penalty vector $p_{base} = (140, 135, 130, 120, 150, 175)^\top$.

\begin{itemize}
    \item \textbf{Procurement Cost Variation $C(k)$:} We define $C(k) = k \cdot c_{base}~T$. As $k$ increases, the system reflects higher fixed costs, requiring higher arrival volumes across multiple periods to justify adding capacity. The tested vectors are:
    \begin{itemize}
        \item $C(1) = (50, 60, 40, 70, 60, 80, 70)^\top$ (Baseline)
        \item $C(2) = (100, 120, 80, 140, 120, 160, 140)^\top$
        \item $C(4) = (200, 240, 160, 280, 240, 320, 280)^\top$
        \item $C(8) = (400, 480, 320, 560, 480, 640, 560)^\top$
        \item $C(16) = (800, 960, 640, 1120, 960, 1280, 1120)^\top$
    \end{itemize}
    \item \textbf{Penalty Variation $P(1/k)$:} We define $P(1/k) = p_{base} / k$. This reduces the underage cost for lost demand, testing the efficacy of the decision-corrected rate when the penalty for under-procurement is less severe.
\end{itemize}

\subsection{Demand Dynamics ($F$ and $S$)}
We vary the parameters for the stochastic demand generation process, specifically the day-over-day growth and the overall volume of patient arrivals. The underlying logic for the Poisson rates and trend factors is detailed in Appendix \ref{append:numerical_experiments}.
\begin{itemize}
    \item \textbf{Trend Factor $F(f)$:} We test $f \in \{1.00, 1.05, 1.10, 1.15, 1.20\}$. This measures the impact of different daily multiplicative growth trajectories within the planning week.
    \item \textbf{Demand Scale $S(s)$:} We vary the scaling factor $s \in \{3, 5, 10, 15, 20\}$ used to multiply the base Poisson arrival rates. This shifts the system from a lower-volume to a higher-volume environment.
\end{itemize}

\subsection{Resource Usage Profiles for the Flexible Resource ($A$)}
This regime focuses on the resource usage  parameters of the flexible resource pool ($m=6$), which is the resource capable of serving the two distinct patient classes ($n=5$ and $n=6$) in the non-decomposable part of the network. In the baseline configuration, this resource uses exactly one unit of capacity for either patient type, represented by the vector $(A_{6,7}, A_{6,8}) = (1.0, 1.0)$ in the capacity requirement matrix $A$. 

By changing these non-zero entries, we vary the relative capacity requirement of using this flexible resource. These variations allow us to observe how the decision-corrected method adapts the suggested arrival rates when a resource is more or less efficient at specific tasks:
\begin{itemize}
    \item $A(1) = (1.0, 1.0)$ (Baseline)
    \item $A(2) = (0.7, 1.0)$
    \item $A(3) = (1.0, 0.7)$
    \item $A(4) = (0.7, 1.3)$
    \item $A(5) = (1.3, 1.5)$
\end{itemize}

Figure \ref{fig:boxplot_c} shows that the decision-corrected method achieves cost savings as large as 62.5\% versus the mean plug-in benchmark, 45.8\% versus the stationary SAA benchmark, and 32.1\% versus the nonstationary SAA benchmark in regimes with lower procurement costs ($C(1)$). As the unit cost of capacity increases, all considered benchmark methods gradually reduce this performance gap. The stationary SAA benchmark outperforms the mean plug-in benchmark at lower cost levels but, notably, performs worse than the mean plug-in benchmark at higher cost levels ($C(4)$, $C(8)$, and $C(16)$). This occurs because the SAA pools all training days and cannot adapt its capacity decision to the trending demand of the target Monday, a mismatch that becomes costly when procurement is expensive. The nonstationary SAA benchmark, by contrast, holds up well across all cost levels and does not significantly fall below the mean plug-in benchmark in any regime, since it accounts for the weekly trend through its SARIMA forecast. At the most extreme cost level ($C(16)$), the high price of resources forces the mean plug-in benchmark, nonstationary SAA, and the decision-corrected method towards similar minimal capacity levels, but the decision-corrected method still outperforms the stationary SAA benchmark by 15.3\%. Importantly, in no regime does the decision-corrected method perform significantly worse than the benchmarks. This trend suggests that the decision-corrected method is most valuable in environments where resource costs allow for flexible capacity planning, rather than in prohibitive cost regimes where minimizing procurement is the dominant strategy.

\begin{figure}[htbp]
    \centering
    \includegraphics[width=0.85\textwidth]{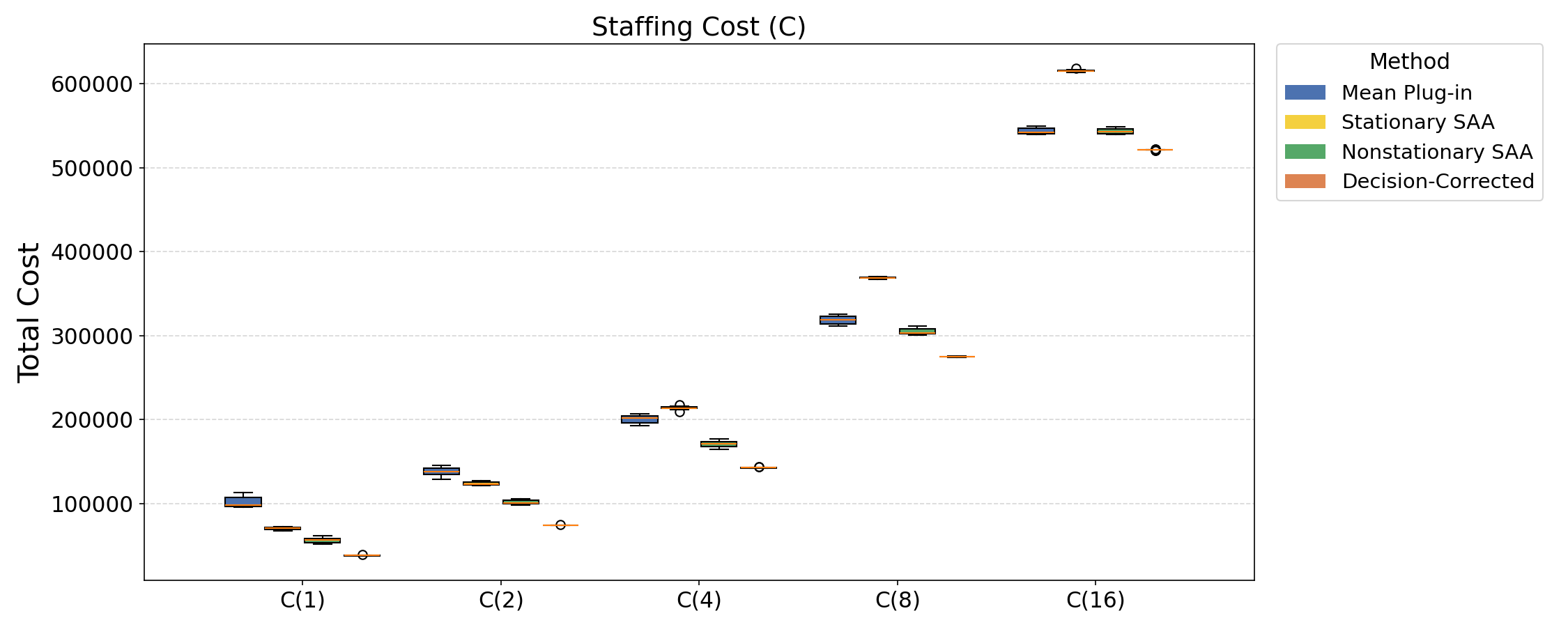}
    \caption{Total Cost Distribution under Procurement Cost Variations ($C$)}
    \label{fig:boxplot_c}
\end{figure}

Figure \ref{fig:boxplot_p} shows that the decision-corrected method delivers its strongest performance in regimes with high abandonment penalties, achieving cost savings of 62.5\% versus the mean plug-in benchmark, 45.8\% versus the stationary SAA benchmark, and 32.1\% versus the nonstationary SAA benchmark at the baseline level ($P(1)$). As the penalty for lost demand decreases, the economic pressure to match capacity lowers, and all considered benchmark methods gradually reduce this performance gap. Similar to the procurement cost variations, the stationary SAA benchmark performs worse than the mean plug-in benchmark in the lowest-penalty regimes ($P(1/4)$, $P(1/8)$, and $P(1/16)$), again because its pooling leads to misallocated capacity that is not offset by a sufficiently large penalty for unmet demand by the Mean Plug-in. The nonstationary SAA benchmark, by contrast, does not significantly fall below the mean plug-in benchmark in any regime. In the most lenient regime ($P(1/16)$), the negligible cost of lost demand drives the mean plug-in benchmark, nonstationary SAA, and the decision-corrected method toward similar minimal capacity levels, though the decision-corrected method still achieves a 15.3\% improvement over the stationary SAA benchmark. Crucially, the decision-corrected method does not perform significantly worse than the benchmarks in any regime. This suggests that while the method is most critical for service environments where unfulfilled demand carries heavy consequences, it remains a robust choice even when the cost of abandonment is minimal.

\begin{figure}[htbp]
    \centering
    \includegraphics[width=0.85\textwidth]{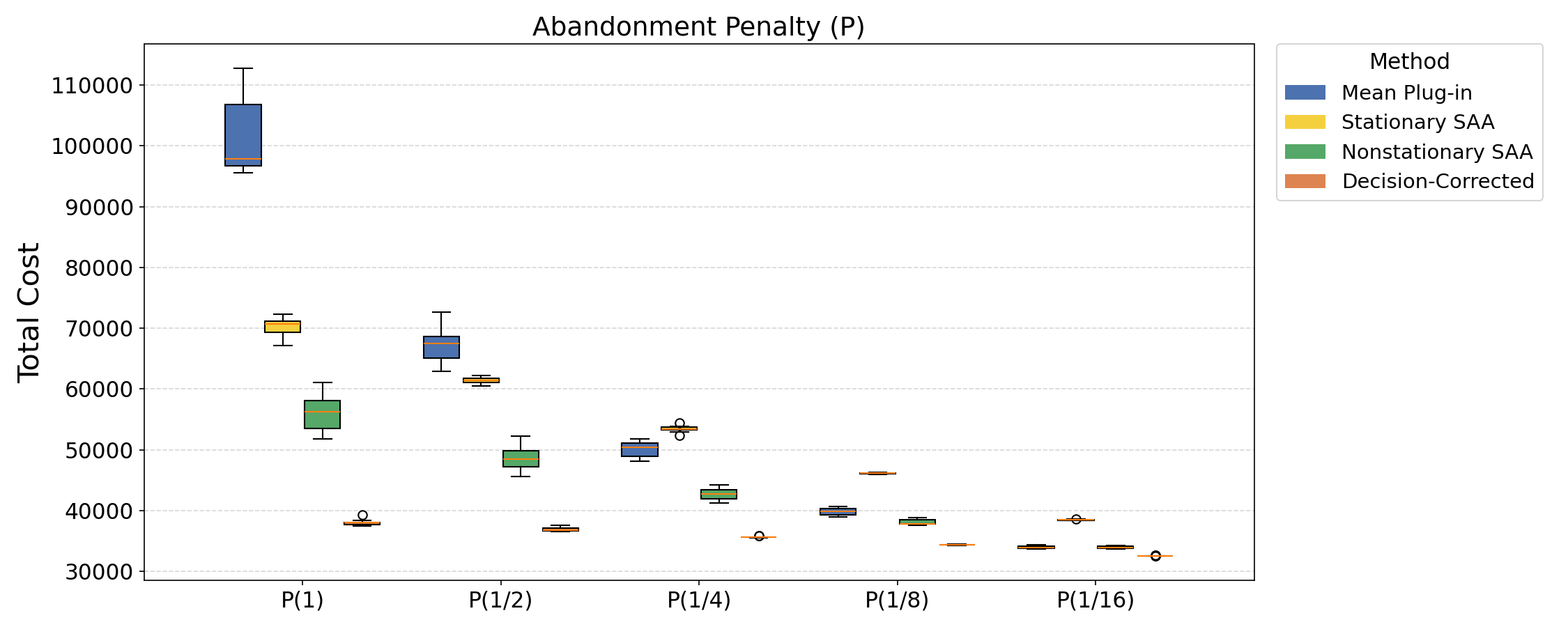}
    \caption{Total Cost Distribution under Lost Demand Penalty Variations ($P$)}
    \label{fig:boxplot_p}
\end{figure}

Figure \ref{fig:boxplot_trend} shows that the decision-corrected method keeps a performance advantage across all tested levels of daily demand growth, with cost savings above 50\% versus the mean plug-in benchmark and ranging from 13.6\% to 48.7\% versus the nonstationary SAA benchmark, with saving increasing as the trend factor increases.

The trend factor variation is particularly informative for understanding the stationary SAA benchmark. Under stationarity ($F(1.00)$), the stationary SAA benchmark and the decision-corrected method achieve nearly identical costs, confirming that when demand does not trend, pooling all historical data and solving the stationary SAA directly is a near-optimal strategy. However, as the trend factor increases, the stationary SAA benchmark degrades sharply. At $F(1.20)$, it starts performing worse than the mean plug-in benchmark. The nonstationary SAA benchmark, by contrast, handles trend better than the stationary SAA benchmark and stays the closest to the decision-corrected method across all non-zero trend factors. The decision-corrected method holds up well across all trend levels against the benchmark methods. As the trend steepens (e.g. at $F(1.20)$), the benchmark methods show a wider spread in total costs, while the decision-corrected method shows relatively lower variance across trials. This implies that in environments with strong daily growth patterns and high variability in arrivals, the decision-focused correction may offer improved consistency compared to the benchmarks.

\begin{figure}[htbp]
    \centering
    \includegraphics[width=0.85\textwidth]{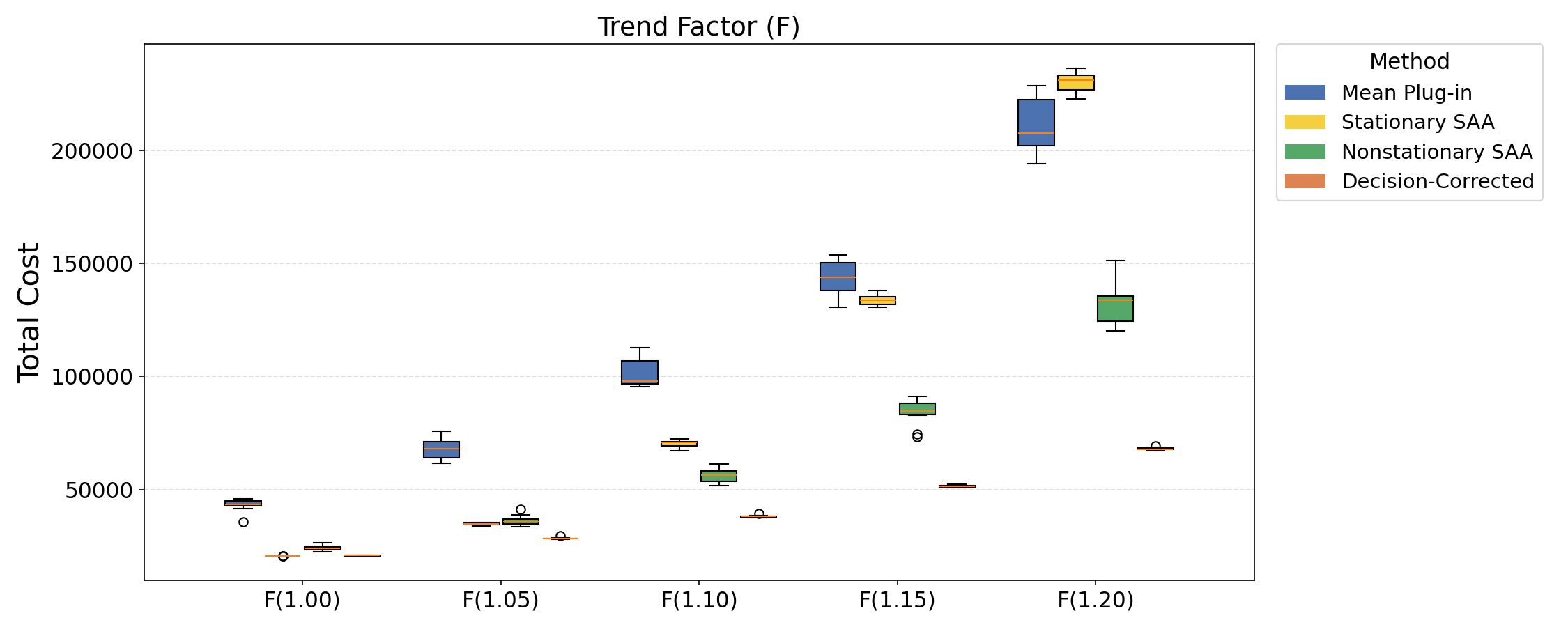}
    \caption{Total Cost Distribution under Trend Factor Variations ($F$)}
    \label{fig:boxplot_trend}
\end{figure}

Figure \ref{fig:boxplot_scale} shows that the decision-corrected method achieves cost savings of nearly 79\% versus the mean plug-in benchmark, 35.5\% versus the stationary SAA benchmark, and 46.9\% versus the nonstationary SAA benchmark in regimes with lower demand scales ($S(3)$), where the discreteness of arrivals has a noticeable impact on system performance. As the magnitude of customer arrivals increases, the mean plug-in benchmark gradually closes this gap, consistent with the theoretical expectation that fluid models become asymptotically optimal in heavy traffic. The stationary SAA benchmark shows a more stable gap, with the decision-corrected method keeping a 35.5-49.8\% improvement across all demand scales. The nonstationary SAA benchmark closes its gap more steadily as demand grows, narrowing from 46.9\% at $S(3)$ to 23.6\% at $S(20)$. Even at the highest demand level ($S(20)$), the decision-corrected method keeps a substantial advantage over all benchmark methods. This indicates that while the cost savings versus the mean plug-in and the nonstationary SAA  benchmarks are largest in lower-volume environments, the decision-focused correction remains a reliable strategy even as the system approaches high-volume limits.

\begin{figure}[htbp]
    \centering
    \includegraphics[width=0.85\textwidth]{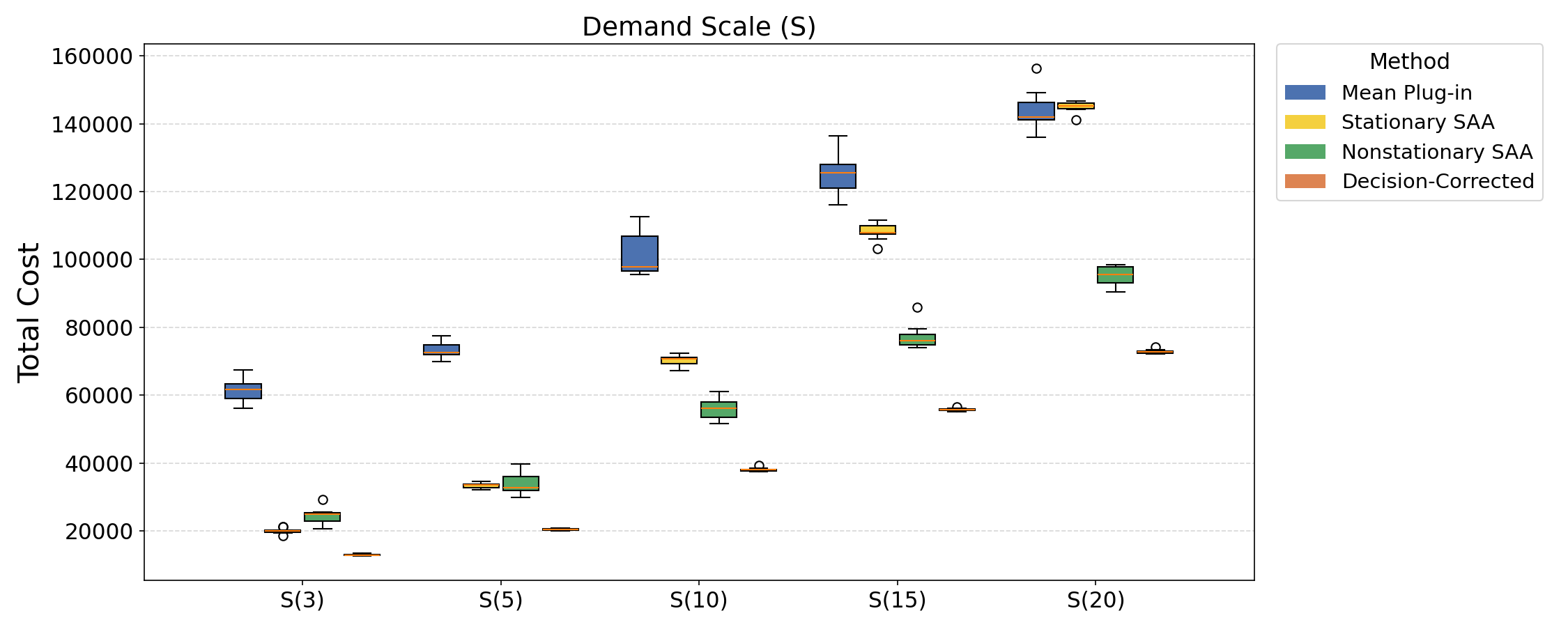}
    \caption{Total Cost Distribution under Demand Scale Variations ($S$)}
    \label{fig:boxplot_scale}
\end{figure}

Figure \ref{fig:boxplot_A} shows that the decision-corrected method keeps a performance advantage across the tested resource usage profiles. Regardless of the specific service rate combinations for the flexible resource, the decision-corrected approach consistently achieves substantial cost savings relative to all three benchmarks, ranging from 47.7\%-62.5\% versus the mean plug-in benchmark, 23.8\%-45.8\% versus the stationary SAA benchmark, and 2.7\%-32.1\% versus the nonstationary SAA benchmark. The nonstationary SAA benchmark comes closest to the decision-corrected method in the $A(3)$ regime, where its total cost is within 2.7\% of the decision-corrected method. The decision-corrected method also shows notably lower variability across trials compared to all three benchmarks. In no regime does the decision-corrected method perform significantly worse than the benchmarks. This shows that the decision-focused correction offers an advantage in reliability, keeping cost efficiency and lower variance across the tested resource usage scenarios.

\begin{figure}[htbp]
    \centering
    \includegraphics[width=0.85\textwidth]{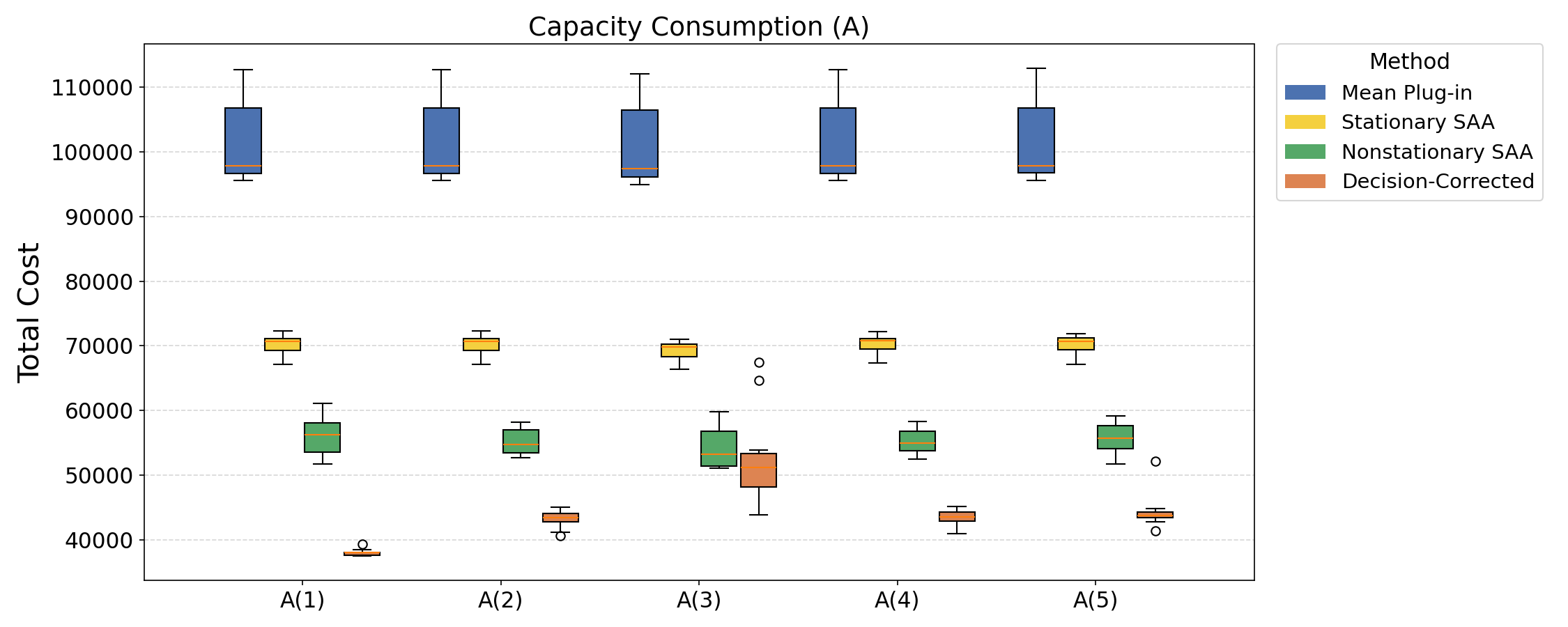}
    \caption{Total Cost Distribution under Resource Usage Profile Variations ($A$)}
    \label{fig:boxplot_A}
\end{figure}

Table~\ref{tab:cost_breakdown} provides a detailed breakdown of procurement and lost demand costs at $N = 20$ training weeks for all parameter regimes. The last columns report the percentage improvement of the decision-corrected method over each benchmark. Several patterns emerge from the table. First, the decision-corrected method achieves positive improvements over the stationary SAA benchmark in all but one setting ($F(1.00)$, where the two methods are essentially tied), and over the nonstationary SAA benchmark in every setting, with the gap narrowing most at $A(3)$, $C(16)$, and $P(1/16)$. Second, the stationary SAA benchmark sometimes incurs higher total costs than the mean plug-in benchmark (e.g., at $C(4)$, $C(8)$, $C(16)$, $P(1/4)$, $P(1/8)$, $P(1/16)$, $F(1.20)$), driven by elevated lost demand costs that result from its trend-blind capacity decisions. The nonstationary SAA benchmark does not show this behavior and stays at or below the mean plug-in benchmark across all regimes, since it accounts for the weekly trend through its SARIMA forecast.

\begin{table}[htbp]
\centering
\caption{Breakdown of Costs: Procurement vs. Lost Demand ($N=20$)}
\label{tab:cost_breakdown}
\footnotesize
\resizebox{\textwidth}{!}{
\begin{tabular}{l rrr rrr rrr rrr rrr}
\toprule
 & \multicolumn{3}{c}{\textbf{Mean Plug-in}} & \multicolumn{3}{c}{\textbf{Stationary SAA}} & \multicolumn{3}{c}{\textbf{Nonstationary SAA}} & \multicolumn{3}{c}{\textbf{Decision-Corrected}} & \multicolumn{3}{c}{\textbf{Improv. (\%)}} \\
\cmidrule(lr){2-4} \cmidrule(lr){5-7} \cmidrule(lr){8-10} \cmidrule(lr){11-13} \cmidrule(lr){14-16}
\textbf{Setting} & \textbf{Proc.} & \textbf{\shortstack{Lost\\Demand}} & \textbf{Total} & \textbf{Proc.} & \textbf{\shortstack{Lost\\Demand}} & \textbf{Total} & \textbf{Proc.} & \textbf{\shortstack{Lost\\Demand}} & \textbf{Total} & \textbf{Proc.} & \textbf{\shortstack{Lost\\Demand}} & \textbf{Total} & \textbf{vs.\ Mean} & \textbf{vs.\ Stat.} & \textbf{vs.\ Nonstat.} \\
\midrule
\multicolumn{16}{l}{\textit{Staffing Cost ($C$)}} \\
$C(1)$ & 28,522 & 72,848 & 101,370 & 30,240 & 39,992 & 70,232 & 31,435 & 24,607 & 56,041 & 35,808 & 2,244 & 38,052 & \textbf{62.46\%} & \textbf{45.82\%} & \textbf{32.10\%} \\
$C(2)$ & 56,265 & 81,695 & 137,959 & 57,337 & 66,488 & 123,824 & 60,311 & 41,320 & 101,630 & 69,846 & 4,231 & 74,077 & \textbf{46.31\%} & \textbf{40.18\%} & \textbf{27.11\%} \\
$C(4)$ & 110,421 & 89,743 & 200,164 & 106,551 & 107,211 & 213,762 & 117,057 & 53,752 & 170,809 & 134,144 & 8,504 & 142,648 & \textbf{28.73\%} & \textbf{33.27\%} & \textbf{16.49\%} \\
$C(8)$ & 215,833 & 102,761 & 318,593 & 190,363 & 178,463 & 368,827 & 221,748 & 83,204 & 304,952 & 254,486 & 20,283 & 274,769 & \textbf{13.76\%} & \textbf{25.50\%} & \textbf{9.90\%} \\
$C(16)$ & 412,398 & 131,140 & 543,538 & 314,733 & 300,728 & 615,461 & 407,360 & 135,910 & 543,270 & 465,986 & 55,309 & 521,295 & \textbf{4.09\%} & \textbf{15.30\%} & \textbf{4.04\%} \\
\addlinespace
\multicolumn{16}{l}{\textit{Penalty ($P$)}} \\
$P(1)$ & 28,522 & 72,848 & 101,370 & 30,240 & 39,992 & 70,232 & 31,435 & 24,607 & 56,041 & 35,808 & 2,244 & 38,052 & \textbf{62.46\%} & \textbf{45.82\%} & \textbf{32.10\%} \\
$P(1/2)$ & 28,192 & 39,254 & 67,446 & 28,679 & 32,672 & 61,351 & 30,435 & 18,128 & 48,563 & 34,870 & 2,035 & 36,905 & \textbf{45.28\%} & \textbf{39.85\%} & \textbf{24.01\%} \\
$P(1/4)$ & 27,605 & 22,436 & 50,041 & 26,637 & 26,806 & 53,443 & 29,265 & 13,437 & 42,701 & 33,538 & 2,125 & 35,663 & \textbf{28.73\%} & \textbf{33.27\%} & \textbf{16.48\%} \\
$P(1/8)$ & 26,979 & 12,845 & 39,824 & 23,795 & 22,310 & 46,105 & 27,719 & 10,399 & 38,118 & 31,813 & 2,534 & 34,347 & \textbf{13.75\%} & \textbf{25.50\%} & \textbf{9.89\%} \\
$P(1/16)$ & 25,775 & 8,196 & 33,971 & 19,670 & 18,796 & 38,467 & 25,460 & 8,494 & 33,954 & 29,124 & 3,457 & 32,581 & \textbf{4.09\%} & \textbf{15.30\%} & \textbf{4.04\%} \\
\addlinespace
\multicolumn{16}{l}{\textit{Trend Factor ($F$)}} \\
$F(1.00)$ & 15,708 & 27,406 & 43,114 & 19,576 & 1,031 & 20,608 & 17,907 & 6,110 & 24,016 & 19,648 & 1,106 & 20,754 & \textbf{51.86\%} & -0.71\% & \textbf{13.58\%} \\
$F(1.05)$ & 21,076 & 46,715 & 67,792 & 23,788 & 11,039 & 34,827 & 23,687 & 12,599 & 36,286 & 26,463 & 1,906 & 28,368 & \textbf{58.15\%} & \textbf{18.55\%} & \textbf{21.82\%} \\
$F(1.10)$ & 28,522 & 72,848 & 101,370 & 30,240 & 39,992 & 70,232 & 31,435 & 24,607 & 56,041 & 35,808 & 2,244 & 38,052 & \textbf{62.46\%} & \textbf{45.82\%} & \textbf{32.10\%} \\
$F(1.15)$ & 38,616 & 105,145 & 143,761 & 38,886 & 95,085 & 133,971 & 41,879 & 42,150 & 84,028 & 48,268 & 3,245 & 51,513 & \textbf{64.17\%} & \textbf{61.55\%} & \textbf{38.70\%} \\
$F(1.20)$ & 50,960 & 160,134 & 211,094 & 49,727 & 180,588 & 230,316 & 54,602 & 77,995 & 132,597 & 63,197 & 4,793 & 67,991 & \textbf{67.79\%} & \textbf{70.48\%} & \textbf{48.72\%} \\
\addlinespace
\multicolumn{16}{l}{\textit{Demand Scale ($S$)}} \\
$S(3)$ & 7,422 & 54,251 & 61,674 & 9,739 & 10,330 & 20,069 & 9,311 & 15,065 & 24,376 & 11,759 & 1,187 & 12,946 & \textbf{79.01\%} & \textbf{35.49\%} & \textbf{46.89\%} \\
$S(5)$ & 13,455 & 59,886 & 73,341 & 15,655 & 17,766 & 33,421 & 15,639 & 18,200 & 33,839 & 18,539 & 1,914 & 20,452 & \textbf{72.11\%} & \textbf{38.80\%} & \textbf{39.56\%} \\
$S(10)$ & 28,522 & 72,848 & 101,370 & 30,240 & 39,992 & 70,232 & 31,435 & 24,607 & 56,041 & 35,808 & 2,244 & 38,052 & \textbf{62.46\%} & \textbf{45.82\%} & \textbf{32.10\%} \\
$S(15)$ & 43,930 & 81,643 & 125,572 & 44,756 & 63,408 & 108,164 & 47,191 & 29,888 & 77,079 & 52,526 & 3,271 & 55,797 & \textbf{55.57\%} & \textbf{48.41\%} & \textbf{27.61\%} \\
$S(20)$ & 59,707 & 84,315 & 144,023 & 59,310 & 85,718 & 145,028 & 63,382 & 31,911 & 95,294 & 69,624 & 3,223 & 72,847 & \textbf{49.42\%} & \textbf{49.77\%} & \textbf{23.55\%} \\
\addlinespace
\multicolumn{16}{l}{\textit{Consumption Profile ($A$)}} \\
$A(1)$ & 28,522 & 72,848 & 101,370 & 30,240 & 39,992 & 70,232 & 31,435 & 24,607 & 56,041 & 35,808 & 2,244 & 38,052 & \textbf{62.46\%} & \textbf{45.82\%} & \textbf{32.10\%} \\
$A(2)$ & 28,522 & 72,848 & 101,370 & 30,240 & 39,992 & 70,232 & 31,494 & 23,697 & 55,192 & 33,804 & 9,384 & 43,189 & \textbf{57.40\%} & \textbf{38.51\%} & \textbf{21.75\%} \\
$A(3)$ & 28,236 & 72,689 & 100,925 & 29,801 & 39,518 & 69,319 & 31,070 & 23,194 & 54,264 & 32,078 & 20,729 & 52,807 & \textbf{47.68\%} & \textbf{23.82\%} & \textbf{2.68\%} \\
$A(4)$ & 28,528 & 72,848 & 101,377 & 30,336 & 39,971 & 70,307 & 31,600 & 23,701 & 55,301 & 33,941 & 9,493 & 43,434 & \textbf{57.16\%} & \textbf{38.22\%} & \textbf{21.46\%} \\
$A(5)$ & 28,530 & 72,878 & 101,408 & 30,338 & 39,918 & 70,256 & 31,601 & 24,111 & 55,712 & 33,780 & 10,650 & 44,430 & \textbf{56.19\%} & \textbf{36.76\%} & \textbf{20.25\%} \\
\bottomrule
\end{tabular}}
\end{table}

\end{APPENDICES}
 
\end{document}